\documentclass[aop,MSNbibl,seceqn,dvips]{arximspdf}
\usepackage{mathrsfs,stmaryrd}
\usepackage{graphicx}

\doi{10.1214/08-AOP419}
\volume{37}
\issue{2}
\pubyear{2009}
\firstpage{742}
\lastpage{789}

\makeatletter
\newtheorem{theorem}{Theorem}[section]
\newtheorem{lemma}[theorem]{Lemma}
\newtheorem{proposition}[theorem]{Proposition}
\newtheorem{corollary}[theorem]{Corollary}
\newproclaim{remark}{Remark}

\def\ref@fmt#1{#1}
\def\ee{e}

\def\underV{{\underline{V}}}
\setattribute{keyword}{AMS}{AMS 2000 subject classification.}
\makeatother

\begin{document}
\begin{frontmatter}

\title{Minimal position and critical martingale convergence in branching
random walks, and~directed polymers on disordered trees}
\runtitle{Minimal position, branching random walks}

\begin{aug}
\author[A]{\fnms{Yueyun} \snm{Hu}\ead[label=e1]{yueyun@math.univ-paris13.fr}\corref{}} \and
\author[B]{\fnms{Zhan} \snm{Shi}\ead[label=e2]{zhan@proba.jussieu.fr}}
\runauthor{Y. Hu and Z. Shi}
\affiliation{Universit\'e Paris XIII and Universit\'e Paris VI}
\address[A]{D\'epartement de Math\'ematiques\\
Universit\'e Paris XIII\\
99 avenue J-B Cl\'ement\\
F-93430 Villetaneuse\\
France\\
\printead{e1}} 
\address[B]{Laboratoire de Probabilit\'es\\
\quad et Mod\`eles Al\'eatoires\\
Universit\'e Paris VI\\
4 place Jussieu\\
F-75252 Paris Cedex 05\\
France\\
\printead{e2}}
\end{aug}

\received{\smonth{4} \syear{2007}}
\revised{\smonth{2} \syear{2008}}

%
\begin{abstract}
We establish a second-order almost sure limit theorem for the minimal position
in a one-dimensional super-critical branching random walk, and also
prove a
martingale convergence theorem which answers a question of
Biggins and Kyprianou [\textit{Electron. J. Probab.}
\textbf{10} (2005) 609--631]. Our method applies,
furthermore, to the study of directed polymers on a disordered tree.
In particular, we give a rigorous proof of a phase transition phenomenon
for the partition function (from the point of view of convergence in
probability), already described by Derrida and Spohn
[\textit{J. Statist. Phys.} \textbf{51} (1988) 817--840].
Surprisingly, this phase transition phenomenon disappears
in the sense of upper almost sure limits.
\end{abstract}

%
\pdfkeywords{60J80, Branching random walk,
minimal position,
martingale convergence,
spine,
marked tree,
directed polymer on a tree}

\end{frontmatter}

\section{Introduction}\label{s:intro}

\subsection{Branching random walk and martingale convergence}\label{subs:martingale-cvg}

We consider a branching random walk on the real line $\mathbb{R}$. Initially,
a particle sits at the origin. Its children form the first generation;
their displacements from the origin correspond to a point process on
the line. These children have children of their own (who form the
second generation), and behave---relative to their respective
positions---like independent copies of the initial particle. And so on.

We write $|u|=n$ if an individual $u$ is in the $n$th generation, and
denote its position by $V(u)$. [In particular, for the initial ancestor
$e$, we have $V(e)=0$.] We assume throughout the paper that, for some
$\delta>0$, $\delta_+>0$ and $\delta_->0$,
%
\begin{eqnarray}
\mathbf{E} \Biggl\{  \Biggl( \sum_{|u|=1} 1 \Biggr) ^{ 1+ \delta}   \Biggr\} &<& \infty,
\label{hyp2}\\
\mathbf{E} \Biggl\{ \sum_{|u|=1} \ee^{-(1+\delta_+) V(u)}  \Biggr\}
+ \mathbf{E} \Biggl\{ \sum_{|u|=1} \ee^{\delta_- V(u)}  \Biggr\}&<& \infty,
\label{hyp3}
\end{eqnarray}
\noindent here $\mathbf{E}$ denotes expectation with respect to
$\mathbf{P}$, the law of the branching random walk.

Let us define the (logarithmic) moment generating function
\[
\psi(t) := \log\mathbf{E} \Biggl\{ \sum_{|u|=1} \ee^{-t V(u)}
\Biggr\} \in (-\infty, \infty], \qquad t\ge0.
\]
By (\ref{hyp3}), $\psi(t) <\infty$ for $t\in[-\delta_-, 1+\delta_+]$.
Following Biggins and Kyprianou \cite{biggins-kyprianou05}, we assume
%
\begin{equation}\label{hyp}
\psi(0) >0, \qquad\psi(1) = \psi'(1) =0.
\end{equation}

\noindent Since the number of particles in each generation forms a
Galton--Watson tree, the assumption $\psi(0) >0$ in (\ref{hyp}) says
that this Galton--Watson tree is super-critical.

In the study of the branching random walk, there is a fundamental
martingale, defined as follows:
%
\begin{equation}\label{W}
W_n := \sum_{|u|=n} \ee^{-V(u)}, \qquad n =0, 1,2, \ldots\
\Biggl( \sum_\varnothing:= 0\Biggr).
\end{equation}
Since $W_n \ge0$, it converges almost surely.

When $\psi'(1)<0$, it is proved by Biggins and Kyprianou
\cite{biggins-kyprianou97} that there exists a sequence of constants $(a_n)$
such that $\frac{W_n}{a_n}$ converges in probability to a nondegenerate
limit which is (strictly) positive upon the survival of the system.
This is called the Seneta--Heyde norming in \cite{biggins-kyprianou97}
for branching random walk, referring to Seneta \cite{seneta} and
Heyde \cite{heyde} on the rate of convergence in the classic
Kesten--Stigum theorem for Galton--Watson processes.

The case $\psi'(1)=0$ is more delicate. In this case, it is known
(Lyons \cite{lyons}) that $W_n \to0$ almost surely. The following
question is raised in Biggins and Kyprianou
\cite{biggins-kyprianou05}: are there deterministic normalizers $(a_n)$ such
that $\frac{W_n}{a_n}$ converges?

We aim at answering this question.
\begin{theorem}\label{t:main}
Assume (\ref{hyp2}), (\ref{hyp3}) and (\ref{hyp}).
There exists a deterministic positive sequence
$(\lambda_n)$ with $0<\liminf_{n\to\infty}{\lambda_n\over n^{1/2}} \le
\limsup_{n\to\infty}{\lambda_n\over n^{1/2}} <\infty$, such
that, conditionally on the system's survival,
$\lambda_n W_n$ converges in distribution
to $\mathscr{W}  $, with $\mathscr{W}>0$ a.s. The
distribution of $\mathscr{W} $ is given in (\ref{wcv}).
\end{theorem}

The limit $\mathscr{W} $ in Theorem \ref{t:main} turns out to
satisfy a functional equation. Such functional equations are known to
be closely related to (a discrete version of) the
Kolmogorov--Petrovski--Piscounov (KPP) traveling wave equation; see
Kyprianou \cite{kyprianou} for more details.

The almost sure behavior of $W_n$ is described in Theorem \ref{t:Mn}
below. The two theorems together give a clear image of the asymptotics
of $W_n$.

\subsection{The minimal position in the branching random walk}\label{subs:leftmost}

A natural question in the study of branching random walks is about
$\inf_{|u|=n} V(u)$, the position of the leftmost individual in the
$n$th generation. In the literature the concentration (in terms of
tightness or even weak convergence) of $\inf_{|u|=n} V(u)$ around its
median/quantiles has been studied by many authors. See, for example,
Bachmann \cite{bachmann} and Bramson and Zeitouni
\cite{bramson-zeitouni}, as well as Section \ref{s:derrida-spohn-moment}
of the survey paper by Aldous
and Bandyopadhyay \cite{aldous-b}. We also mention the recent paper of
Lifshits \cite{lifshits}, where an example of a branching random walk is
constructed such that $\inf_{|u|=n} V(u) - \operatorname{median}
(\{\inf_{|u|=n} V(u)\})$ is tight but does not converge weakly.

We are interested in the asymptotic speed of $\inf_{|u|=n} V(u)$.
Under assumption~(\ref{hyp}), it is known that, conditionally on the
system's survival,
%
\begin{eqnarray}
\frac{1}{n} \inf_{|u|=n} V(u) &\to&0 \qquad\mbox{a.s.},
\label{hammersley-kingman-biggins}\\
\inf_{|u|=n} V(u) &\to&+\infty \qquad\mbox{a.s.}
\label{V->infty}
\end{eqnarray}
The ``law of large numbers'' in
(\ref{hammersley-kingman-biggins}) is a classic result, and can be found in
Hammersley \cite{hammersley}, Kingman \cite{kingman} and Biggins
\cite{biggins}. The system's transience to the right,
stated in (\ref{V->infty}), follows from the fact that $W_n \to0$, a.s.

A refinement of (\ref{hammersley-kingman-biggins}) is obtained by
McDiarmid \cite{mcdiarmid}. Under the additional assumption
$\mathbf{E}\{ (\sum_{|u|=1} 1 ) ^2   \} <\infty$, it is proved in \cite{mcdiarmid}
that, for some constant $c_1<\infty$ and conditionally on the system's survival,
\[
\limsup_{n\to\infty} {1\over\log n}   \inf_{|u|=n} V(u) \le c_1
\qquad\mbox{a.s.}
\]
We intend to determine the exact rate at which
$\inf_{|u|=n}V(u)$ goes to infinity.
\begin{theorem}\label{t:leftmost}
Assume (\ref{hyp2}), (\ref{hyp3}) and (\ref{hyp}).
Conditionally on the system's survival, we have
%
\begin{eqnarray}
\limsup_{n\to\infty}   {1\over\log n} \inf_{|u|=n} V(u)
&=& {3\over2} \qquad\mbox{a.s.},
\label{V-limsup-as}\\
\liminf_{n\to\infty}   {1\over\log n} \inf_{|u|=n} V(u)
&=& {1\over2} \qquad\mbox{a.s.},
\label{V-liminf-as}\\
\lim_{n\to\infty}   {1\over\log n} \inf_{|u|=n} V(u)
&=&{3\over2} \qquad\hbox{in probability.}
\label{V-proba}
\end{eqnarray}
\end{theorem}
\begin{remark*}
(i) The most interesting part of Theorem
\ref{t:leftmost} is (\ref{V-limsup-as})--(\ref{V-liminf-as}). It reveals
the presence of fluctuations of $\inf_{|u|=n} V(u)$ on the logarithmic
level, which is in contrast with known results of Bramson
\cite{bramson78bis} and Dekking and Host \cite{dekking-host} stating
that,
for a class of branching random walks, ${1\over\log\log n}
\inf_{|u|=n} V(u)$ converges almost surely to a finite and positive constant.
{\smallskipamount=0pt
\begin{longlist}[(iii)]
\item[(ii)] Some brief comments on (\ref{hyp}) are in order. In general
[i.e., without assuming $\psi(1) = \psi'(1) =0$], the law of large
numbers (\ref{hammersley-kingman-biggins}) reads\break
${1\over n} \inf_{|u|=n} V(u) \to c$, a.s.
(conditionally on the system's survival),
where $c:= \inf\{ a\in\mathbb{R}\dvtx g(a) \ge0\}$, with
$g(a) := \inf_{t\ge0} \{ ta + \psi(t) \}$. If
%
\begin{equation}\label{t*}
t^* \psi'(t^*) = \psi(t^*)
\end{equation}
for some $t^*\in(0,   \infty)$, then the branching random
walk associated with the point process
$\widehat{V}(u) := t^* V(u) +\psi(t^*)|u|$ satisfies
(\ref{hyp}). That is, as long as (\ref{t*})
has a solution [which is the case, e.g., if $\psi(1)=0$ and
$\psi'(1)>0$], the study will boil down to the case (\ref{hyp}).

It is, however, possible that (\ref{t*}) has no solution. In such a
situation, Theorem~\ref{t:leftmost} does not apply. For example, we
have already mentioned a class of branching random walks exhibited in
Bramson \cite{bramson78bis} and Dekking and Host \cite{dekking-host},
for which $\inf_{|u|=n} V(u)$ has an exotic $\log\log n$ behavior.
\item[(iii)] Under suitable assumptions, Addario--Berry \cite{addario-berry}
obtains a very precise asymptotic estimate of $\mathbf{E}[\inf_{|u|=n}
V(u)]$, which implies (\ref{V-proba}).
\item[(iv)] In the case of branching Brownian motion, the analogue of
(\ref{V-proba}) was proved by Bramson \cite{bramson78}, by means of some
powerful explicit analysis.
\end{longlist}}%
\end{remark*}

\subsection{Directed polymers on a disordered tree}\label{subs:derrida-spohn}

The following model is borrowed from the well-known paper of Derrida
and Spohn \cite{derrida-spohn}: Let $\mathbb{T}$ be a rooted Cayley
tree; we
study all self-avoiding walks ($={}$directed polymers) of $n$ steps
on~$\mathbb{T}$ starting from the root. To each edge of the tree is attached
a random variable ($={}$potential). We assume that these random variables
are independent and identically distributed. For each walk $\omega$,
its energy $E(\omega)$ is the sum of the potentials of the edges
visited by the walk. So the partition function is
\[
Z_n := \sum_{\omega} \ee^{-\beta E(\omega)},
\]
where the sum is over all self-avoiding walks of $n$ steps on
$\mathbb{T}$, and $\beta>0$ is the inverse temperature.

More generally, we take $\mathbb{T}$ to be a Galton--Watson tree, and observe
that the energy $E(\omega)$ corresponds to (the partial sum of) the
branching random walk described in the previous sections. The
associated partition function becomes
%
\begin{equation}\label{Wnbeta}
W_{n,\beta} := \sum_{|u|=n} \ee^{-\beta V(u)}, \qquad\beta>0.
\end{equation}
Clearly, when $\beta=1$, $W_{n,1}$ is just the $W_n$ defined
in (\ref{W}).

If $0<\beta<1$, the study of $W_{n,\beta}$ boils down to the case
$\psi'(1)<0$, which was investigated by Biggins and Kyprianou
\cite{biggins-kyprianou97}. In particular, conditionally on the system's
survival, ${W_{n,\beta} \over\mathbf{E}\{ W_{n,\beta}\} }$
converges almost
surely to a (strictly) positive random variable.

We study the case $\beta\ge1$ in the present paper.
\begin{theorem}\label{t:Mn}
Assume (\ref{hyp2}), (\ref{hyp3}) and
(\ref{hyp}). Conditionally on the system's survival, we have
%
\begin{equation}\label{Wn-a.s.}
W_n = n^{-1/2 + o(1)} \qquad\mbox{a.s.}
\end{equation}
\end{theorem}
\begin{theorem}\label{t:derrida-spohn}
Assume (\ref{hyp2}), (\ref{hyp3}) and (\ref{hyp}), and let
$\beta>1$. Conditionally on the system's survival, we have
%
\begin{eqnarray}
\limsup_{n\to\infty}{\log W_{n, \beta}\over\log n}
&=& - {\beta\over2} \qquad\mbox{a.s.},
\label{Wnbeta-limsup-as}\\
\liminf_{n\to\infty}{\log W_{n, \beta}\over\log n}
&=& - {3\beta\over2} \qquad\mbox{a.s.},
\label{Wnbeta-liminf-as}\\
W_{n, \beta} &=& n^{-3\beta/2 + o(1)} \qquad\hbox{in probability.}
\label{Wnbeta-proba}
\end{eqnarray}
\end{theorem}

Again, the most interesting part in Theorem \ref{t:derrida-spohn} is
(\ref{Wnbeta-limsup-as}) and (\ref{Wnbeta-liminf-as}), which describes a
new fluctuation phenomenon. Also, there is no phase transition any more
for $W_{n,\beta}$ at $\beta=1$ from the point of view of upper almost
sure limits.

The remark on (\ref{hyp}), stated after Theorem \ref{t:leftmost},
applies to Theorems \ref{t:Mn} and \ref{t:derrida-spohn} as well.

An important step in the proof of Theorems \ref{t:Mn} and
\ref{t:derrida-spohn} is to estimate all small moments of $W_n$ and
$W_{n,\beta}$, respectively. This is done in the next theorems.
\begin{theorem}\label{t:tension}
Assume (\ref{hyp2}), (\ref{hyp3}) and (\ref{hyp}).
For any $\gamma\in[0,   1)$, we have
%
\begin{equation}\label{tension}
0<\liminf_{n\to\infty}
\mathbf{E}\{ (n^{1/2}W_n)^\gamma\} \le \limsup_{n\to\infty}
\mathbf{E}\{ (n^{1/2}W_n)^\gamma\} < \infty.
\end{equation}
\end{theorem}
\begin{theorem}\label{t:derrida-spohn-moment}
Assume (\ref{hyp2}), (\ref{hyp3}) and (\ref{hyp}), and let
$\beta>1$. For any $0<r<{1\over\beta}$, we have
%
\begin{equation}\label{E(Wnbeta)}
\mathbf{E} \{ W_{n,\beta}^r  \} =
n^{-3r\beta/2+o(1)}, \qquad n\to\infty.
\end{equation}
\end{theorem}

The rest of the paper is as follows. In Section \ref{s:spine} we
introduce a change-of-measures formula (Proposition
\ref{p:change-proba}) in terms of spines on marked trees. This formula will
be of frequent use throughout the paper. Section \ref{s:Mn-lower-tail}
contains a few preliminary results of the lower tail probability of the
martingale $W_n$. The proofs of the theorems are organized as follows:
\begin{itemize}
\item Section \ref{s:proof(1.7)}: upper bound in part (\ref{V-liminf-as}) of Theorem \ref{t:leftmost}.
\item Section \ref{s:derrida-spohn-moment}: Theorem \ref{t:derrida-spohn-moment}.
\item Section \ref{s:tension}: Theorem \ref{t:tension}.
\item Section \ref{s:Mn}: Theorem \ref{t:Mn}, as well as parts
(\ref{Wnbeta-liminf-as}) and (\ref{Wnbeta-proba}) of Theorem
\ref{t:derrida-spohn}.
\item Section \ref{s:leftmost}: (the rest of) Theorem \ref{t:leftmost}.
\item Section \ref{s:derrida-spohn}: part
(\ref{Wnbeta-limsup-as}) of Theorem \ref{t:derrida-spohn}.
\item Section \ref{s:proof-t:main}: Theorem \ref{t:main}.
\end{itemize}

Section \ref{s:proof(1.7)} relies on ideas borrowed from Bramson
\cite{bramson78}, and does not require the preliminaries in Sections
\ref{s:spine} and \ref{s:Mn-lower-tail}.

Sections \ref{s:derrida-spohn-moment} and \ref{s:tension} are the
technical part of the paper, where a common idea is applied in two
different situations.

Throughout the paper we write
\[
q:= \mathbf{P}\{ \mbox{the system's extinction}\} \in[0,   1).
\]
The letter $c$ with a subscript denotes finite and (strictly)
positive constants. We also use the notation $\sum_\varnothing:= 0$,
$\prod_\varnothing:= 1$, and $0^0 := 1$. Moreover, we use
$a_n \sim b_n$, $n\to\infty$, to denote $\lim_{n\to\infty} {a_n\over b_n} =1$.

\section{Marked trees and spines}\label{s:spine}

This section is devoted to a change-of-measures result (Proposition
\ref{p:change-proba}) on marked trees in terms of spines. The material
of this section has been presented in the literature in various forms;
see, for example, Chauvin, Rouault and Wakolbinger \cite{chauvin-rouault-wakolbinger},
Lyons, Pemantle and Peres \cite{lyons-pemantle-peres},
Biggins and Kyprianou \cite{biggins-kyprianou04} and
Hardy and Harris \cite{hardy-harris}.

There is a one-to-one correspondence between branching random walks and
marked trees. Let us first introduce some notation. We label
individuals in the branching random walk by their line of descent, so
if $u=i_1\cdots i_n \in\mathscr{U}:= \{ \varnothing\} \cup
\bigcup_{k=1}^\infty(\mathbb{N}^*)^k$ (where $\mathbb{N}^* := \{ 1,2,\ldots\}$), then $u$
is the $i_n$th child of the $i_{n-1}$th child of$\ldots$ of the
$i_1$th child of the initial ancestor $e$. It is sometimes convenient
to consider an element $u\in\mathscr{U}$ as a ``word'' of length
$|u|$, with $\varnothing$ corresponding to $e$. We identify an
individual $u$ with its corresponding word.

If $u$, $v\in\mathscr{U}$, we denote by $uv$ the concatenated word,
with $u\varnothing= \varnothing u= u$.

Let $\overline{\mathscr{U}}:= \{ (u, V(u))\dvtx u\in\mathscr{U},
V\dvtx \mathscr{U} \to\mathbb{R}\}$. Let $\Omega$ be Neveu's space of marked
trees, which consists of all the subsets $\omega$ of $\overline
{\mathscr{U}}$ such that the first component of $\omega$ is a tree.
[Recall that a tree $t$ is a subset of $\mathscr{U}$ satisfying:
(i) $\varnothing\in t$; (ii) if $uj \in t$ for some $j\in\mathbb{N}^*$,
then $u\in t$; (iii) if $u\in t$, then $uj\in t$ if and only if
$1\le j\le\nu_u(t)$ for some nonnegative integer $\nu_u(t)$.]

Let $\mathbb{T}\dvtx \Omega\to\Omega$ be the identity application.
According to Neveu \cite{neveu86}, there exists a probability $\mathbf{P}$ on
$\Omega$ such that the law of $\mathbb{T}$ under $\mathbf{P}$ is the law of
the branching random walk described in the \hyperref[s:intro]{Introduction}.

Let us make a more intuitive presentation. For any $\omega\in\Omega$, let
%
\begin{eqnarray}
\mathbb{T}^{\mathrm{GW}} (\omega)
&:=& \mbox{the set of individuals ever born in $\omega$},
\label{TGW}\\
\mathbb{T}(\omega) &:=&  \{ (u, V(u)),
u\in\mathbb{T}^{\mathrm{GW}}(\omega),
V \mbox{ such that } (u, V(u))\in\omega \} .
\label{T}
\end{eqnarray}
[Of course, $\mathbb{T}(\omega) = \omega$.] In words,
$\mathbb{T}^{\mathrm{GW}}$ is a Galton--Watson tree, with
the population members as
the vertices, whereas the {\it marked tree} $\mathbb{T}$ corresponds
to the branching random walk. It is more convenient to write
(\ref{T}) in an informal way:
\[
\mathbb{T}= \{ (u, V(u)),   u\in\mathbb{T}^{\mathrm{GW}} \} .
\]

For any $u\in\mathbb{T}^{\mathrm{GW}}$, the {\it shifted} Galton--Watson
subtree generated by $u$ is
%
\begin{equation}\label{shifted-GW}
\mathbb{T}^{\mathrm{GW}}_u :=  \{ x\in\mathscr{U}\dvtx
ux\in\mathbb{T}^{\mathrm{GW}}  \}.
\end{equation}
[By shifted, we mean that $\mathbb{T}^{\mathrm{GW}}_u$ is also
rooted at $e$.] For any $x\in\mathbb{T}^{\mathrm{GW}}_u$, let
%
\begin{eqnarray}
|x|_u &:=& |ux| -|u|,
\label{|x|_u}\\
V_u(x) &:=& V(ux)-V(u).
\label{V_u}
\end{eqnarray}
As such, $|x|_u$ stands for the (relative) generation of $x$
as a vertex of the Galton--Watson tree $\mathbb{T}^{\mathrm{GW}}_u$, and
$(V_u(x),   x\in\mathbb{T}^{\mathrm{GW}}_u)$ the branching
random walk
which corresponds to the \textit{shifted marked subtree}
\[
\mathbb{T}_u :=  \{ (x, V_u(x)),   x\in\mathbb{T}^{\mathrm{GW}}_u\} .
\]

Let $\mathscr{F}_n := \sigma\{ (u, V(u)),   u\in\mathbb{T}^{\mathrm{GW}}$,
$|u| \le n\}$, which is the sigma-field induced by the first
$n$ generations of the branching random walk. Let
$\mathscr{F}_\infty$ be the sigma-field induced by the whole branching random walk.

Assume now $\psi(0)>0$ and $\psi(1)=0$. Let $\mathbf{Q}$ be a
probability on $\Omega$ such that, for any $n\ge1$,
%
\begin{equation}\label{Q}
\mathbf{Q}|_{{\mathscr F}_n} := W_n \bullet \mathbf{P}|_{{\mathscr F}_n} .
\end{equation}
Fix $n\ge1$. Let $w_n^{(n)}$ be a random variable taking
values in $\{ u\in\mathbb{T}^{\mathrm{GW}},   |u|=n\}$ such
that, for any $|u|=n$,
%
\begin{equation}\label{u*}
\mathbf{Q} \bigl\{ w_n^{(n)} = u    \big| \mathcal{F}_\infty \bigr\}
= {\ee^{-V(u)} \over W_n}.
\end{equation}
We write $\llbracket e,   w_n^{(n)}\rrbracket = \{ e=: w_0^{(n)} ,
w_1^{(n)}, w_2^{(n)}, \ldots, w_n^{(n)}\}$ for the shortest path in
$\mathbb{T}^{\mathrm{GW}}$ relating the root $e$ to $w_n^{(n)}$, with
$|w_k^{(n)}|=k$ for any $1\le k\le n$.

For any individual $u \in\mathbb{T}^{\mathrm{GW}} \setminus\{e\}$,
let $\overleftarrow{u}$ be the parent of $u$ in
$\mathbb{T}^{\mathrm{GW}}$, and
\[
\Delta V(u) := V(u) - V(\overleftarrow{u}) .
\]
For $1\le k\le n$, we write
%
\begin{equation}\label{Ck}
\mathscr{I}_k^{(n)} :=  \bigl\{ u\in \mathbb{T}^{\mathrm{GW}}\dvtx
|u| =k, \overleftarrow{u} = w_{k-1}^{(n)},  u\not= w_k^{(n)} \bigr\} .
\end{equation}
In words, $\mathscr{I}_k^{(n)}$ is the set of children of
$w_{k-1}^{(n)}$ except $w_k^{(n)}$ or, equivalently, the set of the
brothers of $w_k^{(n)}$, and is possibly empty. Finally, let us
introduce the following sigma-field:
%
\begin{equation}\label{Gn}
\mathscr{G}_n := \sigma \Biggl \{\sum_{x\in\mathscr{I}_k^{(n)}}
\delta_{\Delta V(x)},   V\bigl(w_k^{(n)}\bigr),
w_k^{(n)},   \mathscr{I}_k^{(n)}, 1\le k\le n  \Biggr\},
\end{equation}
where $\delta$ denotes the Dirac measure.

\begin{figure}

\includegraphics{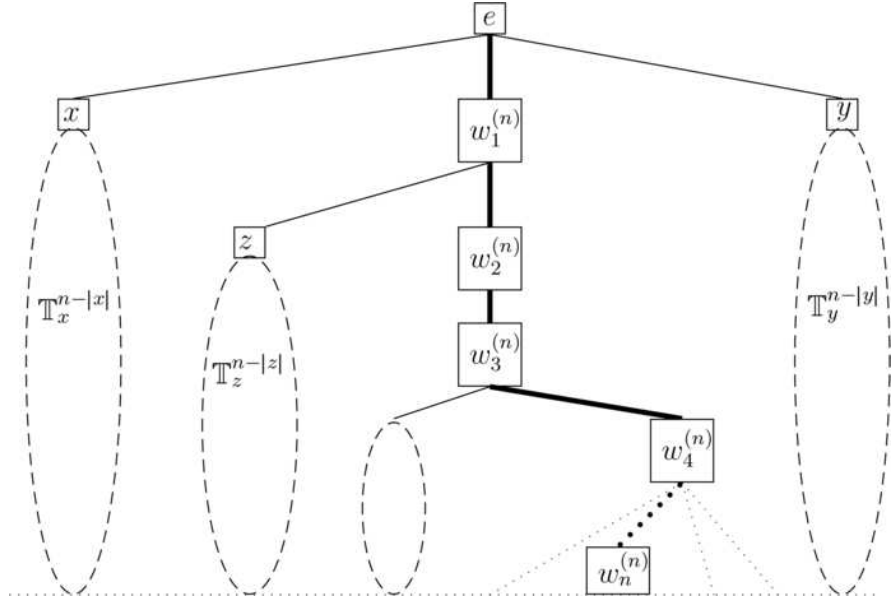}

\caption{Spine; The truncated
shifted subtrees ${\mathbb T}^{n-|x|}_x$,${\mathbb T}^{n-|y|}_y$,
${\mathbb T}^{n-|z|}_z,\ldots$ are actually rooted at $e$.\label{fig-spine}}
\vspace*{-6pt}
\end{figure}

The promised change-of-measures result is as follows. For any marked
tree $\mathbb{T}$, we define its truncation $\mathbb{T}^n$ at level
$n$ by $\mathbb{T}^n :=\{ (x, V(x)),   x\in\mathbb{T}^{\mathrm{GW}},   |x|\le n\}$;
see Figure~\ref{fig-spine}.
\begin{proposition}\label{p:change-proba}
Assume $\psi(0)>0$ and $\psi(1)=0$, and fix
$n\ge1$. Under probability $\mathbf{Q}$,
\begin{longlist}
\item the random variables
$(\sum_{x\in\mathscr{I}_k^{(n)}}\delta_{\Delta V(x)},   \Delta
V(w_k^{(n)}))$, $1\le k\le n$, are i.i.d.,
distributed as $(\sum_{x\in\mathscr{I}_1^{(1)}}
\delta_{\Delta V(x)},   \Delta V(w_1^{(1)}))$;
\item conditionally on $\mathscr{G}_n$,
the truncated shifted marked subtrees $\mathbb{T}_x^{n-|x|}$, for
\mbox{$x\in\bigcup_{k=1}^n \mathscr{I}_k^{(n)}$}, are
independent; the conditional distribution of
$ \mathbb{T}_x^{n-|x|}$ (for any
$x\in\bigcup_{k=1}^n \mathscr{I}_k^{(n)}$) under $\mathbf{Q}$,
given $\mathscr{G}_n$, is identical to the
distribution of $\mathbb{T}^{n-|x|}$ under~$\mathbf{P}$.
\end{longlist}
\end{proposition}

Throughout the paper, let $((S_i,  \sigma_i),   i\ge1)$ be such
that $(S_i-S_{i-1},   \sigma_i)$, for $i\ge1$ (with $S_0=0$), are
i.i.d. random vectors under $\mathbf{Q}$ and distributed as
$(V(w_1^{(1)}),  \# \mathscr{I}_1^{(1)})$.
\begin{corollary}\label{c:change-proba}
Assume $\psi(0)>0$ and $\psi(1)=0$, and fix $n\ge1$.
\begin{longlist}
\item Under $\mathbf{Q}$, $((V(w_k^{(n)}), \# \mathscr{I}_k^{(n)}),
1\le k\le n)$ is distributed as ($(S_k, \sigma_k)$,   \mbox{$1\le k\le n$}). In particular,
under $\mathbf{Q}$, $(V(w_k^{(n)}),   1\le k\le n)$ is
distributed as $(S_k,   1\le k\le n)$.
\item For any measurable function $F\dvtx \mathbb{R}\to\mathbb{R}_+$,
%
\begin{equation}\label{centeredRW}
\mathbf{E}_\mathbf{Q} \{ F(S_1)  \} =
\mathbf{E} \Biggl\{ \sum_{|u|=1}\ee^{-V(u)} F(V(u)) \Biggr\}.
\end{equation}
\end{longlist}
In particular, we have $\mathbf{E}_\mathbf{Q}\{ S_1\} =0$
under (\ref{hyp3}) and (\ref{hyp}).
\end{corollary}

Corollary \ref{c:change-proba} follows immediately from Proposition
\ref{p:change-proba}, and can be found in several papers (e.g.,
Biggins and Kyprianou \cite{biggins-kyprianou05}).

We present two collections of probability estimates for $(S_n)$ and for
($V(u)$, \mbox{$|u|=1$}), respectively. They are simple consequences of
Proposition \ref{p:change-proba}, and will be of frequent use in the
rest of the paper.
\begin{corollary}\label{c:S}
Assume (\ref{hyp3}) and (\ref{hyp}). Then
%
\begin{eqnarray}
\mathbf{E}_\mathbf{Q}\{ \ee^{a S_1} \}
&<&\infty \qquad\forall|a| \le c_2 ,
\label{exp-moment-S1}\\
\mathbf{Q}\{ |S_n| \ge x\} &\le& 2\exp\biggl(- c_3
\min\biggl\{ x,   {x^2\over n} \biggr\} \biggr)
\nonumber\\[-8pt]\label{Petrov}
\\[-8pt]
\eqntext{\qquad\forall n\ge1,   \forall x\ge0,}
\\
\mathbf{Q} \biggl\{ \min_{1\le k\le n} S_k >0  \biggr\}
&\sim& {c_4\over n^{1/2}}, \qquad n\to\infty,
\label{Bingham}\\
\sup_{n\ge1} n^{1/2}  \mathbf{E}_\mathbf{Q}\{
\ee^{b \min_{0\le i\le n} S_i} \} &<& \infty \qquad\forall b\ge0,
\label{E[exp(min)]}
\end{eqnarray}
where $c_2 := \min\{ \delta_+, 1+ \delta_- \}$.
Furthermore, for any $C\ge c>0$, we have
%
\begin{eqnarray}\label{tail-|Sj-Sk|}
\qquad
\mathbf{Q} \biggl\{   \max_{0\le j,   k\le n,
|j-k| \le c \log n} |S_j-S_k| \ge C\log n \biggr\}
\le2c   n^{-(c_3 C -1)}\log n
\nonumber\\[-8pt]
\\[-8pt]
\eqntext{\qquad \forall n\ge2.}
\end{eqnarray}
\end{corollary}
\begin{corollary}\label{c:V}
Assume (\ref{hyp2}), (\ref{hyp3}) and (\ref{hyp}). Let $0<a\le1$. Then
%
\begin{eqnarray}
\mathbf{E}_\mathbf{Q}
\Biggl\{  \Biggl( \sum_{|u|=1} \ee^{-aV(u)} \Biggr)^{\rho(a)}  \Biggr\}
&<& \infty,
\label{existence-moment-Wn}\\
\mathbf{Q} \biggl\{ \sup_{|u|=1} |V(u)| \ge x  \biggr\}
&\le& c_5   \ee^{-c_6   x} \qquad \forall x\ge0,
\label{tail-sup|V|}
\end{eqnarray}
with $\rho(a) := {\delta  \delta_+ \over
1+a\delta+ \delta_+}$, where $\delta$ and
$\delta_+$ are the constants in (\ref{hyp2})
and (\ref{hyp3}), respectively.
\end{corollary}
\begin{pf*}{Proof of Corollary \ref{c:S}}
By Corollary \ref{c:change-proba} (ii),
$\mathbf{E}_\mathbf{Q}\{ \ee^{a S_1}\} = \break \mathbf{E}\{ \sum_{|u|=1}
\ee^{(a-1) V(u)}\}$, which, according to (\ref{hyp3}), is
finite as long as $|a| \le c_2$. This proves (\ref{exp-moment-S1}).

Once we have the exponential integrability in (\ref{exp-moment-S1})
for $(S_n)$, standard probability estimates for sums of i.i.d. random
variables yield (\ref{Petrov}), (\ref{Bingham}) and
(\ref{E[exp(min)]}); see Petrov \cite{petrov}'s Theorem 2.7, Bingham
\cite{bingham} and Kozlov \cite{kozlov}'s Theorem A, respectively.

To check (\ref{tail-|Sj-Sk|}), we observe that the probability term on
the left-hand side of (\ref{tail-|Sj-Sk|}) is bounded by $\sum_{0\le
j< k\le n,   k-j\le c\log n} \mathbf{Q}\{ |S_{k-j}| \ge C \log n\}$.
By (\ref{Petrov}), $\mathbf{Q}\{ |S_{k-j}| \ge C \log n\} \le2 n^{-c_3 C}$ for
$k-j\le c\log n$. This yields (\ref{tail-|Sj-Sk|}).
\end{pf*}
\begin{pf*}{Proof of Corollary \ref{c:V}}
Write $\rho:= \rho (a)$. We have\break
$\mathbf{E}_\mathbf{Q}\{ ( \sum_{|u|=1} \ee^{-aV(u)} )^{\rho}\}
= \mathbf{E}_\mathbf{Q}\{ W_{1,a}^{\rho}\}
= \mathbf{E}\{ W_{1,a}^{\rho}W_{1,1}\}$. Let
$N:=\sum_{|u|=1} 1$. By H\"older's inequality,
$W_{1,a} \le W_{1,1+\delta_+}^{a/(1+\delta_+)} N^{(1-a+\delta_+)/(1+\delta_+)}$,
whereas
$W_{1,1} \le\break W_{1,1+\delta_+}^{1/(1+\delta_+)} N^{\delta_+/(1+\delta_+)}$.
Therefore, by means of another application of H\"older's
inequality, $\mathbf{E}\{ W_{1,a}^{\rho} W_{1,1}\}
\le [\mathbf{E}(W_{1,1+\delta_+})]^{(1+ a\rho)/(1+\delta_+)}
[\mathbf{E}(N^{1+\delta})]^{(\delta_+ - a\rho)/(1+\delta_+)}$,
which is finite [by (\ref{hyp3}) and~(\ref{hyp2})].
This implies (\ref{existence-moment-Wn}).

To prove (\ref{tail-sup|V|}), we write
$A:= \{ \sup_{|u|=1} |V(u)|\ge x \}$. By
Chebyshev's inequality, $\mathbf{P}(A) \le c_7 \ee^{-c_8 x}$,
where $c_7 := \mathbf{E}(\sum_{|u|=1} \ee^{c_8|V(u)|})<\infty$ as
long as $0<c_8 \le\min\{ \delta_-,   1+\delta_+\}$
[by (\ref{hyp3})]. Thus, $\mathbf{Q}(A) = \mathbf{E}\{ \sum_{|u|=1}
\ee^{-V(u)}   \mathbf{1}_A\} \le c_9  [\mathbf{P}(A)]^{\rho(1)/[1+\rho(1)]}$,
where $c_9:= [\mathbf{E}
\{ (\sum_{|u|=1} \ee^{-V(u)})^{1+\rho(1)} \}]^{1/(1+\rho(1))}<\infty$.
Now (\ref{tail-sup|V|}) follows from (\ref{existence-moment-Wn}),
with $c_6 := {c_8\rho(1)\over1+\rho(1)}$.
\end{pf*}

\section{Preliminary: small values of $W_n$}\label{s:Mn-lower-tail}

This preliminary section is devoted to the study of the small values of
$W_n$. Throughout the section, we assume (\ref{hyp2}), (\ref{hyp3})
and (\ref{hyp}). We define two important events:
%
\begin{eqnarray}
\mathscr{S} &:=&  \{  \mbox{the system's ultimate survival}   \},
\label{S}\\
\mathscr{S}_n &:=&  \{\mbox{the system's survival after $n$ generations}  \}
=  \{ W_n >0 \}.
\label{Sn}
\end{eqnarray}
Clearly, $\mathscr{S}\subset\mathscr{S}_n$. Recall (see,
e.g., Harris \cite{harris}, page~16) that, for some constant $c_{10}$ and all
$n\ge1$,
%
\begin{equation}\label{survival-survival}
\mathbf{P} \{ \mathscr{S}_n \setminus\mathscr{S} \}
\le\ee^{-c_{10} n} .
\end{equation}
Here is the main result of the section.
\begin{proposition}\label{p:tail-Mn}
Assume (\ref{hyp2}), (\ref{hyp3}) and (\ref{hyp}).
For any $\varepsilon>0$, there exists $\vartheta>0$
such that, for all sufficiently large $n$,
%
\begin{equation}\label{tail-Mn}
\mathbf{P} \{ n^{1/2}W_n < n^{-\varepsilon}
   |   \mathscr{S}  \} \le n^{-\vartheta} .
\end{equation}
\end{proposition}

The proof of Proposition \ref{p:tail-Mn} relies on Neveu's
multiplicative martingale. Recall that under assumption (\ref{hyp}),
there exists a nonnegative random variable $\xi^*$, with
$\mathbf{P}\{ \xi ^*>0\}>0$, such that
%
\begin{equation}\label{xi*}
\xi^*     \stackrel{\mathit{law}}{=}
\sum_{|u|=1} \xi_u^* \ee^{-V(u)} ,
\end{equation}
where, given $\{ (u, V(u)),   |u|=1\}$, $\xi_u^*$ are
independent copies of $\xi^*$, and ``$\stackrel{\mathit{law}}{=}$'' stands
for identity in distribution. Moreover, there is uniqueness of the
distribution of $\xi^*$ up to a scale change (see Liu \cite{liu00});
in the rest of the paper we take the version of $\xi^*$ as the unique
one satisfying $\mathbf{E}\{ \ee^{-\xi^*} \} = {1\over2}$.

Let us introduce the Laplace transform of $\xi^*$:
%
\begin{equation}\label{phi*}
\varphi^*(t) := \mathbf{E} \{ \ee^{-t \xi^*} \}, \qquad t\ge0.
\end{equation}
Let
%
\begin{equation}\label{Mn*}
W_n^* := \prod_{|u|=n} \varphi^*\bigl(\ee^{-V(u)}\bigr),
\qquad n\ge1.
\end{equation}
The process $(W_n^*,   n\ge1)$ is also a martingale
(Liu \cite{liu00}). Following Neveu \cite{neveu88}, we call $W_n^*$
an associated ``{\it multiplicative martingale}.''

The martingale $W_n^*$ being bounded, it converges almost surely (when
$n\to\infty$) to, say, $W_\infty^*$. Let us recall from Liu
\cite{liu00} (see also Kyprianou \cite{kyprianou}) that, for some $c^* >0$,
%
\begin{eqnarray}
\log{1\over W_\infty^*}
&\stackrel{\mathit{law}}{=}& \xi^*,
\label{M*infty}\\
\log \biggl( {1\over\varphi^*(t)}  \biggr)
&\sim& c^*   t \log \biggl( {1\over t} \biggr) , \qquad t \to0.
\label{phi(0)}
\end{eqnarray}
We first prove the following lemma:
\begin{lemma}\label{l:phi*}
Assume (\ref{hyp2}), (\ref{hyp3}) and (\ref{hyp}).
There exist $\kappa>0$ and $a_0\ge1$ such
that
%
\begin{eqnarray}
\mathbf{E} \{ (W_\infty^*)^a   |   W_\infty^*<1 \}
&\le& a^{-\kappa}, \qquad\forall a\ge a_0,
\label{M*infty-moment}\\
\mathbf{E} \{ (W_n^*)^a   {\bf1}_{\mathscr{S}_n}  \}
&\le& a^{-\kappa} + \ee^{-c_{10} n},
\qquad\forall n\ge1,     \forall a\ge a_0.
\label{Mn*-moment}
\end{eqnarray}
\end{lemma}
\begin{pf}
We are grateful to John
Biggins for fixing a mistake in the original proof.


We first prove (\ref{M*infty-moment}). In view of (\ref{M*infty}), it
suffices to show that
%
\begin{equation}\label{xi*-Laplace}
\mathbf{E} \{ \ee^{-a \xi^*}   |   \xi^*>0 \} \le a^{-\kappa},
\qquad a\ge a_0.
\end{equation}

Let $q \in[0,   1)$ be the system's extinction probability. Let
$N:=\sum_{|u|=1}1$. It is well known for Galton--Watson trees that $q$ is
the unique solution of $\mathbf{E}(q^N)=q$ (for $q \in[0,   1)$);
see, for example, Harris \cite{harris}, page~7. By (\ref{xi*}), $\varphi^*(t) =
\mathbf{E}\{ \prod_{|u|=1}  \varphi^* (t \ee^{-V(u)})\}$.
Therefore, by
(\ref{phi*}), $\mathbf{P}\{ \xi^* =0\} = \varphi^* (\infty)
= \lim_{t\to \infty} \mathbf{E}\{ \prod_{|u|=1} \varphi^*
(t \ee^{-V(u)})\}$, which, by
dominated convergence, is${}= \break \mathbf{E}\{ (\varphi^*(\infty))^N\}
=\mathbf{E}\{(\mathbf{P}\{ \xi^* =0\})^N\}$. Since
$\mathbf{P}\{ \xi^* =0\} <1$, this yields $\mathbf{P} \{ \xi^* =\break 0\} =q$.

Following Biggins and Grey \cite{biggins-grey}, we note that, for any
$t\ge0$,
\[
\mathbf{E} \{ \ee^{-t \xi^*}  \}
= q + (1-q) \mathbf{E} \{ \ee^{-t\xi^*}   |   \xi^*>0 \} .
\]
Let $\widehat{\xi}$ be a random variable such that
$\mathbf{E}\{\ee^{-t \widehat{\xi}} \}
= \mathbf{E}\{ \ee^{-t \xi^*}   |   \xi^*>0\}$
for any $t\ge0$. Let $Y$ be a random variable independent of
everything else, such that $\mathbf{P}\{ Y=0\} =q
= 1-\mathbf{P}\{Y=1\}$. Then $\xi ^*$ and $Y\widehat{\xi}$
have the same law and, by (\ref{xi*}), so
do $\xi^*$ and $\sum_{|u|=1} \ee^{-V(u)} Y_u \widehat{\xi}_u$,
where, given $\{ u,   |u|=1\}$, $(Y_u, \widehat{\xi}_u)$ are
independent copies of $(Y, \widehat{\xi})$, independent of
$\{ V(u),  |u|=1\}$. Since $\{\sum_{|u|=1}
\ee^{-V(u)} Y_u \widehat{\xi}_u >0\} = \{ \sum_{|u|=1} Y_u >0\}$, this leads to
\[
\mathbf{E}\{ \ee^{-t \widehat{\xi}} \}
= \mathbf{E} \Biggl\{\ee^{-t \sum_{|u|=1}
\ee^{-V(u)} Y_u \widehat{\xi}_u}    \Big|   \sum_{|u|=1}
Y_u >0  \Biggr\}, \qquad t\ge0.
\]

Let $\widehat{\varphi}(t) := \mathbf{E}\{
\ee^{-t \widehat{\xi}} \}$,
$t\ge0$. Then for any $t\ge0$ and $c>0$,
\[
\widehat{\varphi}(t) = \mathbf{E} \Biggl\{ \prod_{|u|=1}
\widehat{\varphi}\bigl(t\ee^{-V(u)} Y_u\bigr)
\bigg |   \sum_{|u|=1} Y_u >0  \Biggr\}
\le\mathbf{E}
\Biggl\{ [ \widehat{\varphi}(t\ee^{-c}) ]^{N_c}   \Big|
\sum_{|u|=1} Y_u >0  \Biggr\} ,
\]
where $N_c := \sum_{|u|=1}   {\bf1}_{\{ Y_u=1,
|V(u)|\le c\} }$. By monotone convergence,
$\lim_{c\to\infty}\mathbf{E}\{N_c  |\break  \sum_{|u|=1} Y_u >0\}
= \mathbf{E}\{ \sum_{|u|=1} Y_u|  \sum_{|u|=1} Y_u >0\} >1$
[because $\mathbf{P}\{ \sum_{|u|=1} Y_u \ge\break 2\}>0$ by
assumption (\ref{hyp})]. We can therefore choose and fix a constant
$c>0$ such that $\mathbf{E}\{ N_c  |  \sum_{|u|=1} Y_u >0\}>1$. By writing
$\widehat{f}(s) :=\mathbf{E}\{ s^{N_c}   |   \sum_{|u|=1} Y_u>0 \}$, we have
\[
\widehat{\varphi}(t) \le\widehat{f}( \widehat{\varphi}
(t \ee^{-c})), \qquad\forall t\ge0.
\]
Iterating the inequality yields that, for any $t\ge0$ and any
$n\ge1$,
%
\begin{equation}\label{psi>h(psi)}
\qquad\quad
\mathbf{E}\{ \ee^{-t \widehat{\xi}} \} \le
\widehat{f}^{  (n)}
(\mathbf{E}\{ \ee^{-t\ee^{-nc} \widehat{\xi}} \} ),
\qquad \mbox{that is, }
\mathbf{E}\{ \ee^{-t\ee^{nc} \widehat{\xi}} \} \le
\widehat{f}^{  (n)}
(\mathbf{E}\{ \ee^{-t \widehat{\xi}} \} ),
\end{equation}
where $\widehat{f}^{  (n)}$ denotes the $n$th iterate of
$\widehat{f}$. It is well known for Galton--Watson trees (Athreya and
Ney \cite{athreya-ney}, Section I.11) that, for any $s\in[0,   1)$,
$\lim_{n\to\infty} \gamma^{-n} \*\widehat{f}^{  (n)} (s)$ converges
to a finite limit, with $\gamma:= (\widehat{f}  )'(0)
\le\mathbf{P}\{ \sum_{|u|=1} Y_u =1   | \break  \sum_{|u|=1} Y_u> 0 \} <1$. Therefore,
(\ref{psi>h(psi)}) yields (\ref{xi*-Laplace}), and thus (\ref{M*infty-moment}).

It remains to check (\ref{Mn*-moment}). Let $a\ge1$. Since
$((W_n^*)^a,   n\ge0)$ is a bounded submartingale,
$\mathbf{E}\{ (W_n^*)^a  {\bf1}_{\mathscr{S}_n} \} \le\mathbf{E}
\{ (W_\infty^*)^a{\bf1}_{\mathscr{S}_n} \}$.
Recall that $W_\infty^* \le1$; thus,
\[
\mathbf{E} \{ (W_n^*)^a   {\bf1}_{\mathscr{S}_n}  \}
\le\mathbf{E} \{(W_\infty^*)^a   {\bf1}_\mathscr{S} \}
+ \mathbf{P} \{\mathscr{S}_n \setminus \mathscr{S} \} .
\]
By (\ref{survival-survival}),
$\mathbf{P}\{ \mathscr{S}_n\setminus\mathscr{S}\}
\le\ee^{-c_{10} n}$. To estimate $\mathbf{E}\{ (W_\infty^*)^a
{\bf1}_\mathscr{S}\}$, we identify $\mathscr{S}$ with
$\{ W_\infty^*<1\}$: on the one hand,
$\mathscr{S}^c \subset\{ W_n^*=1, \mbox{ for all sufficiently large }n\}
\subset\{ W_\infty^*=1\}$; on the other hand, by (\ref{M*infty}),
$\mathbf{P}\{ W_\infty^*<1\} = \mathbf{P}\{ \xi^* >0\} = 1-q =\mathbf{P}(\mathscr{S})$.
Therefore, $\mathscr{S}= \{ W_\infty^*<1\}$, $\mathbf{P}$-a.s.
Consequently, $\mathbf{E}\{(W_\infty^*)^a
{\bf1}_\mathscr{S}\} = \mathbf{E}\{ (W_\infty ^*)^a
{\bf1}_{\{W_\infty^*<1\} } \}$, which, according to (\ref{M*infty-moment}), is
bounded by $a^{-\kappa}$, for $a\ge a_0$. Lemma \ref{l:phi*} is
proved.
\end{pf}

We are now ready for the proof of Proposition \ref{p:tail-Mn}.\vadjust{\goodbreak}
\begin{pf*}{Proof of Proposition \ref{p:tail-Mn}}
Let $c_{11}>0$
be such that $\mathbf{P}\{ \xi^*\le c_{11}\} \ge{1\over2}$. Then
$\varphi^*(t)= \mathbf{E}\{ \ee^{-t \xi^*}\}
\ge\ee^{-c_{11}t}\mathbf{P}\{ \xi^*\le
c_{11}\} \ge{1\over2} \ee^{-c_{11}t}$ and, thus,
$\log({1\over\varphi^*(t)}) \le c_{11}t + \log2$.
Together with (\ref{phi(0)}),
this yields, on the event $\mathscr{S}_n$,
\begin{eqnarray*}
\log \biggl( {1\over W_n^*} \biggr)
&=& \sum_{|u|=n} \log
\biggl ( {1\over \varphi^*(\ee^{-V(u)})} \biggr)
\\
&\le& \sum_{|u|=n} {\bf1}_{\{ V(u)\ge1\} }
c_{12} V(u)\ee^{-V(u)} + \sum_{|u|=n}
{\bf1}_{\{ V(u)< 1\} }
\bigl(c_{11}\ee^{-V(u)} + \log2 \bigr) .
\end{eqnarray*}
Since $W_n = \sum_{|u|=n} \ee^{-V(u)}$, we obtain,
on $\mathscr{S}_n$, for any $\lambda\ge1$,
%
\begin{equation}\label{Mn*Mn}
\log \biggl( {1\over W_n^*} \biggr)
\le c_{13}   \lambda W_n + c_{12} \sum_{|u|=n}
{\bf1}_{\{ V(u)\ge\lambda\} }V(u)\ee^{-V(u)} ,
\end{equation}
where $c_{13} := c_{11} + c_{12} + \ee\log2$. Note that
$c_{12}$ and $c_{13}$ do not depend on $\lambda$.

Let $0<y\le1$. Since $\mathscr{S}\subset\mathscr{S}_n$, it follows
that, for $c_{14}:= c_{12} + c_{13}$,
%
\begin{eqnarray}\label{tail-Mn1}
\mathbf{P} \{ \lambda W_n < y    |   \mathscr{S}_n \}
&\le& \mathbf{P} \biggl\{ \log \biggl( {1\over W_n^*}  \biggr)
< c_{14}y    \big|   \mathscr{S}_n   \biggr\}
\nonumber \\
&&{} + \mathbf{P} \Biggl\{ \sum_{|u|=n} {\bf1}_{\{ V(u)\ge
\lambda\} } V(u) \ee^{-V(u)} \ge y    \Big| \mathscr{S}_n  \Biggr\}
\\
& =&\!: \mathrm{RHS}_{(\ref{tail-Mn1})}^1
+ \mathrm{RHS}_{(\ref{tail-Mn1})}^2,
\nonumber
\end{eqnarray}
with obvious notation.

Recall that $\mathbf{P}(\mathscr{S}_n)
\ge\mathbf{P}(\mathscr{S})=1-q$. By Chebyshev's inequality,
\[
\mathrm{RHS}_{(\ref{tail-Mn1})}^1 \le\ee^{c_{14}}
\mathbf{E} \{(W_n^*)^{1/y}   |   \mathscr{S}_n   \}
\le{\ee^{c_{14}}\over1-q}
\mathbf{E} \{ (W_n^*)^{1/y}   {\bf1}_{\mathscr{S}_n}  \} .
\]
By (\ref{Mn*-moment}), for $n\ge1$ and $0<y\le{1\over
a_0}$, with $c_{15} := \ee^{c_{14}}/(1-q)$,
%
\begin{equation}\label{RHS1<}
\mathrm{RHS}_{(\ref{tail-Mn1})}^1 \le
c_{15}  (y^\kappa+ \ee^{-c_{10} n}  ) .
\end{equation}

To estimate $\mathrm{RHS}_{(\ref{tail-Mn1})}^2$, we observe that
\begin{eqnarray*}
\mathrm{RHS}_{(\ref{tail-Mn1})}^2
&\le& {1\over1-q}  \mathbf{P} \Biggl\{ \sum_{|u|=n}
{\bf1}_{\{ V(u)\ge\lambda\} } V(u) \ee^{-V(u)}\ge y  \Biggr\}
\\
&\le& {1\over(1-q)y} \mathbf{E} \Biggl\{ \sum_{|u|=n}
{\bf1}_{\{ V(u)\ge\lambda\} } V(u) \ee^{-V(u)} \Biggr\}
\\
&=& {1\over(1-q)y} \mathbf{E}_\mathbf{Q} \Biggl\{ \sum_{|u|=n}
{\bf1}_{\{ V(u)\ge\lambda\} }{V(u) \ee^{-V(u)}\over W_n}  \Biggr\}
\\
&=& {1\over(1-q)y} \mathbf{E}_\mathbf{Q}\bigl\{ V\bigl({w_n^{(n)}}\bigr)
{\bf1}_{\{ V({w_n^{(n)}})\ge\lambda\} }\bigr\} .
\end{eqnarray*}
By Corollary \ref{c:change-proba}(i),
$\mathbf{E}_\mathbf{Q}\{V({w_n^{(n)}})   {\bf1}_{\{ V({w_n^{(n)}})\ge\lambda\} }\}
= \mathbf{E}_\mathbf{Q}
\{ S_n   {\bf1}_{\{ S_n \ge\lambda\} } \} \le\break (\mathbf{E}_\mathbf{Q}
\{ S_n^2\})^{1/2}
 (\mathbf{Q}\{ S_n \ge\lambda\})^{1/2}$, which, by
(\ref{Petrov}), is bounded by\break $c_{16}  n   \exp( - c_3 \min\{ \lambda,
  {\lambda^2\over n} \})$. Accordingly, $\mathrm{RHS}_{(\ref
{tail-Mn1})}^2 \le{c_{17}  n\over y}\exp( - c_3 \min\{ \lambda,
{\lambda^2\over n} \})$. Together with (\ref{tail-Mn1}) and
(\ref{RHS1<}), it yields that, for $0<y\le{1\over a_0}$,
\[
\mathbf{P} \{ \lambda W_n < y    |   \mathscr{S}_n\}
\le c_{15}( y^\kappa+ \ee^{-c_{10} n}  )
+ {c_{17}  n\over y}
\exp \biggl( - c_3 \min\biggl\{ \lambda,   {\lambda^2\over n} \biggr\}  \biggr) .
\]
Let $\lambda:= n^{1/2} y^{-\kappa/2}$. The inequality
becomes, for $0<y\le{1\over a_0}$ and $n\ge1$,
\begin{eqnarray*}
&& \mathbf{P} \bigl\{ n^{1/2} W_n < y^{(\kappa+2)/2}
   |   \mathscr{S}_n \bigr\}
\\
&&\qquad \le c_{15}
  (y^\kappa+\ee^{-c_{10} n}  )
  +{c_{17}   n \over y}
\exp \biggl( - c_3 {\min\{n^{1/2} y^{\kappa/2},   1\}
\over y^\kappa} \biggr).
\end{eqnarray*}
This readily yields Proposition \ref{p:tail-Mn}.
\end{pf*}
\begin{remark*}
%
Under the additional assumption that $\{ u,   |u|=1\}$ contains at
least two elements almost surely, it is possible (Liu \cite{liu01}) to
improve (\ref{M*infty-moment}): $\mathbf{E}\{ (W_\infty^*)^a   |
W_\infty^*<1 \} \le\exp\{-a^{\kappa_1}\}$ for some $\kappa_1>0$
and all sufficiently large $a$, from which one can deduce the stronger
version of Proposition \ref{p:tail-Mn}: for any $\varepsilon>0$,
there exists $\vartheta_1>0$ such that $\mathbf{P}\{ n^{1/2}W_n <
n^{-\varepsilon}   |   \mathscr{S}  \} \le\exp(- n^{\vartheta_1})$ for
all sufficiently large~$n$.
\end{remark*}

We complete this section with the following estimate which will be
useful in the proof of Theorem \ref{t:tension}.
\begin{lemma}\label{l:log(1/M)}
Assume (\ref{hyp2}), (\ref{hyp3}) and (\ref{hyp}).
For any $0<s<1$,
%
\begin{equation}\label{E(log(1/M))}
\sup_{n\ge1} \mathbf{E} \biggl\{  \biggl( \log
{1\over W_n^*} \biggr)^{  s}  \biggr\} < \infty.
\end{equation}
\end{lemma}
\begin{pf}
Let $x>1$. By Chebyshev's inequality,
$\mathbf{P}\{ \log({1\over W_n^*}) \ge x\} =\break \mathbf{P}\{
\ee^x W_n^* \le1\} \le\ee  \mathbf{E}\{
\ee^{- \ee^x W_n^*}\}$. Since
$W_n^*$ is a martingale, it follows from Jensen's inequality that
$\mathbf{E}\{ \ee^{- \ee^x W_n^*}\} \le\mathbf{E}\{
\ee^{- \ee^x W_\infty^*}\} \le
\mathbf{P}\{ W_\infty^* \le\ee^{-x/2}\} +\break \exp(-\ee^{x/2})$.
Therefore,
%
\begin{equation}\label{log(1/Mn)}
\mathbf{P} \biggl\{ \log\biggl({1\over W_n^*}\biggr) \ge x  \biggr\}
\le\ee  \mathbf{P}\{ W_\infty^* \le\ee^{-x/2}\}
+ \exp ( 1-\ee^{x/2}  ) .
\end{equation}

On the other hand, by integration by parts,
$\int_0^\infty\ee^{-ty}\mathbf{P}(\xi^* \ge y)  \,dy
= {1-\mathbf{E}(\ee^{-t \xi^*}) \over t} = {1-\varphi^*(t) \over t}$,
which, according to (\ref{phi(0)}), is $\le c_{18}
\log( {1\over t})$ for $0<t\le{1\over2}$. Therefore, for $a\ge2$,
$c_{18} \log a \ge\int_0^\infty\ee^{-y/a}
\mathbf{P}(\xi^*\ge y)   \,dy \ge\int_0^a \ee^{-y/a} \mathbf{P}(\xi^* \ge a)  \,dy
= (1-\ee^{-1}) a\mathbf{P}(\xi^* \ge a)$. That is,
$\mathbf{P}(\xi^* \ge a) \le{c_{18} \over 1-\ee^{-1}}
{\log a \over a}$ or, equivalently, $\mathbf{P}(W_\infty^* \le\ee^{-a})
\le{c_{18} \over1-\ee^{-1}}   {\log a \over a}$, for
$a\ge 2$. Substituting this in (\ref{log(1/Mn)}) gives that, for any $x\ge4$,
\[
\mathbf{P} \biggl\{ \log\biggl({1\over W_n^*}\biggr) \ge x  \biggr\}
\le{2\ee  c_{18} \over1-\ee^{-1}}  {\log(x/2)\over x}
+ \exp ( 1-\ee^{x/2} ) .
\]
Lemma \ref{l:log(1/M)} follows immediately.
\end{pf}

\section[Proof of Theorem 1.2: upper bound in (1.8)]{Proof of Theorem \protect\ref{t:leftmost}: upper bound in
(\protect\ref{V-liminf-as})}\label{s:proof(1.7)}

Assume (\ref{hyp2}), (\ref{hyp3}) and~(\ref{hyp}). This section
is devoted to proving the upper bound in (\ref{V-liminf-as}):
conditionally on the system's survival,
%
\begin{equation}\label{V-liminf-as-ub}
\liminf_{n\to\infty}
{1\over\log n} \inf_{|u|=n} V(u) \le {1\over2}, \qquad\hbox{a.s.}
\end{equation}

The proof borrows some ideas from Bramson \cite{bramson78}. We fix
$-\infty<a<b<\infty$ and $\varepsilon>0$. Let $\ell_1\le\ell_2
\le2\ell_1$ be integers; we are interested in the asymptotic case
$\ell_1\to\infty$. Consider $n \in[\ell_1,   \ell_2]\cap\mathbb{Z}$.
Let $0<c_{19}<1$ be a constant, and let
\[
g_n(k) := \min \{ c_{19}   k^{1/3},   c_{19} (n-k)^{1/3} + a
\log\ell_1,   n^\varepsilon \}, \qquad0\le k \le n.
\]
Let $\mathbb{L}_n$ be the set of individuals
$x\in\mathbb{T}^{\mathrm{GW}}$ with $|x|=n$ such that
\[
g_n(k)\le V(x_k)\le c_{20}   k, \qquad\forall0\le k \le n
\quad \mbox{and} \quad
a \log\ell_1 \le V(x) \le b\log\ell_1,
\]

\noindent where $x_0:=e, x_1,\ldots, x_n:= x$ are the vertices
on the shortest path in $\mathbb{T}^{\mathrm{GW}}$ relating the
root $e$
and the vertex $x$, with $|x_k|=k$ for any $0\le k\le n$. We consider
the measurable event
\[
F_{\ell_1,\ell_2} := \bigcup_{n=\ell_1}^{\ell_2}
\bigcup_{|x|=n}\{ x\in\mathbb{L}_n \} .
\]
We start by estimating the first moment of
$\# F_{\ell_1,\ell_2}$:
$\mathbf{E}(\# F_{\ell_1,\ell_2}) =\break  \sum_{n=\ell_1}^{\ell_2}
\mathbf{E}\{ \sum_{|x|=n} {\bf1}_{ \{ x\in\mathbb{L}_n \} } \}$.
Since $\mathbf{E}\{\sum_{|x|=n} {\bf1}_{ \{ x\in\mathbb{L}_n \} } \}
= \mathbf{E}_\mathbf{Q}\{ \sum_{|x|=n}{\ee^{-V(x)} \over W_n}
\ee^{V(x)} \times\break {\bf1}_{ \{ x\in \mathbb{L}_n\} } \} =\mathbf{E}
_\mathbf{Q}\{ \ee^{V(w_n^{(n)})} {\bf1}_{ \{ w_n^{(n)}\in\mathbb{L}_n \} } \}$,
we can apply Corollary \ref{c:change-proba} to see that
\begin{eqnarray*}
\mathbf{E}(\# F_{\ell_1,\ell_2}) &=& \sum_{n=\ell_1}^{\ell_2}
\mathbf{E}_\mathbf{Q} \bigl\{ \ee^{S_n}
{\bf1}_{ \{ g_n(k)\le S_k \le c_{20}   k,
\forall0\le k\le n, a \log\ell_1 \le S_n \le b\log\ell_1 \} }  \bigr\}
\\
&\ge& \sum_{n=\ell_1}^{\ell_2} \ell_1^a
\mathbf{Q} \{ g_n(k)\le S_k \le c_{20}   k, \forall0\le k\le n,
a \log\ell_1\le S_n \le b\log\ell_1 \}.
\end{eqnarray*}
We choose (and fix) the constants $c_{19}$ and $c_{20}$ such
that $\mathbf{Q}\{ c_{19}< S_1 < c_{20}\} >0$. Then,
\footnote{An easy way to
see why $-{3\over2}$ should be the correct exponent for the
probability is to split the event into three pieces: the first piece
involving $S_k$ for $0\le k\le{n\over3}$, the second piece for
${n\over3} \le k\le{2n\over3}$, and the third piece for ${2n\over3}
\le k\le n$. The probability of the first piece is $n^{-(1/2)+o(1)}$
(it is essentially the probability of $S_k$ being positive for $1\le
k\le{n\over3}$, because conditionally on this, $S_k$ converges
weakly, after a suitable normalization, to a Brownian meander; see
Bolthausen \cite{bolthausen}). Similarly, the probability of the third
piece is $n^{-(1/2)+o(1)}$. The second piece essentially says that
after ${n\over3}$ steps, the random walk should lie in an interval of
length of order $\log n$; this probability is also $n^{-(1/2)+o(1)}$.
Putting these pieces together yields the claimed exponent $-{3\over2}$.

For a rigorous proof, the upper bound---not required here---is easier
since we can only look at the event that the walk stays positive during
$n$ steps (with the same condition upon the random variable~$S_n$),
whereas the lower bound needs some tedious but elementary writing,
based on the Markov property. Similar arguments are used for the random
walk $(S_k)$ in several other places in the paper, without further
mention.}
the probability $\mathbf{Q}\{ \cdot\}$ on the right-hand
side is $\ell_1^{-(3/2) + o(1)}$, for $\ell_1\to\infty$. Accordingly,
%
\begin{equation}\label{E(F)}
\mathbf{E}(\# F_{\ell_1,\ell_2}) \ge
(\ell_2-\ell_1+1) \ell_1^{a-(3/2) + o(1)}.
\end{equation}

We now proceed to estimate the second moment of
$\# F_{\ell_1,\ell_2}$. By definition,
\begin{eqnarray*}
\mathbf{E}[(\# F_{\ell_1,\ell_2})^2]
&=& \sum_{n=\ell_1}^{\ell_2} \sum_{m=\ell_1}^{\ell_2}
\mathbf{E} \Biggl\{ \sum_{|x|=n} \sum_{|y|=m}
{\bf1}_{ \{ x\in\mathbb{L}_n,   y\in\mathbb{L}_m\} }  \Biggr\}
\\
&\le& 2\sum_{n=\ell_1}^{\ell_2}
\sum_{m=n}^{\ell_2} \mathbf{E} \Biggl\{ \sum_{|x|=n} \sum_{|y|=m}
{\bf1}_{ \{ x\in\mathbb{L}_n,   y\in\mathbb{L}_m \} } \Biggr\} .
\end{eqnarray*}
We look at the double sum $\sum_{|x|=n} \sum_{|y|=m}$ on
the right-hand side. By considering~$z$, the youngest common ancestor
of $x$ and $y$, and writing $k:= |z|$, we arrive at\looseness=-1
\[
\sum_{|x|=n} \sum_{|y|=m} {\bf1}_{ \{ x\in\mathbb{L}_n,  y\in \mathbb{L}_m \} }
= \sum_{k=0}^n \sum_{|z|=k} \sum_{(u,v)} {\bf1}_{ \{ zu\in\mathbb{L}_n,
  zv\in\mathbb{L}_m \} } ,
\]
where the double sum $\sum_{(u,v)}$ is over $u$,
$v\in \mathbb{T}^{\mathrm{GW}}_z$ such that $|u|_z= n-k$ and $|v|_z = m-k$ and
that the unique common ancestor of $u$ and $v$ in
$\mathbb{T}^{\mathrm{GW}}_z$ is the root. Therefore,
\begin{eqnarray*}
\mathbf{E}[(\# F_{\ell_1,\ell_2})^2] &\le&
2\sum_{n=\ell_1}^{\ell_2} \sum_{m=n}^{\ell_2} \sum_{k=0}^n
\mathbf{E}
\Biggl \{ \sum_{|z|=k} \sum_{(u,v)} {\bf1}_{ \{ zu\in\mathbb{L}_n,
  zv\in \mathbb{L}_m \} }  \Biggr\}
\\
&=&\!: 2\sum_{n=\ell_1}^{\ell_2}
 \sum_{m=n}^{\ell_2} \sum_{k=0}^n \Lambda_{k,n,m} .
\end{eqnarray*}

We estimate $\Lambda_{k,n,m}$ according to three different situations.

\textit{First situation}: $0\le k\le\lfloor n^\varepsilon\rfloor$. Let
$V_z(u) := V(zu)- V(z)$ as in Section \ref{s:spine}. We have $0\le
g_n(k) \le V(z) \le c_{20}   n^\varepsilon$, and $V(zu_i) \ge0$ for
$0\le i\le n-k$ and $V(zu_{n-k}) \le b\log\ell_1$, where
$u_0:= e,u_1,\ldots, u_{n-k}$ are the vertices on the shortest path in
$\mathbb{T}^{\mathrm{GW}}_z$ relating the root $e$ and the vertex $u$,
with $|u_i|_z=i$ for any $0\le i\le n-k$.\vadjust{\goodbreak} Therefore,
$V_z(u_i) \ge -c_{20}   n^\varepsilon$ for $0\le i\le n-k$, and
$V_z(u) \le b\log\ell_1$. Accordingly,
\[
\Lambda_{k,n,m} \le\mathbf{E} \Biggl\{ \sum_{|z|=k} \sum_{v\in
\mathbb{T}^{\mathrm{GW}}_z,
|v|_z = m-k} {\bf1}_{ \{ zv\in\mathbb{L}_m \} } B_{n-k} \Biggr\} ,
\]
where
\begin{eqnarray*}
B_{n-k} :\!\!&=& \mathbf{E} \Biggl\{ \sum_{|x|= n-k}
{\bf1}_{ \{ V(x_i) \ge-c_{20}n^\varepsilon,
  \forall0\le i\le n-k, V(x) \le b\log\ell_1\} }  \Biggr\}
\\
&=& \mathbf{E}_\mathbf{Q}
\bigl\{ \ee^{V(w_{n-k}^{(n-k)})}
{\bf1}_{ \{ V(w_i^{(n-k)}) \ge- c_{20}   n^\varepsilon,
\forall0\le i\le n-k, V(w_{n-k}^{(n-k)}) \le b\log\ell_1\} }  \bigr\}
\\
&=& \mathbf{E}_\mathbf{Q} \bigl\{ \ee^{S_{n-k}}
{\bf1}_{ \{ S_i \ge - c_{20}   n^\varepsilon,
\forall0\le i\le n-k, S_{n-k} \le b\log\ell_1 \} }  \bigr\}
\\
&\le& \ell_1^b   \mathbf{Q} \{ S_i \ge
- c_{20}   n^\varepsilon, \forall0\le i\le n-k,
S_{n-k} \le b \log\ell_1  \}
\\
&\le& \ell_1^{b-(3/2) + \varepsilon+ o(1)}
\le\ell_1^{b- (3/2) + 2\varepsilon}.
\end{eqnarray*}
Therefore,
\begin{eqnarray*}
\Lambda_{k,n,m} &\le& \ell_1^{b- (3/2) + 2\varepsilon} \mathbf{E}
\Biggl\{ \sum_{|z|=k} \sum_{v\in\mathbb{T}^{\mathrm{GW}}_z,
|v|_z = m-k}{\bf1}_{ \{ zv\in\mathbb{L}_m \} }  \Biggr\}
\\
&=& \ell_1^{b- (3/2) + 2\varepsilon}\mathbf{E}
\Biggl\{ \sum_{|x|=m} {\bf1}_{ \{ x\in\mathbb{L}_m \} } \Biggr\}
\end{eqnarray*}
and, thus,
%
\begin{equation}\label{Fj-1}
\quad
\sum_{n=\ell_1}^{\ell_2} \sum_{m=n}^{\ell_2}
\sum_{k=0}^{\lfloor n^\varepsilon\rfloor}
\Lambda_{k,n,m} \le \ell_1^{b- (3/2) + 2\varepsilon}(\ell_2-\ell_1+1)
(\ell_2^\varepsilon+1) \mathbf{E}(\# F_{\ell_1,\ell_2}).
\end{equation}

\textit{Second situation}: $\lfloor n^\varepsilon\rfloor+1 \le k\le\min
\{ m-\lfloor n^\varepsilon\rfloor,   n\}$. In this situation, since
$V(z) \ge\max\{ g_m(k),   g_n(k)\} \ge c_{19}   n^{\varepsilon/3}$,
we simply have $V_z(u) \le b\log\ell_1 - c_{19}n^{\varepsilon/3}$.
Exactly as in the first situation, we get
\[
\Lambda_{k,n,m} \le\mathbf{E} \Biggl\{ \sum_{|x|=m} {\bf1}_{ \{ x\in
\mathbb{L}_m\} }  \Biggr\}  \mathbf{E}
\Biggl\{ \sum_{|x|= n-k} {\bf1}_{ \{ V(x) \le b\log
\ell_1 - c_{19}   n^{\varepsilon/3}\} }  \Biggr\} .
\]
The second $\mathbf{E}\{ \cdot\}$ on the right-hand side is
\[
= \mathbf{E}_\mathbf{Q} \bigl\{ \ee^{S_{n-k}} {\bf1}_{ \{
S_{n-k} \le b \log\ell_1 - c_{19}   n^{\varepsilon/3}\} } \bigr\}
\le\ell_1^b   \ee^{-c_{19}   n^{\varepsilon/3}}
\]
and, thus,
%
\begin{equation}\label{Fj-2}
\qquad
\sum_{n=\ell_1}^{\ell_2} \sum_{m=n}^{\ell_2}
\sum_{k=\lfloor n^\varepsilon\rfloor+1}^{\min\{ m-\lfloor n^\varepsilon\rfloor,
n\}} \Lambda_{k,n,m} \le\ell_1^b
\ee^{- c_{19}   \ell_1^{\varepsilon/3}}
(\ell_2-\ell_1+1) \ell_2 \mathbf{E}(\# F_{\ell_1,\ell_2}).
\end{equation}

\textit{Third and last situation}:
$m-\lfloor n^\varepsilon\rfloor+1 \le k\le n$
(this situation may happen only if $m\le n+\lfloor
n^\varepsilon\rfloor-1$). This time\vadjust{\goodbreak}
$V(z) \ge g_m(k) \ge a\log\ell _1$ and, thus,
$V_z(u) \le(b-a)\log\ell_1$; consequently,
\begin{eqnarray*}
\Lambda_{k,n,m} &\le& \mathbf{E}
\Biggl\{ \sum_{|x|=m} {\bf1}_{ \{ x\in\mathbb{L}_m \} }  \Biggr\}
\mathbf{E} \Biggl\{ \sum_{|x|= n-k} {\bf1}_{ \{ V(x) \le(b-a)\log
\ell_1\} } \Biggr\}
\\
&\le& \ell_1^{b-a}
\mathbf{E} \Biggl\{ \sum_{|x|=m} {\bf1}_{ \{ x\in\mathbb{L}_m \} }  \Biggr\} .
\end{eqnarray*}
Therefore, in case $m\le n+\lfloor n^\varepsilon\rfloor-1$,
\begin{eqnarray*}
\sum_{n=\ell_1}^{\ell_2} \sum_{m=n}^{\ell_2}
\sum_{k=m-\lfloor n^\varepsilon\rfloor+1}^n
\Lambda_{k,n,m} &\le& \sum_{n=\ell_1}^{\ell_2}
\sum_{m=n}^{n+\lfloor n^\varepsilon\rfloor-1}
\ell_2^\varepsilon\ell_1^{b-a}\mathbf{E} \Biggl\{ \sum_{|x|=m}
{\bf1}_{ \{ x\in\mathbb{L}_m \} }  \Biggr\}
\\
&\le& \sum_{m=\ell_1}^{\ell_2}
\sum_{n=m-2\lfloor m^\varepsilon\rfloor}^m
\ell_2^\varepsilon\ell_1^{b-a} \mathbf{E} \Biggl\{ \sum_{|x|=m}
{\bf1}_{ \{ x\in\mathbb{L}_m \} }  \Biggr\}
\\
&\le& 2\ell_2^{2\varepsilon}
\ell_1^{b-a}   \mathbf{E}(\# F_{\ell_1,\ell_2}) .
\end{eqnarray*}
Combining this with (\ref{Fj-1}) and (\ref{Fj-2}), and
since
\[
\mathbf{E}[(\# F_{\ell_1,\ell_2})^2]
\le2\sum_{n=\ell_1}^{\ell_2}\sum_{m=n}^{\ell_2} \sum_{k=0}^n \Lambda_{k,n,m},
\]
we obtain
\begin{eqnarray*}
{\mathbf{E}[(\# F_{\ell_1,\ell_2})^2] \over[\mathbf{E}
(\# F_{\ell_1,\ell_2})]^2} &\le&
\bigl(2\ell_1^{b-(3/2)+2\varepsilon} (\ell_2-\ell_1+1)(\ell_2^\varepsilon+1)
\\
&&\hspace*{2pt}{}
+ 2\ell_1^b \ee^{-c_{19}   \ell_1^{\varepsilon/3}}
(\ell_2-\ell_1+1) \ell_2 +4\ell_2^{2\varepsilon}
\ell_1^{b-a}\bigr)
(\mathbf{E}(\# F_{\ell_1,\ell_2}))^{-1}.
\end{eqnarray*}
Since $\ell_2\le2\ell_1$, we have
$2\ell_1^{b-(3/2)+2\varepsilon}(\ell_2^\varepsilon+1) + 2\ell_1^b
\ee^{-c_{19}   \ell_1^{\varepsilon/3}} \ell_2
\le\ell_1^{b-(3/2)+4\varepsilon}$ for all sufficiently large $\ell_1$. On
the other hand, $\mathbf{E}(\# F_{\ell_1,\ell_2})
\ge(\ell_2-\ell_1+1)\ell_1^{a-(3/2) -\varepsilon}$ by (\ref{E(F)}) (for large
$\ell_1$). Therefore, when $\ell_1$ is large, we have
\[
{\mathbf{E}[(\# F_{\ell_1,\ell_2})^2] \over[\mathbf{E}
(\# F_{\ell_1,\ell_2})]^2}\le{\ell_1^{b-(3/2)+4\varepsilon}
(\ell_2-\ell_1+1) +\ell_1^{b-a+3\varepsilon}\over(\ell_2-\ell_1+1)
\ell_1^{a-(3/2)-\varepsilon}} .
\]
By the Paley--Zygmund inequality, $\mathbf{P}\{ F_{\ell_1,\ell_2}
\not= \varnothing\} \ge{1\over4} {[\mathbf{E}(\# F_{\ell_1,\ell_2})]^2
\over\mathbf{E}[(\# F_{\ell_1,\ell_2})^2]}$; thus,
%
\begin{equation}\label{McDiarmid1}
\qquad
\mathbf{P} \biggl\{ \min_{\ell_1 \le|x| \le\ell_2} V(x)
\le b\log\ell_1  \biggr\} \ge{1\over4}
{(\ell_2-\ell_1+1) \ell_1^{a-(3/2) -\varepsilon}
\over\ell_1^{b-(3/2)+4\varepsilon}
(\ell_2-\ell_1+1) +\ell_1^{b-a+3\varepsilon}} .
\end{equation}
Of course, we can make $a$ close to $b$, and $\varepsilon$
close to 0, to see that, for any $b\in\mathbb{R}$ and $\varepsilon>0$, all
sufficiently large $\ell_1$ and all $\ell_2\in[\ell_1,
2\ell_1]\cap\mathbb{Z}$,
%
\begin{equation}\label{McDiarmid2}
\mathbf{P} \biggl\{ \min_{\ell_1 \le|x| \le\ell_2} V(x)
\le b\log\ell_1  \biggr\} \ge {\ell_2-\ell_1+1\over
\ell_1^\varepsilon(\ell_2-\ell_1+1)
+\ell_1^{(3/2)-b+\varepsilon}} .
\end{equation}
[This is our basic estimate for the minimum of $V(x)$. In
Section \ref{s:derrida-spohn-moment} we are going to apply
(\ref{McDiarmid2}) to $\ell_2 :=\ell_1$.]

We now let $b>{1\over2}$ and take the subsequence $n_j := 2^j$,
$j\ge j_0$ (with a sufficiently large integer $j_0$). By (\ref{McDiarmid2})
(and possibly by changing the value of $\varepsilon$),
\[
\mathbf{P} \biggl\{ \min_{n_j \le|x| \le n_{j+1}} V(x) \le b\log n_j
 \biggr\} \ge n_j^{-\varepsilon}.
\]
Let $\tau_j := \inf\{ k\dvtx  \# \{ u\dvtx   |u|=k\} \ge
n_j^{2\varepsilon}\}$. Then we have, for $j\ge j_0$,
\begin{eqnarray*}
&& \mathbf{P} \biggl\{ \tau_j <\infty,
\min_{\tau_j + n_j \le|x| \le\tau_j + n_{j+1}}
V(x) > \max_{|y|=\tau_j} V(y) + b\log n_j  \biggr\}
\\
&&\qquad \le \biggl( \mathbf{P} \biggl\{
\min_{n_j \le|x| \le n_{j+1}} V(x)> b\log n_j \biggr\}
 \biggr)^{\lfloor n_j^{2\varepsilon}\rfloor}
\\
&&\qquad \le (1- n_j^{-\varepsilon})^{\lfloor
n_j^{2\varepsilon}\rfloor} ,
\end{eqnarray*}
which is summable in $j$. By the Borel--Cantelli lemma,
almost surely for all large $j$, we have either $\tau_j=\infty$, or
$\min_{\tau_j + n_j \le|x| \le\tau_j + n_{j+1}} V(x) \le\max
_{|y|=\tau_j} V(y) + b\log n_j$.

By the well-known law of large numbers for the branching random walk
(Hammersley \cite{hammersley}, Kingman \cite{kingman} and Biggins
\cite{biggins}), of which (\ref{hammersley-kingman-biggins}) was a
special case, there exists a constant $c_{21}>0$ such that
${1\over n} \max_{|y|=n} V(y) \to c_{21}$ almost surely upon the system's
survival. In particular, upon survival, $\max_{|y|=n} V(y)\le2c_{21}  n$
almost surely for all large $n$. Consequently, upon the system's
survival, almost surely for all large $j$, we have either
$\tau_j=\infty$, or $\min_{\tau_j + n_j \le|x| \le\tau_j + n_{j+1}}
V(x) < 2c_{21} \tau_j + b\log n_j$.

Recall that the number of particles in each generation forms a
Galton--Watson tree, which is super-critical under assumption
(\ref{hyp}) (because $m:=\break  \mathbf{E}\{ \sum_{|u|= 1}1\} >1$). In particular,
conditionally on the system's survival, ${\# \{ u:   |u|= k\} \over
m^k}$ converges a.s.\ to a (strictly) positive random variable when
$k\to\infty$, which implies $\tau_j \sim2\varepsilon{\log n_j
\over\log m}$ a.s.\ ($j\to\infty$). As a consequence, upon the
system's survival, we have, almost surely for all large $j$,
\[
\min_{n_j \le|x| \le2n_{j+1}} V(x) \le{5\varepsilon c_{21}
\over \log m}   \log n_j + b \log n_j.
\]
Since $b$ can be as close to ${1\over2}$ as possible, this
readily yields (\ref{V-liminf-as-ub}).

\section[Proof of Theorem 1.6]{Proof of Theorem \protect\ref{t:derrida-spohn-moment}}\label{s:derrida-spohn-moment}

Before proving Theorem \ref{t:derrida-spohn-moment}, we need three estimates.

The first estimate, stated as Proposition
\ref{p:mcdiarmid-integrability}, was proved
by McDiarmid \cite{mcdiarmid}
under the additional assumption $\mathbf{E}\{ (\sum_{|u|=1}1)^2\}<\infty$.
\begin{proposition}\label{p:mcdiarmid-integrability}
Assume (\ref{hyp2}), (\ref{hyp3}) and
(\ref{hyp}). There exists $c_{22}>0$ such that,
for any $\varepsilon>0$, we can find
$c_{23} = c_{23}(\varepsilon) >0$ satisfying
%
\begin{equation}\label{mcdiarmid-integrability}
\mathbf{E} \biggl\{ \exp \biggl( c_{22} \inf_{|x|=n} V(x)
 \biggr)   {\bf1}_{\mathscr{S}_n}  \biggr\}
\le c_{23}  n^{{(3+\varepsilon)/2} c_{22}}, \qquad n\ge1.
\end{equation}
\end{proposition}
\begin{remark*}
Since $W_n \ge\exp[ - \inf_{|x|=n} V(x)]$,
it follows from (\ref{mcdiarmid-integrability}) and\break  H\"older's
inequality that, for any $0\le s<c_{22}$ and $\varepsilon>0$,
%
\begin{equation}\label{Mn-moment-neg}
\mathbf{E} \biggl\{ {1\over W_n^s}   {\bf1}_{\mathscr{S}_n} \biggr\}
\le c_{23}^{s/c_{22}} n^{{(3+\varepsilon)/2} s}, \qquad n\ge1.
\end{equation}
This estimate will be useful in the proof of Theorem
\ref{t:tension} in Section \ref{s:tension}.
\end{remark*}
\begin{pf*}{Proof of Proposition \ref{p:mcdiarmid-integrability}}
In the proof we write, for any $k\ge0$,
\[
\underV_k := \inf_{|u|=k} V(u).
\]

Taking $\ell_2=\ell_1$ in (\ref{McDiarmid2}) gives that, for any
$\varepsilon>0$ and all sufficiently large $\ell$ (say,
$\ell\ge\ell_0$), we have
$\mathbf{P}\{ \underV_\ell\le{3\over2}\log\ell\}\ge\ell^{-\varepsilon}$;
thus, $\mathbf{P}\{ \underV_\ell> {3\over2}\log\ell\}
\le1- \ell^{-\varepsilon} \le\exp(- \ell^{-\varepsilon})$,
$\forall\ell\ge\ell_0$. For any $r\in\mathbb{R}$ and integers
$k\ge1$ and $n>\ell\ge\ell_0$, we have
\begin{eqnarray*}
&& \mathbf{P} \bigl\{ \underV_n >\tfrac{3}{2} \log\ell+ r \bigr\}
\\
&&\qquad \le\mathbf{P} \bigl\{ \# \{ u\dvtx  |u|=n-\ell,   V(u)\le r\} < k \bigr\}
+ \bigl ( \mathbf{P} \bigl\{ \underV_\ell > \tfrac{3}{2}\log\ell \bigr\}  \bigr)^k
\\
&&\qquad \le \mathbf{P} \bigl\{ \# \{ u\dvtx   |u|=n-\ell,   V(u)
\le r\} < k \bigr\} + \exp(- \ell^{-\varepsilon}k) .
\end{eqnarray*}
By Lemma 1 of McDiarmid \cite{mcdiarmid}, there exist
$c_{24}>0$, $c_{25}>0$ and $c_{26}>0$ such that, for any $j\ge1$,
$\mathbf{P} \{ \# \{ u\dvtx |u|=j,   V(u) \le c_{24} j\} \le\ee^{c_{25} j}
\} \le q + \ee^{-c_{26} j}$, $q$ being as before the probability of
extinction. We choose $j:= \lfloor{r\over c_{24}} \rfloor$ and
$\ell:= n- \lfloor{r\over c_{24}} \rfloor$ to see that,
for all $n\ge\ell_0$ and all $0\le r\le c_{24}(n-\ell_0)$,
\[
\mathbf{P} \bigl\{ \underV_n > \tfrac{3}{2}\log n + r  \bigr\} \le q
+ \ee^{-c_{26} \lfloor r/c_{24}\rfloor}
+ \exp \bigl(- n^{-\varepsilon}
\lfloor\ee^{c_{25} \lfloor r/c_{24}\rfloor}\rfloor \bigr).
\]
Noting that $\{ \underV_n > {3\over2}\log n + r \} \cap \mathscr{S}
_n^c = \mathscr{S}_n^c$ and that $\mathbf{P}\{\mathscr{S}_n^c\} \ge
q- \ee^{-c_{10} n}$ [see~(\ref{survival-survival})],
we obtain, for $0\le r\le c_{24}(n-\ell_0)$,
%
\begin{eqnarray}\label{E(infV):1er}
&& \mathbf{P} \bigl\{ \underV_n > \tfrac{3}{2}\log n + r,
\mathscr{S}_n  \bigr\}
\nonumber\\[-8pt]
\\[-8pt]
&&\qquad \le\ee^{-c_{10} n}+
\ee^{-c_{26} \lfloor r/c_{24}\rfloor} +
\exp \bigl(- n^{-\varepsilon} \lfloor\ee^{c_{25}
\lfloor r/c_{24}\rfloor} \rfloor \bigr).
\nonumber
\end{eqnarray}
This implies that, for any $0< c_{27}< \min\{ {c_{26}\over
c_{24}},   {2c_{10} \over c_{24}}\}$, there exists a constant
\mbox{$c_{28}>0$} such that $\mathbf{E}\{ \ee^{c_{27}  \underV_n}
  {\bf1}_{\{{3/2} \log n < \underV_n \le{c_{24}/2} n\} \cap\mathscr
{S}_n} \}\le c_{28}   n^{c_{29}}$, with $c_{29} := ({3\over2} + {c_{24} \over
c_{25}}\varepsilon) c_{27}$. Therefore,
%
\begin{equation}\label{E(infV):2e}
\mathbf{E} \bigl\{ \ee^{c_{27}  \underV_n}
{\bf1}_{\{\underV_n \le{c_{24}/2} n\} \cap
\mathscr{S}_n}  \bigr\}
\le c_{30}  n^{c_{29}}, \qquad n\ge1,
\end{equation}
where $c_{30} := c_{28} +1$.

On the other hand, letting $\delta_- >0$ be as in (\ref{hyp3}), we
have $\ee^{\delta_-   \underV_n}   {\bf1}_{\mathscr{S}_n}
\le\sum_{|u|=n} \ee^{\delta_- V(u)}$. Since $\psi(-\delta_-)
:=\log\mathbf{E}\{
\sum_{|u|=n} \ee^{\delta_- V(u)} \} <\infty$ by (\ref{hyp3}), we
can choose and fix $c_{31}>0$ sufficiently large (in particular,
$c_{31} > {c_{24}\over2}$) such that, for any $x\ge c_{31}$,
\[
\mathbf{P}\{ \underV_n> x n,   \mathscr{S}_n\}
\le\ee^{-\delta_- x n + \psi
(-\delta_-) n} \le\ee^{-\delta_- x   n/2},
\qquad\forall n\ge1.
\]
Therefore, for any $c_{32} < {\delta_-\over2}$, we have
%
\begin{equation}\label{E(infV):3e}
\sup_{n\ge1} \mathbf{E} \bigl\{ \ee^{c_{32}  \underV_n}
{\bf1}_{\{ \underV_n > c_{31}  n\} \cap
\mathscr{S}_n}  \bigr\} <\infty.
\end{equation}

Finally, (\ref{E(infV):1er}) also implies that, for $n\ge\ell_0$,
\begin{eqnarray*}
\mathbf{P} \biggl\{ \underV_n >{c_{24}\over2}  n,   \mathscr{S}_n
 \biggr\} &\le& \ee
^{-c_{10} n}+ \ee^{-c_{26}  \lfloor{n/2} - {3/ (2c_{24})}\log
n\rfloor}
\\
&&{} + \exp \bigl(- n^{-\varepsilon} \lfloor\ee^{c_{25}
\lfloor{n/2} - {3/(2c_{24})}\log n \rfloor}\rfloor \bigr).
\end{eqnarray*}
Therefore, for any $c_{33} < \min\{ {c_{10}\over c_{31}},
{c_{26}\over2c_{31}}\}$,
\[
\sup_{n\ge1} \mathbf{E} \bigl\{ \ee^{c_{33}  \underV_n}
{\bf1}_{\{{c_{24}/2}   n < \underV_n \le c_{31}  n\} \cap\mathscr{S}_n}
 \bigr\}<\infty,
\]
which, combined with (\ref{E(infV):2e}) and (\ref{E(infV):3e}),
completes the proof of Proposition \ref{p:mcdiarmid-integrability},
with $c_{22} := \min\{ c_{27},   c_{32},   c_{33}\}$.
\end{pf*}
\begin{lemma}\label{l:max-ineq}
Let $X_1, X_2,\ldots, X_N$ be independent
nonnegative random variables, and let
$T_N := \sum_{i=1}^N X_i$. For any nonincreasing function
$F\dvtx (0,   \infty) \to\mathbb{R}_+$, we have
\[
\mathbf{E} \bigl\{ F(T_N)   {\bf1}_{ \{ T_N>0\} } \bigr\}
\le\max_{1\le i\le N}\mathbf{E} \{ F(X_i)   |   X_i >0  \} .
\]
Moreover,
\[
\mathbf{E} \bigl\{ F(T_N)   {\bf1}_{ \{ T_N>0\} } \bigr\} \le
\sum_{i=1}^N b^{i-1} \mathbf{E}
\bigl\{F(X_i)   {\bf1}_{ \{X_i >0\} }  \bigr\} ,
\]
where $b:= \max_{1\le i\le N} \mathbf{P}\{ X_i=0\}$.
\end{lemma}
\begin{pf}
Let $\tau:= \min\{ i\ge1:   X_i >0\}$ (with $\min\varnothing:= \infty$). Then
$\mathbf{E} \{ F(T_N)  \times\break  {\bf1}_{ \{ T_N>0\} } \}
= \sum_{i=1}^N \mathbf{E} \{ F(T_N)   {\bf1}_{ \{\tau=i\}}  \}$.
Since $F$ is nonincreasing, we have $F(T_N)   \times\break {\bf1}_{\{ \tau=i\} }
\le F(X_i)   {\bf1}_{ \{ \tau=i\} } = F(X_i)
{\bf1}_{ \{ X_i>0\} }   {\bf1}_{ \{ X_j=0,   \forall j<i\} }$. By
independence, this leads to
\[
\mathbf{E} \bigl\{ F(T_N)   {\bf1}_{ \{ T_N>0\} }  \bigr\}
\le\sum_{i=1}^N \mathbf{E} \bigl\{ F(X_i)   {\bf1}_{ \{ X_i>0\} }  \bigr\}
\mathbf{P} \{ X_j=0,   \forall j<i  \}.
\]
This yields immediately the second inequality of the lemma,
since $\mathbf{P}\{ X_j=0,   \forall j<i \} \le b^{i-1}$.

To prove the first inequality of the lemma, we observe that
$\mathbf{E}\{F(X_i)   {\bf1}_{ \{ X_i>0\} } \} \le\mathbf{P}\{ X_i>0\}
\max_{1\le k\le N} \mathbf{E}\{ F(X_k)   |   X_k >0 \}$. Therefore,
\[
\mathbf{E} \bigl\{ F(T_N)   {\bf1}_{ \{ T_N>0\} }  \bigr\}
\le\max_{1\le k\le N}
\mathbf{E}\{ F(X_k)   |   X_k >0 \} \sum_{i=1}^N
\mathbf{P}\{ X_i>0\} \mathbf{P} \{ X_j=0,   \forall j<i  \} .
\]
The $\sum_{i=1}^N \cdots$ expression on the right-hand side
is $=\sum_{i=1}^N \mathbf{P}\{ X_i>0,   X_j=0, \break   \forall j<i \} =
\sum_{i=1}^N \mathbf{P}\{ \tau=i\} = \mathbf{P}\{ T_N>0\} \le1$. This
yields the first inequality of the lemma.
\end{pf}

To state our third estimate, let $\underline{w}^{(n)} \in\llbracket e,
w_n^{(n)}\rrbracket $ be a vertex such that
%
\begin{equation}\label{minimal-vertex}
V\bigl({\underline{w}^{(n)}}\bigr) = \min_{u\in\llbracket e,
w_n^{(n)}\rrbracket } V(u) .
\end{equation}
[If there are several such vertices, we choose, say, the
oldest.] The following estimate gives a (stochastic) lower bound for
${1\over W_{n,\beta}}$ under $\mathbf{Q}$ outside a ``small'' set. We recall
that $W_{n,\beta}>0$, $\mathbf{Q}$-almost surely (but not necessarily
$\mathbf{P} $-almost surely).
\begin{lemma}\label{l:survival}
Assume (\ref{hyp2}), (\ref{hyp3}) and (\ref{hyp}).
For any $K>0$, there exist $\theta>0$ and
$n_0<\infty$ such that, for any $n\ge n_0$, any
$\beta>0$, and any nondecreasing function $G\dvtx (0,\infty) \to(0,\infty)$,
%
\begin{equation}\label{GG}
\mathbf{E}_\mathbf{Q} \biggl\{ G \biggl(
{\ee^{-\beta V({\underline{w}^{(n)}})} \over
W_{n,\beta}} \biggr) {\bf1}_{E_n}  \biggr\} \le
{1\over1-q} \max_{0\le k<n}
\mathbf{E} \biggl\{ G \biggl({n^{\theta\beta} \over W_{k,\beta}} \biggr)
{\bf1}_{\mathscr{S}_k}  \biggr\} ,
\end{equation}
where $E_n$ is a measurable event such that
\[
\mathbf{Q}\{ E_n\} \ge1 - {1\over n^K}, \qquad n\ge n_0.
\]
\end{lemma}
%
%
%
\begin{pf}
Recall from (\ref{Ck}) that $\mathscr{I}_k^{(n)}$
is the set of the brothers of $w_k^{(n)}$. For any pair
$0\le k< n$, we say that \textit{the level $k$ is $n$-good} if
\[
\mathscr{I}_k^{(n)} \not= \varnothing \quad
\mbox{and}\quad
\mathbb{T}^{GW}_u\mbox{ survives at least
$n-k$ generations},\qquad \forall u\in \mathscr{I}_k^{(n)} ,
\]
where $\mathbb{T}^{\hbox{\scriptsize GW}}_u$ is the shifted
Galton--Watson subtree generated by $u$ [see (\ref{shifted-GW})]. By
$\mathbb{T}^{\mathrm{GW}}_u$ surviving at least $n-k$
generations, we
mean that there exists $v\in\mathbb{T}^{\mathrm{GW}}_u$ such that
$|v|_u=n-k$ [see (\ref{|x|_u}) for notation].

In words, $k$ is $n$-good means any subtree generated by any of the
brothers of $w_k^{(n)}$ has offspring for at least $n-k$ generations.

Let $\mathscr{G}_n$ be the sigma-field defined in (\ref{Gn}). By
Proposition \ref{p:change-proba},
\[
\mathbf{Q} \{ \mbox{$k$ is $n$-good}   |   \mathscr{G}_n \}
= {\bf1}_{ \{ \mathscr{I}_k^{(n)} \not= \varnothing\} }
(\mathbf{P}\{\mathscr{S}_{n-k}\} )^{\# \mathscr{I}_k^{(n)}},
\]
where $\mathscr{S}_n$ denotes the system's survival after
$n$ generations [see (\ref{Sn})]. Since
$\mathbf{P}\{\mathscr{S}_{n-k}\} \ge \mathbf{P}\{ \mathscr{S}\} =1-q$,
whereas $\# \mathscr{I}_k^{(n)}$ and
$\#\mathscr{I}_1^{(1)}$ have the same distribution under $\mathbf{Q}$
(Proposition \ref{p:change-proba}), we have
\[
\mathbf{Q} \{ \mbox{$k$ is $n$-good}  \}
\ge\mathbf{E}_\mathbf{Q} \bigl\{ {\bf1}_{ \{ \# \mathscr{I}_1^{(1)}
\ge1 \} }  ( 1-q)^{\# \mathscr{I}_1^{(1)}}  \bigr\}
= c_{34} \in(0,   1).
\]
As a consequence, for all $1\le\ell<n$, by Proposition
\ref{p:change-proba} again,
\begin{eqnarray*}
\mathbf{Q} \Biggl\{ \bigcup_{k=1}^n
\bigcap_{j:   1\le j\le n,   |j-k| \le\ell}
\{ \mbox{$j$ is not $n$-good} \}  \Biggr\}
&\le& \sum_{k=1}^n \prod_{j:   1\le j\le n,   |j-k| \le\ell}
\mathbf{Q} \{ \mbox{$j$ is not $n$-good} \}
\\
&\le& n(1-c_{34})^{\ell+1},
\end{eqnarray*}
which is bounded by $n \ee^{-c_{34} (\ell+1)}$ (using the
inequality $1-x\le\ee^{-x}$, for $x\ge0$). Let $K>0$. We take
$\ell=\ell(n):=\lfloor c_{35} \log n\rfloor$ with $c_{35}:={K+2\over
c_{34}}$. Let $c_{36} := {K+2\over c_6}$ [where $c_6$ is as in
(\ref{tail-sup|V|})] and $c_{37} := \max\{ {K+2\over c_3},   c_{35}\}$
[$c_3$ being the constant in (\ref{Petrov})]. Let
%
\begin{eqnarray}
\quad
E_n^{(1)} &:=&\bigcap_{k=1}^n
  \bigcup_{j:   1\le j\le n, |j-k| \le\lfloor c_{35} \log n\rfloor}
\{ \mbox{$j$ is $n$-good} \} ,
\label{E1=}\\
E_n^{(2)} &:=&  \biggl\{   \max_{1\le j\le n}
\sup_{u\in\mathscr{I}_j^{(n)}}
\bigl|V(u)-V\bigl(w_{j-1}^{(n)}\bigr)\bigr| \le c_{36}\log n \biggr\},
\label{E2=}\\
E_n^{(3)} &:=& \biggl\{   \max_{0\le j,   k\le n,
|j-k| \le c_{35} \log n}\bigl|V\bigl(w_j^{(n)}\bigr)
- V\bigl(w_k^{(n)}\bigr)\bigr|
\le c_{37}   \log n  \biggr\}.
\label{E3=}
\end{eqnarray}
We have
\[
\mathbf{Q} \bigl\{ E_n^{(1)}  \bigr\}
\ge1- n \ee^{-c_{34} c_{35} \log n} = 1 - {1\over n^{K+1}}.
\]
On the other hand, by Corollary \ref{c:change-proba},
\[
\mathbf{Q}\bigl\{ \bigl(E_n^{(2)}\bigr)^c\bigr\} \le n  \mathbf{Q}
\biggl\{\sup_{u\in\mathscr{I}_1^{(1)}} |V(u)| > c_{36}   \log n  \biggr\} \le
n  \mathbf{Q} \biggl\{
\sup_{|u|=1} |V(u)| > c_{36}   \log n  \biggr\}.
\]
Applying (\ref{tail-sup|V|}) yields that
\[
\mathbf{Q}\bigl\{ E_n^{(2)} \bigr\}
\ge1- c_5   n^{-(c_{36}c_6-1)} = 1- {c_5\over n^{K+1}}.
\]

To estimate $\mathbf{Q}\{ E_n^{(3)} \}$, we note that,
by Corollary \ref{c:change-proba},
\[
\mathbf{Q}\bigl\{ \bigl(E_n^{(3)}\bigr)^c\bigr\}
= \mathbf{Q} \biggl\{   \max_{0\le j, k\le n,
|j-k| \le c_{35} \log n} |S_j-S_k| > c_{37}   \log n \biggr\} ,
\]

\noindent which, in view of (\ref{tail-|Sj-Sk|}), is bounded by
$2c_{35} n^{-(c_3c_{37}-1)}\log n$. Consequently, if
%
\begin{equation}\label{En}
E_n := E_n^{(1)} \cap E_n^{(2)} \cap E_n^{(3)},
\end{equation}
then $\mathbf{Q}\{E_n\} \ge1-{1\over n^K}$ for all large $n$.\vadjust{\goodbreak}

It remains to check (\ref{GG}). By definition,
%
\begin{eqnarray}\label{Wn-decomp}
W_{n,\beta} &=& \sum_{j=1}^n \sum_{u\in\mathscr{I}_j^{(n)}}
\ee^{-\beta V(u)}
\sum_{x\in\mathbb{T}^{\mathrm{GW}}_u,   |x|_u=n-j}
\ee^{-\beta V_u(x)}
+ \ee^{-\beta V(w_n^{(n)})}
\nonumber\\[-8pt]
\\[-8pt]
&\ge& \sum_{j\in\mathscr{L}}
\sum_{u\in\mathscr{I}_j^{(n)}}
\ee^{-\beta V(u)}
\sum_{x\in\mathbb{T}^{\mathrm{GW}}_u,   |x|_u=n-j}
\ee^{-\beta V_u(x)}
\nonumber
\end{eqnarray}
for any $\mathscr{L} \subset\{1,2,\ldots, n\}$. We choose
$\mathscr{L} := \{ 1\le j\le n\dvtx   | j-|{\underline{w}^{(n)}}| | < c_{35} \log n\}$.

On the event $E_n$, for $u\in\mathscr{I}_j^{(n)}$ with some $j\in
\mathscr{L}$, we have $V(u) \le V({\underline{w}^{(n)}}) + (c_{36} +
c_{37}) \log n$. Writing $\theta:= c_{36} + c_{37}$, this leads to
$W_{n,\beta} \ge n^{-\theta\beta} \ee^{-\beta V({\underline{w}^{(n)}})}
\times\break \sum_{j\in\mathscr{L}} \sum_{u\in
\mathscr{I}_j^{(n)}} \xi_u$, where
\[
\xi_u := \sum_{x\in\mathbb{T}^{\mathrm{GW}}_u,   |x|_u=n-j}
\ee^{-\beta V_u(x)}.
\]
Since $\sum_{j\in\mathscr{L}} \sum_{u\in\mathscr{I}_j^{(n)}} \xi_u >0$
on $E_n$, we arrive at
\[
{\ee^{-\beta V({\underline{w}^{(n)}})} \over W_{n,\beta}}
{\bf1}_{E_n} \le{n^{\theta\beta} \over\sum_{j\in\mathscr{L}}
\sum_{u\in\mathscr{I}_j^{(n)}} \xi_u}
{\bf1}_{ \{ \sum_{j\in\mathscr{L}} \sum_{u\in\mathscr{I}_j^{(n)}} \xi_u >0 \} } .
\]
Let $\mathscr{G}_n$ be the sigma-field in (\ref{Gn}). We\vspace*{2pt}
observe that $\mathscr{L}$ and $\mathscr{I}_j^{(n)}$ are
$\mathscr{G}_n$-measurable. Moreover, an application of Proposition
\ref{p:change-proba} tells us that under $\mathbf{Q}$, conditionally on
$\mathscr{G}_n$, the random variables $\xi_u$, $u\in\mathscr{I}_j^{(n)}$,
$j\in\mathscr{L}$, are independent, and are distributed as $W_{n-j,
\beta}$ under $\mathbf{P}$. We are thus entitled to apply
Lemma \ref{l:max-ineq} to see that, if $G$ is nondecreasing,
\begin{eqnarray*}
\mathbf{E}_\mathbf{Q} \biggl\{ G \biggl(
{\ee^{-\beta V({\underline{w}^{(n)}})}
\over W_{n,\beta}}  \biggr) {\bf1}_{E_n}
   \big|   \mathscr{G}_n \biggr\}
&\le& \max_{j\in\mathscr{L}}
\mathbf{E} \biggl\{ G \biggl(
{n^{\theta\beta}\over W_{n-j, \beta}}  \biggr)
\big |   W_{n-j, \beta}>0  \biggr\}
\\
&\le& \max_{0\le k< n} \mathbf{E} \biggl\{ G \biggl(
{n^{\theta\beta}\over W_{k, \beta}}  \biggr)
\big |   W_{k, \beta}>0  \biggr\} .
\end{eqnarray*}
Since $\mathbf{P}\{ W_{k,\beta}>0\} =\mathbf{P}\{
\mathscr{S}_k\} \ge\mathbf{P}\{ \mathscr{S}\} = 1-q$,
this yields Lemma \ref{l:survival}.
\end{pf}

We are now ready for the proof of Theorem
\ref{t:derrida-spohn-moment}. For the sake of clarity, the upper and lower
bounds are proved in distinct parts. Let us start with the upper bound.
\begin{pf*}{Proof of Theorem \ref{t:derrida-spohn-moment}}
\textit{The upper bound}.
We assume (\ref{hyp2}), (\ref{hyp3}) and~(\ref{hyp}), and fix $\beta>1$.

For any $Z\ge0$ which is $\mathscr{F}_n$-measurable, we have
$\mathbf{E}\{W_{n,\beta} Z\} =\break  \mathbf{E}_\mathbf{Q}\{ \sum_{|u|=n}
{\ee^{-\beta V(u)} \over W_n} Z\} =\mathbf{E}_\mathbf{Q}
\{ \sum_{|u|=n} {\bf1}_{\{ w_n^{(n)}=u\} }\ee^{-(\beta-1) V(u)} Z\}$
and, thus,
%
\begin{equation}\label{E(WnZ)}
\mathbf{E}\{ W_{n,\beta} Z\} =
\mathbf{E}_\mathbf{Q}\bigl\{ \ee^{-(\beta-1) V(w_n^{(n)})} Z\bigr\} .
\end{equation}

Let $s\in({\beta-1\over\beta},   1)$, and $\lambda>0$. (We will
choose $\lambda= {3\over2}$.) Then
\begin{eqnarray*}
\mathbf{E} \{ W_{n,\beta}^{1-s}  \}
&\le& n^{-(1-s)\beta\lambda} + \mathbf{E} \bigl\{ W_{n,\beta}^{1-s}
{\bf1}_{ \{ W_{n,\beta} > n^{-\beta\lambda}\} }  \bigr\}
\\
&=& n^{-(1-s)\beta\lambda} + \mathbf{E}_\mathbf{Q} \biggl\{
{\ee^{-(\beta-1)V(w_n^{(n)})} \over
W_{n,\beta}^s}{\bf1}_{ \{ W_{n,\beta} > n^{-\beta\lambda}\} }  \biggr\} .
\end{eqnarray*}
Since $\ee^{-\beta V(w_n^{(n)})} \le W_{n,\beta}$,
we have
${\ee^{-(\beta-1)V(w_n^{(n)})} \over W_{n,\beta}^s} \le
{1\over W_{n,\beta}^{s-(\beta-1)/\beta}}$; thus, on the event
$\{ W_{n,\beta} > n^{-\beta\lambda} \}$, we have
${\ee^{-(\beta-1)V(w_n^{(n)})}\over W_{n,\beta}^s}
\le n^{[\beta s - (\beta-1)] \lambda}$.

Let $K:=[\beta s - (\beta-1)] \lambda+ (1-s) \beta\lambda$, and let
$E_n$ be the event in Lemma \ref{l:survival}. Since
$\mathbf{Q}(E_n^c) \le n^{-K}$ for
all sufficiently large $n$ (see Lemma \ref{l:survival}), we obtain, for large $n$,
%
\begin{eqnarray}\label{pizza}
\mathbf{E} \{ W_{n,\beta}^{1-s}  \}
&\le& n^{-(1-s)\beta\lambda} + n^{[\beta s - (\beta-1)] \lambda- K}
\nonumber \\
&&{} + \mathbf{E}_\mathbf{Q} \biggl\{ {\ee^{-(\beta-1)V(w_n^{(n)})}
\over W_{n,\beta}^s}
{\bf1}_{ \{ W_{n,\beta} > n^{-\beta\lambda}
\} \cap E_n}  \biggr\}
\\
&=& 2n^{-(1-s)\beta\lambda} +\mathbf{E}_\mathbf{Q}
\biggl\{ {\ee^{-(\beta-1)V(w_n^{(n)})}
\over W_{n,\beta}^s}
{\bf1}_{ \{ W_{n,\beta} > n^{-\beta\lambda}
\} \cap E_n}  \biggr\} .
\nonumber
\end{eqnarray}

We now estimate the expectation expression
$\mathbf{E}_\mathbf{Q}\{\cdot\}$ on the
right-hand side. Let $a>0$ and $\varrho> b>0$ be constants such that
$(\beta-1) a> s\beta\lambda+ {3\over2}$ and $[\beta s - (\beta-1)]
b > {3\over2}$. (The choice of $\varrho$ will be made precise later
on.) We recall that $\underline{w}_n^{(n)}\in\llbracket e,   w_n^{(n)}\rrbracket $
satisfies $V({\underline{w}^{(n)}})=
\min_{u\in\llbracket e,   w_n^{(n)}\rrbracket }V(u)$, and consider the following events:
\begin{eqnarray*}
E_{1,n} &:=&  \bigl\{ V\bigl(w_n^{(n)}\bigr) > a \log n \bigr\}
\cup \bigl\{ V\bigl(w_n^{(n)}\bigr) \le- b\log n  \bigr\} ,
\\
E_{2,n} &:=&  \bigl\{ V\bigl({\underline{w}^{(n)}}\bigr)
< -\varrho \log n,   V\bigl(w_n^{(n)}\bigr) > -b\log n\bigr\},
\\
E_{3,n} &:=&  \bigl\{ V\bigl({\underline{w}^{(n)}}\bigr) \ge
- \varrho\log n,   -b\log n <V\bigl(w_n^{(n)}\bigr) \le a \log n \bigr\} .
\end{eqnarray*}
Clearly, $\mathbf{Q}(\bigcup_{i=1}^3 E_{i,n})=1$.

On the event $E_{1,n}\cap\{ W_{n,\beta} > n^{-\beta\lambda} \}$, we
have either $V(w_n^{(n)}) > a \log n$, in which case
${\ee^{-(\beta-1)V(w_n^{(n)})} \over W_{n,\beta}^s}
\le n^{s\beta\lambda- (\beta
-1) a}$, or $V(w_n^{(n)}) \le- b\log n$, in which case we use the
trivial inequality $W_{n,\beta} \ge\ee^{-\beta V(w_n^{(n)})}$
to see
that ${\ee^{-(\beta-1)V(w_n^{(n)})} \over W_{n,\beta}^s} \le
\ee
^{[\beta s - (\beta-1)] V(w_n^{(n)})} \le n^{-[\beta s - (\beta-1)]
b}$ (recalling that $\beta s >\beta-1$). Since $s\beta\lambda-
(\beta-1) a < -{3\over2}$ and $[\beta s - (\beta-1)] b > {3\over 2}$, we obtain
%
\begin{equation}\label{E1}
\mathbf{E}_\mathbf{Q} \biggl\{ {\ee^{-(\beta-1)V(w_n^{(n)})}
\over W_{n,\beta}^s}   {\bf1}_{ E_{1,n} \cap
\{ W_{n,\beta} > n^{-\beta\lambda} \}}\biggr\} \le n^{- 3/2}.
\end{equation}

We now study the integral on $E_{2,n}\cap\{ W_{n,\beta} > n^{-\beta
\lambda} \} \cap E_n$. Since $s>0$, we can choose $s_1>0$ and
$0<s_2\le{c_{22}\over\beta}$ [where $c_{22}$ is the constant in
(\ref{mcdiarmid-integrability})] such that $s=s_1+s_2$. We have, on
$E_{2,n}\cap\{ W_{n,\beta} > n^{-\beta\lambda} \}$,
\begin{eqnarray*}
{\ee^{-(\beta-1)V(w_n^{(n)})} \over W_{n,\beta}^s}
&=& {\ee^{\beta s_2 V({\underline{w}^{(n)}}) - (\beta-1)V(w_n^{(n)})}
\over W_{n,\beta}^{s_1}} {\ee^{-\beta s_2 V({\underline{w}^{(n)}})}
\over W_{n,\beta}^{s_2}}
\\
&\le& n^{-\beta s_2 \varrho+ (\beta-1) b
+ \beta \lambda s_1} {\ee^{-\beta s_2 V({\underline{w}^{(n)}})}
\over W_{n,\beta}^{s_2}} .
\end{eqnarray*}
Therefore, by an application of Lemma \ref{l:survival} to
$G(x) := x^{s_2}$, $x>0$, we obtain, for all sufficiently large $n$,
\begin{eqnarray*}
&& \mathbf{E}_\mathbf{Q} \biggl\{ {\ee^{-(\beta-1)V(w_n^{(n)})}
\over W_{n,\beta}^s} {\bf1}_{E_{2,n}\cap\{ W_{n,\beta}
> n^{-\beta\lambda} \} \cap E_n}
 \biggr\}
\\
&&\qquad \le{n^{-\beta s_2 \varrho+ (\beta-1) b + \beta\lambda
s_1}\over1-q} \max_{0\le k<n} \mathbf{E} \biggl\{ {n^{s_2 \theta\beta
} \over W_{k,\beta}^{s_2} } {\bf1}_{\mathscr{S}_k} \biggr\} .
\end{eqnarray*}
By definition, ${1\over W_{k,\beta}^{s_2}} \le\exp(\beta
s_2 \inf_{|x|=k} V(x))$; thus, by (\ref{mcdiarmid-integrability}),
$\mathbf{E}\{ {n^{s_2 \theta\beta} \over W_{k,\beta}^{s_2} }
{\bf1}_{\mathscr{S}_k} \} \le c_{23}^{\beta s_2/c_{22}} n^{s_2 \theta\beta
+{(3+\varepsilon)/2}\beta s_2}$ for all $0\le k<n$. We choose (and
fix) the constant $\varrho$ so large that $-\beta s_2 \varrho+ (\beta
-1) b + \beta\lambda s_1 + s_2 \theta\beta+ {3+\varepsilon\over
2}\beta s_2 < - {3\over2}$. Therefore, for all large $n$,
%
\begin{equation}\label{E2}
\mathbf{E}_\mathbf{Q} \biggl\{ {\ee^{-(\beta-1)V(w_n^{(n)})}
\over W_{n,\beta}^s} {\bf1}_{E_{2,n}\cap
\{ W_{n,\beta} > n^{-\beta\lambda} \} \cap E_n} \biggr\} \le n^{-3/2}.
\end{equation}

We make a partition of $E_{3,n}$: let $M\ge2$ be an integer, and let
$a_i := -b + {i(a+b)\over M}$, $0\le i\le M$. By definition,
\begin{eqnarray*}
E_{3,n} &=& \bigcup_{i=0}^{M-1}  \bigl\{ V\bigl({\underline{w}^{(n)}}\bigr)
\ge- \varrho\log n,   a_i \log n <V\bigl(w_n^{(n)}\bigr)
\le a_{i+1} \log n \bigr\}
\\
&=&\!\!: \bigcup_{i=0}^{M-1} E_{3,n,i}.
\end{eqnarray*}
Let $0\le i\le M-1$. There are two possible situations. First
situation: $a_i \le\lambda$. In this case, we argue that, on the event
$E_{3,n,i}$, we have $W_{n,\beta} \ge\ee^{-\beta V(w_n^{(n)})}
\ge n^{-\beta a_{i+1}}$ and $\ee^{-(\beta-1)V(w_n^{(n)})}
\le n^{-(\beta-1) a_i}$, thus, ${\ee^{-(\beta-1)V(w_n^{(n)})}
\over W_{n,\beta}^s}
\le n^{\beta s a_{i+1} -(\beta-1) a_i} = n^{\beta s a_i -(\beta-1)
a_i + \beta s(a+b)/M} \le n^{[\beta s-(\beta-1) ] \lambda+ \beta
s(a+b)/M}$. Accordingly, in this situation,
\[
\mathbf{E}_\mathbf{Q}
\biggl\{ {\ee^{-(\beta-1)V(w_n^{(n)})}
\over W_{n,\beta}^s}{\bf1}_{ E_{3,n,i} } \biggr\}
\le n^{[\beta s-(\beta-1)] \lambda+ \beta s(a+b)/M} \mathbf{Q}(E_{3,n,i}) .
\]
Second (and last) situation: $a_i > \lambda$. We have, on
$E_{3,n,i} \cap\{ W_{n,\beta} > n^{-\beta\lambda} \}$, ${\ee
^{-(\beta-1)V(w_n^{(n)})} \over W_{n,\beta}^s} \le n^{\beta\lambda s
- (\beta-1)a_i} \le n^{[\beta s - (\beta-1)]\lambda}$; thus, in this
situation,
\[
\mathbf{E}_\mathbf{Q} \biggl\{ {\ee^{-(\beta-1)V(w_n^{(n)})}
\over W_{n,\beta}^s}
{\bf1}_{ E_{3,n,i} \cap\{ W_{n,\beta} > n^{-\beta\lambda} \} }
\biggr\} \le n^{[\beta s - (\beta-1)]\lambda} \mathbf{Q}(E_{3,n,i}).
\]
We have therefore proved that
\begin{eqnarray*}
&& \mathbf{E}_\mathbf{Q} \biggl\{ {\ee^{-(\beta-1)V(w_n^{(n)})}
\over W_{n,\beta}^s}   {\bf1}_{ E_{3,n} \cap
\{ W_{n,\beta} > n^{-\beta\lambda} \} } \biggr\}
\\
&&\qquad = \sum_{i=0}^{M-1} \mathbf{E}_\mathbf{Q}
\biggl\{{\ee^{-(\beta-1)V(w_n^{(n)})} \over
W_{n,\beta}^s}   {\bf1}_{ E_{3,n,i} \cap
\{ W_{n,\beta} > n^{-\beta\lambda} \} } \biggr\}
\\
&&\qquad \le n^{[\beta s-(\beta-1)]\lambda
+ \beta s(a+b)/M} \mathbf{Q}(E_{3,n}) .
\end{eqnarray*}
By Corollary \ref{c:change-proba},
$\mathbf{Q}(E_{3,n}) =\mathbf{P}\{ \min_{0\le k \le n} S_k
\ge- \varrho\log n,   -b\log n\le S_n \le a\times\break \log n\} = n^{- (3/2)+o(1)}$.
Combining this with (\ref{pizza}), (\ref{E1}) and (\ref{E2}) yields
\[
\mathbf{E} \{ W_{n,\beta}^{1-s}  \} \le2n^{-(1-s)\beta \lambda}
+ 2n^{-3/2} + n^{[\beta s-(\beta-1)] \lambda+ \beta s(a+b)/M -(3/2)
+o(1)} .
\]
We choose $\lambda:= {3\over2}$. Since $M$ can be as large
as possible, this yields the upper bound in Theorem
\ref{t:derrida-spohn-moment} by posing $r:= 1-s$.
\end{pf*}
\begin{pf*}{Proof of Theorem \ref{t:derrida-spohn-moment}}
\textit{The lower bound}.
Assume (\ref{hyp2}), (\ref{hyp3}) and (\ref{hyp}).
Let $\beta>1$ and $s\in(1- {1\over\beta},   1)$. By means
of (\ref{Wn-decomp}) and the elementary inequality
$(a+b)^{1-s} \le a^{1-s} + b^{1-s}$ (for $a\ge0$ and $b\ge0$), we have
\begin{eqnarray*}
W_{n,\beta}^{1-s} &\le& \sum_{j=1}^n \sum_{u\in\mathscr{I}_j^{(n)}}
\ee^{-(1-s)\beta V(u)}
\Biggl ( \sum_{x\in\mathbb{T}^{\mathrm{GW}}_u,
|x|_u=n-j} \ee^{-\beta V_u(x)}  \Biggr)^{1-s}
+ \ee^{- (1-s)\beta V(w^{(n)}_n)}
\\
&=& \sum_{j=1}^n
\ee^{-(1-s)\beta V(w_{j-1}^{(n)})}
\sum_{u\in\mathscr{I}_j^{(n)}}
\ee^{-(1-s)\beta[V(u)-V(w_{j-1}^{(n)})]}
\\
&&\hspace*{114pt}{}\times
\Biggl ( \sum_{x\in\mathbb{T}^{\mathrm{GW}}_u,
|x|_u=n-j} \ee^{-\beta V_u(x)}  \Biggr)^{1-s}
\\
&&{} + \ee^{-(1-s)\beta V(w^{(n)}_n)} .
\end{eqnarray*}
Let $\mathscr{G}_n$ be the sigma-field defined in
(\ref{Gn}), and let
\[
\Xi_j = \Xi_j(n,s,\beta) := \sum_{u\in\mathscr{I}_j^{(n)}}
\ee^{-(1-s)\beta[V(u)-V(w_{j-1}^{(n)})]}, \qquad 1\le j\le n.
\]
Since $V(w_j^{(n)})$ and $\mathscr{I}_j^{(n)}$, for $1\le
j\le n$, are $\mathscr{G}_n$-measurable, it follows from Proposition
\ref{p:change-proba} that
\[
\mathbf{E}_\mathbf{Q} \{ W_{n,\beta}^{1-s}   |
\mathscr{G}_n \} \le
\sum_{j=1}^n \ee^{-(1-s)\beta V(w_{j-1}^{(n)})}  \Xi_j
\mathbf{E}\{W_{n-j,\beta}^{1-s}\}
+ \ee^{-(1-s)\beta V(w^{(n)}_n)}.
\]
Let $\varepsilon>0$ be small, and let
$r:= {3\over 2}(1-s)\beta- \varepsilon$. By means of the already proved upper
bound for $\mathbf{E}(W_{n,\beta}^{1-s})$, this leads to, with
$c_{38}\ge1$,
%
\begin{eqnarray}\label{E(Wnbeta)=}
\quad
&& \mathbf{E}_\mathbf{Q} \{ W_{n,\beta}^{1-s}   |
\mathscr{G}_n \}
\nonumber\\[-8pt]
\\[-8pt]
&&\qquad \le c_{38} \sum_{j=1}^n
\ee^{-(1-s)\beta V(w_{j-1}^{(n)})} (n-j+1)^{-r}
  \Xi_j + \ee^{-(1-s)\beta V(w^{(n)}_n)}.
\nonumber
\end{eqnarray}

Since $\mathbf{E}( W_{n, \beta}^{1-s} ) = \mathbf{E}_\mathbf{Q}\{
{\ee^{-(\beta-1)
V(w_n^{(n)})} \over W_{n,\beta}^s} \}$ [see (\ref{E(WnZ)})], we have,
by Jensen's inequality [noticing that $V(w_n^{(n)})$ is
$\mathscr{G}_n$-measurable],
\[
\mathbf{E} ( W_{n, \beta}^{1-s} ) \ge\mathbf{E}_\mathbf{Q}
\biggl\{ {\ee^{-(\beta-1) V(w_n^{(n)})} \over
\{ \mathbf{E}_\mathbf{Q}( W_{n,\beta}^{1-s}   |
  \mathscr{G}_n )\}^{s/(1-s)} }  \biggr\} ,
\]
which, in view of (\ref{E(Wnbeta)=}), yields
\begin{eqnarray*}
\mathbf{E} ( W_{n,\beta}^{1-s} )
& \ge& {1\over c_{38}^{s/(1-s)}}
\\
&&{}\times
\mathbf{E}_\mathbf{Q} \Biggl\{ \bigl(\ee^{-(\beta-1) V(w_n^{(n)})}\bigr)
\\
&&\hspace*{32pt}{}\times
\Biggl(\Biggl\{\sum_{j=1}^n \ee^{-(1-s)\beta V(w_{j-1}^{(n)})}
(n-j+1)^{-r}  \Xi_j
\\
&&\hspace*{143pt}{} + \ee^{-(1-s)
\beta V(w^{(n)}_n)} \Biggr\}^{s/(1-s)}\Biggr)^{-1}  \Biggr\}.
\end{eqnarray*}
By Proposition \ref{p:change-proba}, if $(S_j-S_{j-1},
\xi_j)$, for $j\ge1$ (with $S_0:=0$), are i.i.d.\ random variables
under $\mathbf{Q}$ and distributed as $(V(w^{(1)}_1),
\sum_{u\in\mathscr{I}_1^{(1)}} \ee^{-(1-s)\beta V(u)})$, then
the $\mathbf{E}_\mathbf{Q}\{\cdot\}$
expression on the right-hand side is
\begin{eqnarray*}
&=& \mathbf{E}_\mathbf{Q} \biggl\{ {\ee^{-(\beta-1) S_n}
\over\{ \sum_{j=1}^n (n-j+1)^{-r}
\ee^{-(1-s) \beta S_{j-1}} \xi_j
+ \ee^{-(1-s)\beta S_n} \}^{s/(1-s)}} \biggr\}
\\
&=& \mathbf{E}_\mathbf{Q} \biggl\{
{\ee^{[\beta s-(\beta-1)]\widetilde{S}_n}
\over\{ \sum_{k=1}^n k^{-r}
\ee^{(1-s) \beta\widetilde{S}_k}
\widetilde{\xi}_k +1\}^{s/(1-s)}}  \biggr\},
\end{eqnarray*}
where
\[
\widetilde{S}_\ell:= S_n - S_{n-\ell}, \qquad
\widetilde{\xi}_\ell:= \xi_{n+1-\ell}, \qquad1\le\ell\le n.
\]
Consequently,
\[
\mathbf{E} ( W_{n,\beta}^{1-s} ) \ge{1\over c_{38}^{s/(1-s)}} \mathbf{E}
_\mathbf{Q} \biggl\{ {\ee^{[\beta s-(\beta-1)]\widetilde{S}_n}
\over\{ \sum_{k=1}^n k^{-r} \ee^{(1-s) \beta\widetilde{S}_k}
\widetilde{\xi}_k +1 \}^{s/(1-s)}}  \biggr\}.
\]

Let $c_{39}>0$ be a constant, and define
\begin{eqnarray*}
E^{\widetilde{S}}_{n,1}
&:=& \bigcap_{k=1}^{\lfloor n^\varepsilon\rfloor-1}
 \{ \widetilde{S}_k \le - c_{39}  k^{1/3}  \}
\cap \bigl\{ - 2 n^{\varepsilon/2} \le
\widetilde{S}_{\lfloor n^\varepsilon\rfloor}
\le- n^{\varepsilon/2} \bigr\},
\\
E^{\widetilde{S}}_{n,2}
&:=& \bigcap_{k=\lfloor n^\varepsilon\rfloor+1}^{n-\lfloor n^\varepsilon
\rfloor-1}  \{ \widetilde{S}_k \le - [ k^{1/3} \wedge(n-k)^{1/3}]  \}
\cap \bigl\{ - 2 n^{\varepsilon/2} \le
\widetilde{S}_{n-\lfloor n^\varepsilon\rfloor}
\le- n^{\varepsilon/2} \bigr\},
\\
E^{\widetilde{S}}_{n,3}
&:=& \bigcap_{k=n-\lfloor n^\varepsilon\rfloor+1}^{n-1}
\biggl \{ \widetilde{S}_k \le{3\over2}\log n \biggr\}
\cap \biggl\{ {3-\varepsilon\over2}\log n \le
\widetilde{S}_n\le{3\over2} \log n \biggr\} .
\end{eqnarray*}
Let $\rho:= \rho((1-s)\beta)$ be the constant in Corollary
\ref{c:V}, and let
\begin{eqnarray*}
E^{\widetilde{\xi}}_{n,1}
&:=& \bigcap_{k=1}^{\lfloor n^\varepsilon\rfloor}
 \{ \widetilde{\xi}_k \le n^{2\varepsilon/\rho}  \},
\\
E^{\widetilde{\xi}}_{n,2}
&:=& \bigcap_{k=\lfloor n^\varepsilon\rfloor+1}^{
n-\lfloor n^\varepsilon\rfloor}
 \{ \widetilde{\xi}_k \le \ee^{n^{\varepsilon/4}}  \},
\\
E^{\widetilde{\xi}}_{n,3}
&:=& \bigcap_{k=n-\lfloor n^\varepsilon\rfloor+1}^n
 \{ \widetilde{\xi}_k \le n^{2\varepsilon/\rho}  \}.
\end{eqnarray*}
On $\bigcap_{i=1}^3 (E^{\widetilde{S}}_{n,i} \cap
E^{\widetilde{\xi}}_{n,i})$, we have $\sum_{k=1}^n k^{-r} \ee
^{(1-s) \beta\widetilde{S}_k} \widetilde{\xi}_k +1 \le c_{40}
n^{2\varepsilon+ (2\varepsilon/\rho)}$, while $\ee^{[\beta s
-(\beta-1)] \widetilde{S}_n} \ge n^{(3-\varepsilon)[\beta s -(\beta
-1)]/2}$ (recalling that $\beta s > \beta-1$). Therefore, with
$c_{41}:= (2 + {2\over\rho}){s\over1-s}$,
%
\begin{eqnarray}\label{E(Wnbeta)>}
\mathbf{E} ( W_{n,\beta}^{1-s} ) &\ge&
(c_{38}c_{40})^{-s/(1-s)} n^{-c_{41}  \varepsilon}
n^{(3-\varepsilon)[\beta s -(\beta-1)]/2}
\nonumber\\[-8pt]
\\[-8pt]
&&{}\times
\mathbf{Q}\Biggl\{ \bigcap_{i=1}^3 (E^{\widetilde{S}}_{n,i} \cap
E^{\widetilde{\xi }}_{n,i})\Biggr\} .
\nonumber
\end{eqnarray}

We need to bound $\mathbf{Q}(\bigcap_{i=1}^3 (E^{\widetilde{S}}_{n,i}
\cap E^{\widetilde{\xi}}_{n,i}))$ from below. Let $\widetilde{S}_0:=0$.
Note that, under $\mathbf{Q}$, $(\widetilde{S}_\ell-\widetilde
{S}_{\ell-1},  \widetilde{\xi}_\ell)$, $1\le\ell\le n$, are i.i.d.,
distributed as $(S_1,   \xi_1)$. For $j\le n$, let $\widetilde
{\mathscr{G}}_j$ be the sigma-field generated by $(\widetilde{S}_k,
  \widetilde{\xi}_k)$, $1\le k\le j$. Then
 $E^{\widetilde{S}}_{n,1}$, $E^{\widetilde{S}}_{n,2}$,
 $E^{\widetilde{\xi}}_{n,1}$ and $E^{\widetilde{\xi}}_{n,2}$ are
$\widetilde{\mathscr{G}}_{n-\lfloor n^\varepsilon\rfloor}$-measurable, whereas
$E^{\widetilde{\xi}}_{n,3}$ is independent of
$\widetilde{\mathscr{G}}_{n-\lfloor n^\varepsilon\rfloor}$. Therefore,
\begin{eqnarray*}
&& \mathbf{Q}\Biggl(\bigcap_{i=1}^3 (E^{\widetilde{S}}_{n,i} \cap
E^{\widetilde{\xi}}_{n,i})   |
\widetilde{\mathscr{G}}_{n-\lfloor n^\varepsilon\rfloor}\Biggr)
\\
&&\qquad
\ge\bigl[\mathbf{Q}\bigl(E^{\widetilde{S}}_{n,3}   |   \widetilde
{\mathscr{G}}_{n-\lfloor n^\varepsilon\rfloor}\bigr)
+ \mathbf{Q}(E^{\widetilde{\xi}}_{n,3}) -1\bigr]
{\bf1}_{E^{\widetilde{S}}_{n,1} \cap
E^{\widetilde{S}}_{n,2} \cap E^{\widetilde{\xi}}_{n,1} \cap
E^{\widetilde{\xi}}_{n,2}} .
\end{eqnarray*}
We have $c_{42}:= \mathbf{E}_\mathbf{Q}( \xi_1^{\rho})<\infty$ [by
(\ref{existence-moment-Wn})]; thus,
$\mathbf{Q}\{ \xi_1>n^{2\varepsilon/\rho}\} \le c_{42}  n^{-2\varepsilon}$,
which entails $\mathbf{Q}(E^{\widetilde
{\xi}}_{n,3}) = (\mathbf{Q}\{ \xi_1\le n^{2\varepsilon/\rho}\})^{\lfloor
n^\varepsilon\rfloor} \ge(1- c_{42}  n^{-2\varepsilon})^{\lfloor
n^\varepsilon\rfloor} \ge1 - c_{43}   n^{-\varepsilon}$. To
estimate $\mathbf{Q}(E^{\widetilde{S}}_{n,3}   |
\widetilde{\mathscr{G}}_{n-\lfloor n^\varepsilon\rfloor})$,
we use the Markov property
to see that, if $\widetilde{S}_{n-\lfloor n^\varepsilon\rfloor} \in
I_n:= [- 2 n^{\varepsilon/2},   - n^{\varepsilon/2}]$, the
conditional probability is (writing $N:= \lfloor n^\varepsilon\rfloor$)
\begin{eqnarray*}
&\ge&\inf_{z\in I_n} \mathbf{Q} \biggl\{ S_i \le{3\over2}\log n -z,\
\forall 1\le i\le N-1,
\\
&&\hspace*{34pt}  {3-\varepsilon\over2}\log n -z \le S_N \le{3\over
2} \log n -z \biggr\} ,
\end{eqnarray*}
which is greater than $N^{-(1/2)+o(1)}$. Therefore,
\[
\mathbf{Q}\bigl(E^{\widetilde{S}}_{n,3}   |   \widetilde{\mathscr
{G}}_{n-\lfloor n^\varepsilon\rfloor}\bigr) + \mathbf{Q}(E^{\widetilde
{\xi}}_{n,3}) -1 \ge n^{-(\varepsilon/2)+o(1)} - c_{43}
n^{-\varepsilon} = n^{-(\varepsilon/2)+o(1)}.
\]
As a consequence,
%
\begin{equation}\label{i=1...6}
\qquad
\mathbf{Q}\Biggl\{ \bigcap_{i=1}^3 (E^{\widetilde{S}}_{n,i}
\cap E^{\widetilde{\xi}}_{n,i})\Biggr\}
\ge n^{-(\varepsilon/2)+o(1)}
\mathbf{Q}(E^{\widetilde{S}}_{n,1} \cap
E^{\widetilde{S}}_{n,2} \cap
E^{\widetilde{\xi}}_{n,1} \cap
E^{\widetilde{\xi}}_{n,2}).
\end{equation}

To estimate $\mathbf{Q}(E^{\widetilde{S}}_{n,1} \cap
E^{\widetilde{S}}_{n,2} \cap E^{\widetilde{\xi}}_{n,1} \cap
E^{\widetilde{\xi}}_{n,2})$, we condition on
$\widetilde{\mathscr{G}}_{\lfloor n^\varepsilon\rfloor}$, and note that
$E^{\widetilde{S}}_{n,1}$ and $E^{\widetilde{\xi}}_{n,1}$ are
$\widetilde{\mathscr{G}}_{\lfloor n^\varepsilon\rfloor}$-measurable, whereas
$E^{\widetilde{\xi}}_{n,2}$ is independent of
$\widetilde{\mathscr{G}}_{\lfloor n^\varepsilon\rfloor}$.\vspace*{1pt} Since
$\mathbf{Q}(E^{\widetilde{S}}_{n,2}  |
\widetilde{\mathscr{G}}_{\lfloor n^\varepsilon
\rfloor}) \ge n^{-(3-\varepsilon)/2+o(1)}$, whereas
$\mathbf{Q}(E^{\widetilde{\xi}}_{n,2}) = [\mathbf{Q}\{
\xi_1 \le\ee^{n^{\varepsilon/4}}\}]^{n-2\lfloor n^\varepsilon\rfloor}
\ge[1- c_{42}  \ee ^{-\rho n^{\varepsilon/4}} ]^{n-2\lfloor n^\varepsilon\rfloor}
\ge1- \ee^{- n^{\varepsilon/5}}$ (for large $n$), we have
\begin{eqnarray*}
\mathbf{Q}\bigl(E^{\widetilde{S}}_{n,1} \cap
E^{\widetilde{S}}_{n,2} \cap
E^{\widetilde{\xi}}_{n,1} \cap
E^{\widetilde{\xi}}_{n,2}   |
\widetilde{\mathscr{G}}_{\lfloor n^\varepsilon\rfloor}\bigr)
&\ge& \bigl[\mathbf{Q}\bigl(E^{\widetilde{S}}_{n,2}  |
\widetilde{\mathscr{G}}_{\lfloor
n^\varepsilon\rfloor}\bigr)
+ \mathbf{Q}(E^{\widetilde{\xi}}_{n,2})-1\bigr]
{\bf1}_{E^{\widetilde{S}}_{n,1} \cap E^{\widetilde{\xi}}_{n,1}}
\\
&\ge& n^{-(3-\varepsilon)/2+o(1)}
{\bf1}_{E^{\widetilde{S}}_{n,1} \cap E^{\widetilde{\xi}}_{n,1}} .
\end{eqnarray*}
Thus, $\mathbf{Q}(E^{\widetilde{S}}_{n,1} \cap
E^{\widetilde{S}}_{n,2} \cap E^{\widetilde{\xi}}_{n,1} \cap
E^{\widetilde{\xi}}_{n,2}) \ge n^{-(3-\varepsilon)/2+o(1)}
\mathbf{Q} (E^{\widetilde{S}}_{n,1} \cap E^{\widetilde{\xi}}_{n,1})$. Going
back to (\ref{i=1...6}), we have
\begin{eqnarray*}
\mathbf{Q}\Biggl\{ \bigcap_{i=1}^3 (E^{\widetilde{S}}_{n,i}
\cap E^{\widetilde{\xi}}_{n,i})\Biggr\}
&\ge& n^{-(3/2)+o(1)} \mathbf{Q}(E^{\widetilde{S}}_{n,1} \cap
E^{\widetilde{\xi}}_{n,1})
\\
&\ge& n^{-(3/2)+o(1)}
[\mathbf{Q}(E^{\widetilde{S}}_{n,1}) +
\mathbf{Q}(E^{\widetilde{\xi}}_{n,1})-1] .
\end{eqnarray*}
We choose the constant $c_{39}>0$ sufficiently small so that
$\mathbf{Q}(E^{\widetilde{S}}_{n,1}) \ge n^{-(\varepsilon/2) + o(1)}$,
whereas $\mathbf{Q}(E^{\widetilde{\xi}}_{n,1}) = \mathbf
{Q}(E^{\widetilde{\xi
}}_{n,3}) \ge1- c_{43}  n^{-\varepsilon}$. Accordingly,
\[
\mathbf{Q}\Biggl\{ \bigcap_{i=1}^3 (E^{\widetilde{S}}_{n,i} \cap
E^{\widetilde{\xi}}_{n,i})\Biggr\}
\ge n^{-(3+\varepsilon)/2+o(1)} , \qquad n\to\infty.
\]
Substituting this into (\ref{E(Wnbeta)>}) yields
\[
\mathbf{E} ( W_{n,\beta}^{1-s} ) \ge n^{-c_{41}\varepsilon
}
n^{(3-\varepsilon)[\beta s -(\beta-1)]/2} n^{-(3+\varepsilon)/2+o(1)}.
\]
Since $\varepsilon$ can be as small as possible, this
implies the lower bound in Theorem
\ref{t:derrida-spohn-moment}.
\end{pf*}

\section[Proof of Theorem 1.5]{Proof of Theorem \protect\ref{t:tension}}\label{s:tension}

The basic idea in the proof of Theorem \ref{t:tension} is the same as
in the proof of Theorem \ref{t:derrida-spohn-moment}. Again, we prove
the upper and lower bounds in distinct parts, for the sake of clarity.
Throughout the section, we assume (\ref{hyp2}), (\ref{hyp3}) and
(\ref{hyp}).
\begin{pf*}{Proof of Theorem \ref{t:tension}: \textit{The upper bound}}
Clearly, $n^{1/2}W_n \le\overline{Y}_{   n}$, where
\[
\overline{Y}_{   n} := \sum_{|u|=n}
\bigl( n^{1/2}\vee V(u)^+ \bigr) \ee^{-V(u)} .
\]
Recall $W_n^*$ from (\ref{Mn*}). Applying (\ref{Mn*Mn}) to
$\lambda=1$, we see that $\overline{Y}_{   n} \ge{1\over c_{44}}
\log({1\over W_n^*})$, with $c_{44} := c_{12} + c_{13}$. Thus,
$\mathbf{P}\{\overline{Y}_{   n} <x,   \mathscr{S}_n\} \le\mathbf{P}\{ \log
({1\over W_n^*}) < c_{44}  x,   \mathscr{S}_n\}
\le\ee^{c_{44}}
\mathbf{E}\{(W_n^*)^{1/x} {\bf1}_{\mathscr{S}_n}\}$,
which, according to (\ref{Mn*-moment}), is bounded by $\ee
^{c_{44}} (x^\kappa+ \ee^{-c_{10} n})$ for $0<x\le{1\over a_0}$.
Thus, for any fixed $c>0$ and $0<s< \min\{ {c_{10}\over c},   \kappa
\}$, we have $\sup_{n\ge1} \mathbf{E}\{ {1\over\overline{Y}_{n}^s}
{\bf1}_{ \{ \overline{Y}_{   n} \ge\ee^{-cn} \} \cap
\mathscr{S}_n} \}
<\infty$. On the other hand, let $c_{31}$ and $c_{32}$ be as in
(\ref{E(infV):3e}); since $\overline{Y}_{   n} \ge\exp\{-\inf_{|u|=n}V(u)\}$,
it follows from (\ref{E(infV):3e}) that $\sup_{n\ge1}
\mathbf{E}\{{1\over\overline{Y}_{   n}^{c_{32}}}
{\bf1}_{ \{ \overline{Y}_{   n} < \ee^{-c_{31} n} \} \cap\mathscr{S}_n} \}<\infty$.
As a consequence,
%
\begin{equation}\label{E(1/Y)}
\sup_{n\ge1} \mathbf{E} \biggl\{ {1\over
\overline{Y}_{   n}^s}   {\bf1}_{\mathscr{S}_n}
 \biggr\} <\infty, \qquad
0<s< \min\biggl\{ c_{32},   {c_{10}\over c_{31}},  \kappa\biggr\} .
\end{equation}

We now fix $0<s<\min\{ {1\over2},   c_{32},   {c_{10}\over c_{31}},
  \kappa\}$. Let $K\ge1$ and let $E_n$ be the event in~(\ref{En}),
satisfying $\mathbf{Q}\{ E_n \} \ge1- n^{-K}$ for $n\ge n_0$. We write
\[
\mathbf{E}\{ (n^{1/2} W_n)^{1-s} \} = \mathbf{E}\{ (n^{1/2}
W_n)^{1-s}   {\bf1}_{E_n} \} + \mathbf{E}\{ (n^{1/2} W_n)^{1-s}   {\bf1}_{E_n^c} \} .
\]
For $n\ge n_0$, $\mathbf{E}\{ W_n^{1-s}   {\bf1}_{E_n^c}\} \le
[\mathbf{E}\{ W_n^{1-2s} \}]^{1/2} [\mathbf{E}\{W_n   {\bf
1}_{E_n^c}\}]^{1/2} = [\mathbf{E}
\{ W_n^{1-2s} \}]^{1/2} \times\break [\mathbf{Q}\{E_n^c\}]^{1/2}
\le[\mathbf{E}\{W_n \}]^{(1/2)-s} n^{-K/2}$, which equals $n^{-K/2}$ (since\vadjust{\goodbreak}
$\mathbf{E}\{W_n\}=1$). Therefore, for $n\to\infty$,
\[
\mathbf{E}\{ (n^{1/2} W_n)^{1-s} \} \le\mathbf{E}
\{ \overline{Y}_{  n}^{1-s}{\bf1}_{E_n} \} + o(1).
\]
Exactly as in (\ref{E(WnZ)}), we have $\mathbf{E}\{
\overline{Y}_{   n}^{1-s}   {\bf1}_{E_n} \}
= \mathbf{E}_\mathbf{Q}\{ (n^{1/2}\vee V({w_n^{(n)}})^+)
\overline{Y}_{   n}^{-s}   {\bf1}_{E_n}\}$. Thus, for $n\to\infty$,
%
\begin{equation}\label{tarte}
\mathbf{E}\{ (n^{1/2} W_n)^{1-s} \}
\le\mathbf{E}_\mathbf{Q}\bigl\{ \bigl( n^{1/2}+ V\bigl({w_n^{(n)}}\bigr)^+ \bigr)
\overline{Y}_{   n}^{-s}   {\bf1}_{E_n} \bigr\}+ o(1).
\end{equation}

For any subset $\mathscr{L}\subset\{ 1,2, \ldots, n\}$, we have
\begin{eqnarray*}
\overline{Y}_{   n} &\ge& \sum_{j\in\mathscr{L}}
\sum_{u\in\mathscr{I}^{(n)}_j}
\sum_{x\in\mathbb{T}^{\mathrm{GW}}_u,   |x|_u=n-j}
\max \{ n^{1/2},   V(x)^+ \} \ee^{-V(x)}
\\
&=& \sum_{j\in\mathscr{L}}
\sum_{u\in\mathscr{I}^{(n)}_j} \ee^{-V(u)}
\sum_{x\in\mathbb{T}^{\mathrm{GW}}_u,   |x|_u=n-j}
\max \{ n^{1/2}, [V(u) + V_u(x)]^+ \} \ee^{-V_u(x)} .
\end{eqnarray*}
Recall that $\underline{w}^{(n)}$ is the oldest vertex in
$\llbracket e,   w_n^{(n)}\rrbracket $ such that
$V({\underline{w}^{(n)}}) =\break \min_{u\in\llbracket e,   w_n^{(n)}\rrbracket } V(u)$.
Let $c_{35}$ be the constant in (\ref{E1=}). We choose
\[
\mathscr{L} := \cases{%
\bigl\{ j\le n\dvtx   \mathscr{I}^{(n)}_j \not=\varnothing,
  \bigl|{\underline{w}^{(n)}}\bigr| <j < \bigl|{\underline{w}^{(n)}}\bigr|
  + c_{35} \log n\bigr\}, \cr
\qquad \mbox{if $n-\bigl|{\underline{w}^{(n)}}\bigr|   \ge2c_{35} \log n$}, \vspace*{3pt}\cr
\bigl\{ j\le n\dvtx   \mathscr{I}^{(n)}_j \not=\varnothing,
\bigl |{\underline{w}^{(n)}}\bigr| - c_{35} \log n < j
< \bigl|{\underline{w}^{(n)}}\bigr|   \bigr\}, \cr
\qquad \mbox{otherwise}.}
\]

On the event $E_n$, it is clear that $\mathscr{L} \not= \varnothing$
and that, for any $u\in\mathscr{I}_j^{(n)}$ (with $j\in\mathscr{L}$),
%
\begin{equation}\label{V(u)-V(w)}
\bigl|V(u) - V\bigl(\underline{w}^{(n)}\bigr)\bigr|
\le c_{45} \log n,
\end{equation}
where $c_{45} := c_{36} + c_{37}$, with $c_{36}$ and $c_{37}$
as in (\ref{E2=}) and (\ref{E3=}), respectively.

We distinguish two possible situations, depending on whether
$V(\underline{w}^{(n)}) \ge- c_{46} \log n$, where
$c_{46} := {1\over s} + c_{45}$. In both situations,
we consider a sufficiently large $n$
and an arbitrary $u\in\mathscr{I}_j^{(n)}$ (with $j\in\mathscr{L}$).

On $\{ V({\underline{w}^{(n)}}) \ge- c_{46} \log n \} \cap E_n$, we
have $\max\{ n^{1/2},   [V(u) + V_u(x)]^+ \} \ge{1\over2} (n^{1/2}
\vee V_u(x)^+)$ [this holds trivially in case $V_u(x) \le n^{1/2}$;
otherwise $[V(u) + V_u(x)]^+ \ge V_u(x) -(c_{46} + c_{45})\log n \ge
{1\over2} V_u(x)^+$] and, thus,
\begin{eqnarray*}
\overline{Y}_{   n} &\ge&{1\over2} \sum_{j\in\mathscr{L}}
\sum_{u\in\mathscr{I}^{(n)}_j} \ee^{-V(u)}
\sum_{x\in\mathbb{T}^{\mathrm{GW}}_u,   |x|_u=n-j}
\max \{ (n-j)^{1/2}, V_u(x)^+ \} \ee^{-V_u(x)}
\\
&=&\!: {1\over2} \sum_{j\in\mathscr{L}}
\sum_{u\in\mathscr{I}^{(n)}_j} \ee^{-V(u)} \xi_u.
\end{eqnarray*}
If, however, $V({\underline{w}^{(n)}}) < - c_{46} \log n$,
then on $E_n$, $V(u) \le V({\underline{w}^{(n)}}) + c_{45} \log n < -
{1\over s} \log n$ and, since $\max\{ n^{1/2},   [V(u) + V_u(x)]^+ \}
\ge n^{1/2}$, we have, in this case,
\begin{eqnarray*}
\overline{Y}_{   n}
&\ge& n^{(1/s)+(1/2)} \sum_{j\in\mathscr{L}}
\sum_{u\in\mathscr{I}^{(n)}_j}
\sum_{x\in\mathbb{T}^{\mathrm{GW}}_u,  |x|_u=n-j}
\ee^{-V_u(x)}
\\
&=&\!: n^{(1/s)+(1/2)}
\sum_{j\in \mathscr{L}} \sum_{u\in\mathscr{I}^{(n)}_j} \eta_u.
\end{eqnarray*}
Therefore, in both situations we have
%
\begin{eqnarray}\label{Yn>}
\overline{Y}_{   n}^{-s}   {\bf1}_{E_n}
&\le& 2^s  \Biggl( \sum_{j\in\mathscr{L}}
\sum_{u\in\mathscr{I}^{(n)}_j} \ee^{-V(u)}
\xi_u \Biggr)^{  -s} {\bf1}_{E_n}
\nonumber\\[-8pt]
\\[-8pt]
&&{} + n^{-(s/2)-1} \Biggl( \sum_{j\in\mathscr{L}}
\sum_{u\in\mathscr{I}^{(n)}_j} \eta_u\Biggr)^{  -s} {\bf1}_{E_n} .
\nonumber
\end{eqnarray}
[Since $\sum_{j\in\mathscr{L}} \sum_{u\in\mathscr{I}_j^{(n)}}
\sum_{x\in\mathbb{T}^{\mathrm{GW}}_u,|x|_u=n-j} 1 >0$
on $E_n$, the $(\cdot)^{-s}$ expressions on the right-hand side are
well defined.]

We claim that there exists $0<s_0<1$ such that,
for any $\varepsilon>0$ and $s\in(0,   s_0)$,
%
\begin{eqnarray}
&& \mathbf{E}_\mathbf{Q} \Biggl\{ \bigl( n^{1/2}+ V\bigl(w_n^{(n)}\bigr)^+ \bigr)
\Biggl ( \sum_{j\in\mathscr{L}}
\sum_{u\in\mathscr{I}^{(n)}_j} \ee^{-V(u)}
\xi_u \Biggr)^{-s}   {\bf1}_{E_n}  \Biggr\}
\nonumber\\[-8pt]\label{key-estimate-1}
\\[-8pt]
&&\qquad \le c_{48},
\nonumber \\
&& \mathbf{E}_\mathbf{Q} \Biggl\{ \bigl( n^{1/2}+ V\bigl({w_n^{(n)}}\bigr)^+ \bigr)
 \Biggl( \sum_{j\in\mathscr{L}}
\sum_{u\in\mathscr{I}^{(n)}_j} \eta_u
\Biggr)^{-s}   {\bf1}_{E_n}  \Biggr\}
\nonumber\\[-8pt]\label{key-estimate-2}
\\[-8pt]
&&\qquad \le c_{47} n^{{1/2} +{(3+\varepsilon)/2} s}.
\nonumber
\end{eqnarray}

We admit (\ref{key-estimate-1}) and (\ref{key-estimate-2}) for the
time being. In view of (\ref{Yn>}), we obtain, for
$0<s< s_*:= \min\{{1\over2},   s_0,   c_{32},
{c_{10}\over c_{31}},   \kappa\}$,
\[
\mathbf{E}_\mathbf{Q} \bigl\{ \bigl( n^{1/2}+ V\bigl({w_n^{(n)}}\bigr)^+ \bigr)
\overline{Y}_{   n}^{-s}   {\bf1}_{E_n}  \bigr\} \le2^s c_{48} + o(1).
\]
Substituting this in (\ref{tarte}), we see that $\sup_{n\ge
1} \mathbf{E}\{ (n^{1/2} W_n)^{1-s} \}<\infty$ for any $s\in(0,   s_*)$.
This yields the last inequality in (\ref{tension}) when $\gamma$ is
close to 1. By Jensen's inequality, it holds for all
$\gamma\in[0,1)$. This will complete the proof of the upper bound in Theorem
\ref{t:tension}.

It remains to check (\ref{key-estimate-1}) and
(\ref{key-estimate-2}). We only present the proof of
(\ref{key-estimate-1}), because the proof of (\ref{key-estimate-2}) is
similar and slightly easier, using (\ref{Mn-moment-neg}) in place of
(\ref{E(1/Y)}).

Recall $\mathscr{G}_n$ from (\ref{Gn}). By Proposition
\ref{p:change-proba}, under $\mathbf{Q}$ and conditionally on $\mathscr{G}_n$,
the random variables $\xi_u$, for $u\in\mathscr{I}_j^{(n)}$ and
$j\in\mathscr{L}$, are independent. We write $\mathscr{L}:= \{ j(1),\ldots, j(N)\}$,
with $j(1)< \cdots<j(N)$. It follows from the second
part of Lemma \ref{l:max-ineq} that
\begin{eqnarray*}
&& \mathbf{E}_\mathbf{Q} \Biggl\{
 \Biggl( \sum_{j\in\mathscr{L}}
\sum_{u\in\mathscr{I}^{(n)}_j} \ee^{-V(u)}
\xi_u \Biggr)^{-s}   {\bf1}_{E_n} \Big|   \mathscr{G}_n  \Biggr\}
\\
&& \qquad \le \sum_{i=1}^N b^{i-1}
\mathbf{E}_\mathbf{Q} \Biggl\{  \Biggl(
\sum_{u\in\mathscr{I}^{(n)}_{j(i)}}
\ee^{-V(u)}\xi_u  \Biggr)^{-s}
{\bf1}_{\{\sum_{u\in\mathscr{I}^{(n)}_{j(i)}}
\ee^{-V(u)}\xi_u >0\} }\Big |   \mathscr{G}_n  \Biggr\},
\end{eqnarray*}
where $b:= \max_{j\in\mathscr{L}} \mathbf{Q}
\{ \sum_{u\in\mathscr{I}^{(n)}_j} \ee^{-V(u)} \xi_u =0  |   \mathscr{G}_n\}$.
We note that $b \le\break \max_{1\le j\le n} \mathbf{P}\{ \mathscr{S}_{n-j}^c\}\le q$,
and that, for any $i\le N$, the $\mathbf{E}_\mathbf{Q}\{\cdot\}$ expression
on the right-hand side is, according to the first part of
Lemma \ref{l:max-ineq}, bounded by
\[
{1\over1-q} \max_{u\in\mathscr{I}^{(n)}_{j(i)}}
\mathbf{E}_\mathbf{Q}
\biggl\{ {\ee^{s V(u)}\over\xi_u^s}
{\bf1}_{ \{\xi_u>0\} } \big|   \mathscr{G}_n  \biggr\}.
\]
By Proposition \ref{p:change-proba},
$\mathbf{E}_\mathbf{Q}\{ {1\over\xi_u^s}   {\bf1}_{ \{\xi_u>0\} }  |
\mathscr{G}_n \} = \mathbf{E}\{{1\over\overline{Y}_{   n-j}^s}
{\bf1}_{ \mathscr{S}_{n-j} } \}$, which is
bounded in $n$ and $j$ [by (\ref{E(1/Y)})]. Summarizing,
we have proved that
\[
\mathbf{E}_\mathbf{Q} \Biggl\{  \Biggl( \sum_{j\in\mathscr{L}}
\sum_{u\in\mathscr{I}^{(n)}_j}
\ee^{-V(u)} \xi_u \Biggr)^{-s}
{\bf1}_{E_n}    \Big|   \mathscr{G}_n \Biggr\}
\le c_{49} \sum_{i=1}^N q^{i-1}
\max_{u\in\mathscr{I}^{(n)}_{j(i)}} \ee^{sV(u)} .
\]
As a consequence, the expression on the left-hand side of
(\ref{key-estimate-1}) is bounded by
$c_{49}   \mathbf{E}_\mathbf{Q}\{ \Lambda_n\}$, where
\begin{eqnarray*}
\Lambda_n :\!\!&=& \bigl(n^{1/2} + V\bigl(w^{(n)}_n\bigr)^+\bigr)
\sum_{i=1}^N q^{i-1}
\max_{u\in\mathscr{I}^{(n)}_{j(i)} }
\ee^{sV(u)}
{\bf1}_{ \{ |V(u) - V(\underline{w}^{(n)})|\le c_{45} \log n\} }
\\
&\le& \widetilde{\Lambda}_n := \bigl(n^{1/2} + V\bigl(w^{(n)}_n\bigr)^+\bigr)
\sum_{i=1}^N q^{i-1} \max_{u\in\mathscr{I}^{(n)}_{j(i)} }
\ee^{sV(u)} .
\end{eqnarray*}
The proof of (\ref{key-estimate-1}) now boils down to
verifying the following estimates: there exists $0<s_0<1$ such that,
for any $s\in(0,   s_0)$,
%
\begin{eqnarray}
\sup_n   \mathbf{E}_\mathbf{Q}\bigl\{ \widetilde{\Lambda}_n
{\bf1}_{ \{ n-|{\underline{w}^{(n)}}| \ge
2 c_{35} \log n \} } \bigr\} &<& \infty  ,
\label{key-estimate-1bis}\\
\lim_{n\to\infty}
\mathbf{E}_\mathbf{Q}\bigl\{ \Lambda_n
{\bf1}_{ \{ n-|{\underline{w}^{(n)}}| <
2 c_{35} \log n \} } \bigr\} &=& 0.
\label{key-estimate-1ter}
\end{eqnarray}

Let us first check (\ref{key-estimate-1bis}). Let $S_0:=0$ and let
$(S_j-S_{j-1},   \sigma_j,   \Delta_j)$, $j\ge1$, be i.i.d.
random variables under $\mathbf{Q}$ and distributed as $(V(w^{(1)}),
 \# \mathscr{I}_1^{(1)},\break  \max_{u\in\mathscr{I}_1^{(1)}} \ee^{sV(u)})$. Let
\[
\underline{S}_n := \min_{0\le i\le n} S_i, \qquad
\vartheta_n :=\inf\{ k\ge0\dvtx   S_k = \underline{S}_n\} .
\]
[The random variable $\vartheta_n$ has nothing to do with
the constant $\vartheta$ in Proposition \ref{p:tail-Mn}.] Writing
LHS\tsub{(\ref{key-estimate-1bis})} for
$\mathbf{E}_\mathbf{Q}\{\widetilde{\Lambda}_n
{\bf1}_{ \{ n-|{\underline{w}^{(n)}}|\ge 2 c_{35} \log n \} }\}$,
it follows from Proposition \ref{p:change-proba} that
\begin{eqnarray*}
\mathrm{LHS}_{(\ref{key-estimate-1bis})}
&=& \mathbf{E}_\mathbf{Q}
\Biggl\{[n^{1/2} + S_n^+] \sum_{i=1}^M q^{i-1}
\ee^{sS_{\ell(i)-1}}\Delta_{\ell(i)}
{\bf1}_{ \{ n-\vartheta_n \ge 2 c_{35} \log n \} }  \Biggr\}
\\
&=& \mathbf{E}_\mathbf{Q} \Biggl\{
[n^{1/2} + S_n^+] \ee^{s\underline{S}_n}
\sum_{i=1}^M q^{i-1}\ee^{s[S_{{\ell(i)}-1}-S_{\ell(0)}]}
\Delta_{\ell(i)}{\bf1}_{ \{ n-\vartheta_n \ge
2 c_{35} \log n \} }  \Biggr\} ,
\end{eqnarray*}
where $\ell(i) := \inf\{ k> \ell(i-1)\dvtx   \sigma_k \ge1\}
$ with $\ell(0):= \vartheta_n$, and $M := \sup\{ i\dvtx
\ell(i) <\vartheta_n + c_{35} \log n\}$.

At this stage, we use a standard trick for random walks: let
$\nu_0:=0$ and let
\[
\nu_i := \inf \biggl\{ k>\nu_{i-1}\dvtx
S_k < \min_{0\le j\le\nu_{i-1}} S_j  \biggr\}, \qquad i\ge1.
\]
In words, $0=\nu_0<\nu_1< \cdots$ are strict descending
ladder times. On the event $\{\nu_k \le n<\nu_{k+1} \}$ (for $k\ge0$),
we have $\vartheta_n =\nu_k$ and $\underline{S}_n = S_{\nu_k}$.
Thus, $\mathrm{LHS}$\tsub{(\ref{key-estimate-1bis})} equals
\begin{eqnarray*}
&& \sum_{k=0}^\infty\mathbf{E}_\mathbf{Q}
\Biggl\{{\bf1}_{ \{ n-\nu_k \ge2 c_{35} \log n \} }
{\bf1}_{ \{\nu_k \le n<\nu_{k+1} \} }
[n^{1/2} + S_n^+] \ee^{sS_{\nu_k}}
\\
&&\hspace*{104pt}{}\times
\sum_{i=1}^M q^{i-1}\ee^{s[S_{\ell(i)-1}-S_{\ell(0)}]}
\Delta_{\ell(i)}  \Biggr\}.
\end{eqnarray*}
For any $k$, we look at the expectation $\mathbf{E}_\mathbf{Q}\{ \cdot\}$
on the right-hand side. By conditioning upon $(S_j,   \sigma_j,
\Delta_j,   1\le j\le\nu_k)$, and since $S_n^+ = [S_{\nu_k} +
(S_n-S_{\nu_k})]^+ \le(S_n-S_{\nu_k})^+ =S_n-S_{\nu_k}$ on
$\{\nu_k \le n<\nu_{k+1} \}$, we obtain
%
\begin{equation}\label{fn}
\mathrm{LHS}_{(\ref{key-estimate-1bis})} \le
\sum_{k=0}^\infty\mathbf{E}_\mathbf{Q}
\bigl\{{\bf1}_{ \{ n-\nu_k \ge2 c_{35} \log n \} }
  \ee^{sS_{\nu_k}} f_n(n-\nu_k)  \bigr\} ,
\end{equation}
where, for any $1\le j\le n$,
\[
f_n(j) := \mathbf{E}_\mathbf{Q} \Biggl\{ {\bf1}_{\{ \nu_1 > j\} }
[n^{1/2} + S_j]
\sum_{i=1}^{M'} q^{i-1} \ee^{s S_{m(i)-1}}
\Delta_{m(i)}\Biggr\},
\]
and $m(i) := \inf\{ k> m(i-1)\dvtx   \sigma_k \ge1\}$ with
$m(0):=0$, and $M' := \sup\{ i\dvtx m(i) < c_{35} \log n\}$. For
brevity, we write $L_n := \sum_{i=1}^{M'} q^{i-1} \ee^{s S_{m(i)-1}}
\Delta_{m(i)} = \sum_{i=1}^\infty q^{i-1} \times\break \ee^{s S_{m(i)-1}}
\Delta_{m(i)}   {\bf1}_{\{ m(i) < c_{35} \log n\} }$ for the moment. By the
Cauchy--Schwarz inequality,
\[
f_n(j) \le[\mathbf{Q}\{ \nu_1>j\}]^{1/2}
\bigl[\mathbf{E}_\mathbf{Q}\{ (n^{1/2} + S_j)^2   |   \nu_1>j \}\bigr]^{1/2}
\bigl[\mathbf{E}_\mathbf{Q}
\bigl\{ L_n^2   {\bf1}_{\{ \nu_1>j\} } \bigr\}\bigr]^{1/2} .
\]
By (\ref{Bingham}), $\mathbf{Q}\{ \nu_1>j\} \le c_{50}  j^{-1/2}$
for some $c_{50}>0$ and all $j\ge1$. On the other hand, $(n^{1/2} +
S_j)^2 \le2(n+S_j^2)$, and it is known (Bolthausen \cite{bolthausen})
that $\mathbf{E}_\mathbf{Q}\{ {S_j^2\over j}   |   \nu_1>j \} \to
c_{51} \in(0,  \infty)$ for $j\to\infty$. Therefore,
$\mathbf{E}_\mathbf{Q}\{(n^{1/2} + S_j)^2   |   \nu_1>j \} \le c_{52}  n$
for some $c_{52}>0$ and all
$n\ge j\ge1$. Accordingly, with $c_{53} := c_{50}^{1/2} c_{52}^{1/2}$,
we have
\[
f_n(j) \le c_{53}  j^{-1/4} n^{1/2}
\bigl[\mathbf{E}_\mathbf{Q}
\bigl\{ L_n^2   {\bf1}_{\{ \nu_1>j \} } \bigr\}\bigr]^{1/2},
\qquad1\le j\le n.
\]
By the Cauchy--Schwarz inequality,
$L_n^2 \le(\sum_{i=1}^\infty q^{i-1}) \sum_{i=1}^\infty q^{i-1}
\ee^{2sS_{{m(i)}-1}} \Delta_{m(i)}^2\times\break  {\bf1}_{ \{ m(i) < c_{35} \log n\}}$.
Therefore, for $j\ge2c_{35} \log n$,
\begin{eqnarray*}
\mathbf{E}_\mathbf{Q}\bigl\{ L_n^2   {\bf1}_{\{ \nu_1>j \} }\bigr\}
&\le& {1\over1-q}
\sum_{i=1}^\infty q^{i-1}\mathbf{E}_\mathbf{Q}
\bigl\{ \ee^{2s S_{m(i)-1}}
\Delta_{m(i)}^2
{\bf1}_{ \{ m(i) < c_{35} \log n\} }
{\bf1}_{\{ \nu_1>j\} } \bigr\}
\\
&\le& {1\over1-q}   \sum_{i=1}^\infty q^{i-1}
\mathbf{E}_\mathbf{Q}
\bigl\{ \ee^{2s S_{m(i)-1}} \Delta_{m(i)}^2
{\bf1}_{ \{ m(i) \le{j/2}\} }{\bf1}_{\{ \nu_1>j\} } \bigr\} .
\end{eqnarray*}
For any $i\ge1$, to estimate the expectation
$\mathbf{E}_\mathbf{Q}\{\cdot\}$ on the right-hand side,
we apply the strong Markov property
at time $m(i)$ to see that
\[
\mathbf{E}_\mathbf{Q}\{ \cdot\}
= \mathbf{E}_\mathbf{Q}\bigl\{ \ee^{2s S_{m(i)-1}} \Delta_{m(i)}^2
{\bf1}_{ \{ m(i) \le{j/2}\} }   {\bf1}_{\{ \nu_1>m(i)\} }
g\bigl(S_{m(i)}, j-m(i)\bigr)\bigr\} ,
\]
where $g(z,k) := \mathbf{Q}\{ z+S_i \ge0,   \forall1\le i\le k\}$
for any $z\ge0$ and $k\ge1$. By (13) of Kozlov \cite{kozlov},
$g(z,k) \le c_{54} (z+1)/k^{1/2}$ for some $c_{54}>0$ and all $z\ge0$
and $k\ge1$. Since $z+1 \le c_{55}\ee^{sz}$ for all $z\ge0$, this
yields, with $c_{56} := {c_{55} \over1-q}$,
\begin{eqnarray*}
\mathbf{E}_\mathbf{Q}\bigl\{ L_n^2   {\bf1}_{\{ \nu_1>j \} } \bigr\}
&\le& c_{56}   \sum_{i=1}^\infty
q^{i-1} \mathbf{E}_\mathbf{Q}
\biggl\{ \ee^{2s S_{m(i)-1}} \Delta_{m(i)}^2
{\bf1}_{ \{ m(i) \le{j/2}\} }
{\ee^{s S_{m(i)}}\over(j-m(i))^{1/2}}  \biggr\}
\\
&\le& {c_{56}\over(j/2)^{1/2}}
\sum_{i=1}^\infty q^{i-1} \mathbf{E}_\mathbf{Q}
\bigl\{ \ee^{2s S_{m(i)-1} + s S_{m(i)}}
\Delta_{m(i)}^2 \bigr\}
\\
&=& {2^{1/2} c_{56}\over j^{1/2}}
\mathbf{E}_\mathbf{Q} \Biggl\{ \sum_{i=1}^\infty
q^{i-1} \ee^{2s S_{m(i)-1} + s S_{m(i)}}
\Delta_{m(i)}^2  \Biggr\} .
\end{eqnarray*}
We observe that $\sum_{i=1}^\infty q^{i-1} \ee^{2s
S_{m(i)-1} + s S_{m(i)}} \Delta_{m(i)}^2 \le\sum_{k=1}^\infty
q^{R(k)-1} \ee^{2s S_{k-1} + s S_k} \Delta_k^2$, where
$R(k):= \# \{1\le j\le k\dvtx \sigma_j \ge1\}$. Therefore, with
$c_{57} :=2^{1/2} c_{56}$,
\begin{eqnarray*}
\mathbf{E}_\mathbf{Q}\bigl\{ L_n^2   {\bf1}_{\{ \nu_1>j \} } \bigr\}
&\le& {c_{57}\over j^{1/2}} \sum_{k=1}^\infty
\mathbf{E}_\mathbf{Q}
\bigl\{ q^{R(k)-1} \ee^{2s S_{k-1} + s S_k}\Delta_k^2\bigr\}
\\
&\le& {c_{57}\over j^{1/2}} \sum_{k=1}^\infty
\bigl[\mathbf{E}_\mathbf{Q}\bigl\{ q^{2[R(k)-1]} \bigr\}\bigr]^{1/2}
[\mathbf{E}_\mathbf{Q}\{\ee^{4s S_{k-1} + 2s S_k}
\Delta_k^4\}]^{1/2}.
\end{eqnarray*}
By definition, $\mathbf{E}_\mathbf{Q}\{ q^{2[R(k)-1]} \} =q^{-2}   r^k$,
with $r:= \mathbf{Q}(\sigma_1=0) + q^2 \mathbf{Q}(\sigma_1\ge1)<1$
[because $q<1$ and $\mathbf{Q}(\sigma_1=0)<1$]. On the other hand,
\begin{eqnarray*}
\mathbf{E}_\mathbf{Q}\{\ee^{4s S_{k-1} + 2s S_k}
\Delta_k^4\} &=& \mathbf{E}_\mathbf{Q}\{\ee^{6s S_{k-1}} \}
\mathbf{E}_\mathbf{Q}\bigl\{ \ee^{2s (S_k-S_{k-1})}\Delta_k^4\bigr\}
\\
&=& [\mathbf{E}_\mathbf{Q}\{\ee^{6s S_1}\}]^{k-1}
\mathbf{E}_\mathbf{Q}\{ \ee^{2s S_1} \Delta_1^4\} .
\end{eqnarray*}
By (\ref{exp-moment-S1}), there exists $s_{\#}>0$
sufficiently small such that
$\mathbf{E}_\mathbf{Q}\{\ee^{6s S_1}\} < {1\over r}$
for all $0<s<s_{\#}$. On the other hand, $\mathbf{E}_\mathbf{Q}\{\Delta_1^8\}
<\infty$ for $0<s<{c_6\over8}$ [by (\ref{tail-sup|V|})], and
$\mathbf{E}_\mathbf{Q}
\{ \ee^{4s S_1} \} <\infty$ for $0<s\le{c_2\over4}$ [by
(\ref{exp-moment-S1})]; thus,
$\mathbf{E}_\mathbf{Q}\{ \ee^{2s S_1}\Delta_1^4\} <\infty
$ for $0<s<\min\{ {c_6\over8},   {c_2\over4}\}$. As a consequence,
for any $0<s<\min\{ s_{\#},   {c_6\over8},   {c_2\over4}\}$, we
have $\mathbf{E}_\mathbf{Q}\{ L_n^2   {\bf1}_{\{ \nu_1>j \} } \}
\le{c_{58}\over j^{1/2}}$, for some $c_{58}>0$ and all
$n\ge j\ge1$ with $j\ge2c_{35}\log n$, which yields
\[
f_n(j) \le c_{53}c_{58}^{1/2} j^{-1/2}n^{1/2}.
\]
Going back to (\ref{fn}), we obtain, for any
$0<s<\min\{s_{\#},   {c_6\over8},   {c_2\over4}\}$
and $c_{59} := c_{53}c_{58}^{1/2}$,
\[
\mathrm{LHS}_{(\ref{key-estimate-1bis})} \le
c_{59}   n^{1/2} \sum_{k=0}^\infty\mathbf{E}_\mathbf{Q}
\biggl\{{\bf1}_{ \{ n-\nu_k \ge2 c_{35} \log n \} }
{\ee^{sS_{\nu_k}} \over(n-\nu_k)^{1/2}} \biggr\}.
\]
By (\ref{Bingham}) again, ${1\over j^{1/2}} \le c_{60}\mathbf{Q}
\{ \nu_1 >j\}$ for all $j\ge1$. Thus, with $c_{61} := c_{59}c_{60}$,
\begin{eqnarray*}
\mathrm{LHS}_{(\ref{key-estimate-1bis})}
&\le& c_{61}  n^{1/2} \sum_{k=0}^\infty
\mathbf{E}_\mathbf{Q} \bigl\{
{\bf1}_{ \{ n-\nu_k \ge2 c_{35} \log n \} }
\ee^{sS_{\nu_k}}  {\bf1}_{\{ \nu_{k+1}>n\} } \bigr\}
\\
&\le& c_{61}   n^{1/2} \sum_{k=0}^\infty
\mathbf{E}_\mathbf{Q} \bigl\{
{\bf1}_{ \{ \nu_k\le n<\nu_{k+1}\} }
\ee^{sS_{\nu_k}}  \bigr\} ,
\end{eqnarray*}
which equals $c_{61}   n^{1/2} \mathbf{E}_\mathbf{Q}
\{\ee^{s \min_{0\le i\le n}S_i}\}$, and, according
to (\ref{E[exp(min)]}), is bounded in
$n$. This completes the proof of (\ref{key-estimate-1bis}).

It remains to check (\ref{key-estimate-1ter}). By definition,
\[
\Lambda_n \le\bigl[n^{1/2} + V\bigl(w^{(n)}_n\bigr)^+\bigr]
n^{s c_{45}} \ee^{sV(\underline{w}^{(n)})} \sum_{i=1}^N q^{i-1} .
\]
Since $\sum_{i=1}^N q^{i-1} \le{1\over1-q}$, this leads
to, by an application of Proposition \ref{p:change-proba},
\[
\mathbf{E}_\mathbf{Q}\bigl\{ \Lambda_n   {\bf1}_{ \{ n-|{\underline
{w}^{(n)}}| < 2 c_{35} \log n \} } \bigr\}
\le{n^{s c_{45}} \over1-q} \mathbf{E}_\mathbf{Q}
\bigl\{ [n^{1/2} +S_n^+] \ee ^{s\underline{S}_n}
{\bf1}_{\{ n-\vartheta_n < 2 c_{35} \log n \} }\bigr\},
\]
where $(S_i)$ is as in Proposition \ref{p:change-proba} and,
as before, $\underline{S}_n := \min_{0\le i\le n} S_i$,
$\vartheta_n:= \inf\{ k\ge0\dvtx S_k = \underline{S}_n\}$.

Let $0<\varepsilon<{1\over2}$; let $A_n:= \{ S_n >n^{{1/2} +
\varepsilon} \}$ and $B_n:= \{ S_n \le n^{{1/2} + \varepsilon}\}
= A_n^c$.

Since $\mathbf{E}_\mathbf{Q}\{ \ee^{aS_1}\} <\infty$ for
$|a|<c_2$ [see (\ref{exp-moment-S1})] and
$\mathbf{Q}(A_n) \le2\exp(-c_3 n^{2\varepsilon})$
[see (\ref{Petrov})], the Cauchy--Schwarz inequality yields $n^{s
c_{45}} \mathbf{E}_\mathbf{Q}\{ [n^{1/2} + S_n^+]\*
\ee^{s\underline{S}_n}   {\bf1}_{A_n}\} \to0$, $n\to\infty$.

On $B_n$, we have $n^{1/2} + S_n^+ \le2n^{{1/2} + \varepsilon}$;
thus, $\mathbf{E}_\mathbf{Q}\{ [n^{1/2} + S_n^+]
\ee^{s\underline{S}_n} \times\break {\bf1}_{B_n
\cap\{ n-\vartheta_n < 2 c_{35} \log n \}} \} \le2n^{{1/2} +
\varepsilon} \mathbf{E}_\mathbf{Q}\{ \ee^{s\underline{S}_n}
{\bf1}_{\{n-\vartheta_n < 2 c_{35} \log n \}} \}$. It is clear that
$\underline{S}_n \le\underline{S}_{\lfloor n/2\rfloor}:= \min_{0\le i\le
n/2}S_i$, and that $\{ n-\vartheta_n < 2 c_{35} \log n \} \subset\{
{n\over2} -\widetilde{\vartheta}_{n/2} < 2 c_{35} \log n \}$, where
$\widetilde{\vartheta}_{n/2} := \min\{ k\ge0\dvtx  \widetilde{S}_k =
\min_{0\le i\le n-\lfloor n/2\rfloor} \widetilde{S}_i\}$, with
$\widetilde{S}_i := S_{i+\lfloor n/2\rfloor}-S_{\lfloor n/2\rfloor}$,
$i\ge0$. Since $\underline{S}_{\lfloor n/2\rfloor}$ and
$\widetilde{\vartheta}_{n/2}$ are independent, we have
$\mathbf{E}_\mathbf{Q}\{ \ee^{s\underline{S}_n} {\bf1}_{\{ n-\vartheta_n < 2
c_{35} \log n \}}\} \le\mathbf{E}_\mathbf{Q}\{
\ee^{s\underline{S}_{\lfloor n/2\rfloor}} \}   \mathbf{Q}
\{ {n\over2} -\widetilde{\vartheta}_{n/2} < 2 c_{35} \log n \}$. By
(\ref{E[exp(min)]}), $\mathbf{E}_\mathbf{Q}\{ \ee^{s
\underline{S}_{\lfloor n/2\rfloor}} \} \le c_{62}   n^{-1/2}$; on the
other hand, $\mathbf{Q}\{{n\over2} -\widetilde{\vartheta}_{n/2} < 2 c_{35} \log n \} \le
c_{63}   {(\log n)^{1/2}\over n^{1/2}}$ (see Feller \cite{feller},
page~398). Therefore, $\mathbf{E}_\mathbf{Q}\{ [n^{1/2} + S_n^+]
\ee^{s\underline{S}_n}\* {\bf1}_{B_n \cap\{ n-\vartheta_n
< 2 c_{35} \log n \}} \} \le c_{64}   n^{-{1/2} + \varepsilon} (\log n)^{1/2}$.

Summarizing, we have proved that, for any $s>0$ and
$0<\varepsilon<{1\over2}$, when $n\to\infty$,
\[
\mathbf{E}_\mathbf{Q}\bigl\{ \Lambda_n
{\bf1}_{ \{ n-|{\underline{w}^{(n)}}| < 2 c_{35} \log n \} } \bigr\}
\le o(1) + {c_{64}\over1-q}   n^{s c_{45} -{1/2} + \varepsilon}
(\log n)^{1/2} ,
\]
which yields (\ref{key-estimate-1ter}), as long as
$0<s<{1\over2c_{45}}$.
\end{pf*}
\begin{pf*}{Proof of Theorem \ref{t:tension}}
\textit{The lower bound.} We start with
\[
n^{1/2} W_n \ge\underline{Y  }_n := \sum_{|u|=n}
\bigl(n^{1/2}\wedge V(u)^+ \bigr) \ee^{-V(u)}.
\]
Let $s\in(0,   1)$. Exactly as in (\ref{E(WnZ)}), we have
%
\begin{equation}\label{tartee}
\mathbf{E}\{ \underline{Y}_n^{  1-s} \} =
\mathbf{E}_\mathbf{Q} \bigl\{ \bigl( n^{1/2}\wedge
V\bigl({w_n^{(n)}}\bigr)^+\bigr) \underline{Y  }_n^{ -s}  \bigr\}.
\end{equation}
By definition,
\begin{eqnarray*}
\underline{Y  }_n &=& \sum_{j=1}^n
\sum_{u\in\mathscr{I}^{(n)}_j} \ee^{-V(u)}
\sum_{x\in\mathbb{T}^{\mathrm{GW}}_u,   |x|_u=n-j}
\min \{ n^{1/2}, [V(u) + V_u(x)]^+ \} \ee^{-V_u(x)}
\\
&&{} + \min \{ n^{1/2},
V(w^{(n)}_n)^+ \} \ee^{-V(w^{(n)}_n)}
\\
&\le& \sum_{j=1}^n \ee^{-V(w_{j-1}^{(n)})}
\sum_{u\in\mathscr{I}^{(n)}_j} \ee^{-\Delta_u}
\sum_{x\in\mathbb{T}^{\mathrm{GW}}_u,   |x|_u=n-j}
\bigl[V\bigl(w_{j-1}^{(n)}\bigr)^+ +\Delta_u^+ + V_u(x)^+\bigr]
\\
&&\hspace*{181pt}{}\times
\ee^{-V_u(x)} + \Theta_n,
\end{eqnarray*}
where $\Delta_u := V(u) - V(w_{j-1}^{(n)})$ [for
$u\in\mathscr{I}^{(n)}_j$], and $\Theta_n := V(w^{(n)}_n)^+
\ee^{-V(w^{(n)}_n)}$.

By means of the elementary inequality $(\sum_i a_i)^{-s}
\ge(\sum_i a_i^s)^{-1}$ and $(\sum_i b_i)^s \le\sum_i b_i^s$ for nonnegative
$a_i$ and $b_i$, we obtain $\underline{Y  }_n^{  -s}
\ge{1\over Z_n}$ on $\mathscr{S}_n$, with $Z_n$ being defined as
\begin{eqnarray*}
&& \sum_j \ee^{-sV\bigl(w_{j-1}^{(n)}\bigr)}
\sum_u \ee^{-s\Delta_u}
\Biggl\{\bigl[\bigl(V(w_{j-1}^{(n)})^+\bigr)^s
+ (\Delta_u^+)^s\bigr]
\Biggl(\sum_x \ee^{-V_u(x)}\Biggr)^s
\\
&&\hspace*{178pt}{}
+ \Biggl[\sum_x V_u(x)^+ \ee^{-V_u(x)}\Biggr]^s \Biggr\}
+ \Theta_n^s ,
\end{eqnarray*}
where $\sum_j := \sum_{j=1}^n$,
$\sum_u := \sum_{u\in\mathscr{I}^{(n)}_j}$, and
$\sum_x := \sum_{x\in\mathbb{T}^{\mathrm{GW}}_u,   |x|_u=n-j}$.
We now condition upon $\mathscr{G}_n$, and
note that $V(w_j^{(n)})$ and $\mathscr{I}^{(n)}_j$ are
$\mathscr{G}_n$-measurable. By Proposition~\ref{p:change-proba},
\begin{eqnarray*}
\mathbf{E}_\mathbf{Q}\{ Z_n   |   \mathscr{G}_n \}
&=& \sum_j\ee^{-sV(w_{j-1}^{(n)})}
\sum_u \ee^{-s\Delta_u}
\bigl\{\bigl(\bigl(V\bigl(w_{j-1}^{(n)}\bigr)^+\bigr)^s + (\Delta_u^+)^s\bigr)
\\
&&\hspace*{115pt}{}\times \mathbf{E}(W_{n-j}^s)
+ \mathbf{E}(U_{n-j}^s) \bigr\} + \Theta_n^s,
\end{eqnarray*}
where, for any $k\ge0$, $U_k := \sum_{|y|=k} V(y)^+ \ee^{-V(y)}$.
By Jensen's inequality, $\mathbf{E}(W_{n-j}^s)\le[\mathbf{E}(W_{n-j})]^s=1$.
On the other hand, by (\ref{phi(0)}), $U_k \le c_{65} \log {1\over W_k^*}$
and, thus, by Lemma \ref{l:log(1/M)}, $\mathbf{E}(U_k^s) \le
c_{65}^s \mathbf{E}\{ [\log{1\over W_k^*}]^s\} \le c_{66}$.
Therefore, the
$\sum_u$ sum on the right-hand side (without $\Theta_n^s$, of course) is
\begin{eqnarray*}
&\le& \sum_u \ee^{-s\Delta_u}
\bigl\{ \bigl(V\bigl(w_{j-1}^{(n)}\bigr)^+\bigr)^s
+ (\Delta_u^+)^s + c_{67} \bigr\}
\\
&=& \bigl[V\bigl(w_{j-1}^{(n)}\bigr)^+\bigr]^s
\sum_u \ee^{-s\Delta_u}
+ \sum_u\ee^{-s\Delta_u}\{ (\Delta_u^+)^s + c_{67} \} .
\end{eqnarray*}
There exists $c_{68} = c_{68}(s)<\infty$ such that $\ee
^{-sa} \{ (a^+)^s + c_{67} \} \le c_{68}   (\ee^{-sa} +
\ee
^{-sa/2})$ for all $a\in\mathbb{R}$. As a consequence,
\begin{eqnarray*}
\mathbf{E}_\mathbf{Q}\{ Z_n   |   \mathscr{G}_n \}
&\le & c_{69}\sum_{j=1}^n \ee^{-sV(w_{j-1}^{(n)})}
\bigl\{ \bigl[V\bigl(w_{j-1}^{(n)}\bigr)^+\bigr]^s +1\bigr\}
\\
&&\hspace*{30pt}{}\times
\sum_{u\in\mathscr{I}^{(n)}_j} [\ee^{-s\Delta_u}
+ \ee^{-{s/2}\Delta_u}] + \Theta_n^s.
\end{eqnarray*}
By Jensen's inequality again, $\mathbf{E}_\mathbf{Q}\{{1\over Z_n}   |
\mathscr{G}_n\} \ge{1\over\mathbf{E}_\mathbf{Q}\{ Z_n  |
\mathscr{G}_n\} }$.
Since $\underline{Y  }_n^{  -s} \ge{1\over Z_n}$ on $\mathscr
{S}_n$, this leads to
\begin{eqnarray*}
&& \mathbf{E}_\mathbf{Q}\{ \underline{Y  }_n^{  -s}   |
\mathscr{G}_n \}
\\
&&\qquad \ge {c_{70} \over\sum_{j=1}^n \ee^{-sV(w_{j-1}^{(n)})} \{
[V(w_{j-1}^{(n)})^+]^s +1\} \sum_{u\in\mathscr{I}^{(n)}_j} [\ee
^{-s\Delta_u}+\ee^{-{s\over2}\Delta_u}] + \Theta_n^s} .
\end{eqnarray*}
We apply Proposition \ref{p:change-proba}: if $(S_j-S_{j-1},
  \eta_j)$, for $j\ge1$ (with $S_0:=0$), are i.i.d. random
variables (under $\mathbf{Q}$) and distributed as $(V(w^{(1)}),
\sum_{u\in \mathscr{I}_1^{(1)}} [\ee^{-s V(u)}
+ \ee^{-{s/2}V(u)}])$, then
\begin{eqnarray*}
&& \mathbf{E}_\mathbf{Q} \bigl\{ \bigl( n^{1/2}\wedge V\bigl(w_n^{(n)}\bigr)^+\bigr)
\underline{Y  }_n^{  -s}  \bigr\}
\\
&&\qquad \ge c_{70}   \mathbf{E}_\mathbf{Q}
\biggl \{ {n^{1/2}\wedge S_n^+ \over
\sum_{j=1}^n \ee^{-sS_{j-1}}
[(S_{j-1}^+)^s + 1 ] \eta_j + \ee^{-sS_n} (S_n^+)^s}  \biggr\}
\\
&&\qquad \ge c_{70}   \mathbf{E}_\mathbf{Q}
\biggl \{ {(n^{1/2}\wedge S_n)
{\bf1}_{\{ \min_{1\le j\le n} S_j>0\} } \over
\sum_{j=1}^n \ee^{-sS_{j-1}} (S_{j-1}^s + 1)
\eta_j + \ee^{-sS_n} S_n^s}  \biggr\} .
\end{eqnarray*}
Note that if $S_j>0$, then $\ee^{-sS_j} [S_j^s + 1 ]
\le c_{71}   \ee^{-tS_j}$ with $t:= {s\over2}$. Therefore, by writing
\[
\mathbf{Q}^{(n)}  \{   \cdot   \} := \mathbf{Q} \biggl\{
  \cdot \big|   \min_{1\le j\le n} S_j>0  \biggr\} ,
\]
and $\mathbf{E}_\mathbf{Q}^{(n)}$ the expectation with
respect to $\mathbf{Q}
^{(n)}$, and $\widehat{\eta}_j := \eta_j +1$ for brevity, we get that
\begin{eqnarray*}
\mathbf{E}_\mathbf{Q} \bigl\{ \bigl( n^{1/2}\wedge V\bigl(w_n^{(n)}\bigr)^+\bigr)
\underline{Y  }_n^{  -s}  \bigr\}
&\ge& c_{72} \mathbf{Q}
\biggl\{ \min_{1\le j\le n} S_j>0  \biggr\}
\mathbf{E}_\mathbf{Q}^{(n)}
\biggl\{ {n^{1/2}\wedge S_n \over
\sum_{j=1}^{n+1} \ee^{-tS_{j-1}} \widehat{\eta}_j}  \biggr\}
\\
&\ge&c_{72} \mathbf{Q} \biggl\{ \min_{1\le j\le n} S_j>0  \biggr\}
\mathbf{E}_\mathbf{Q}^{(n)} \biggl \{ {\varepsilon n^{1/2}
{\bf1}_{ \{ S_n > \varepsilon n^{1/2}\} } \over
\sum_{j=1}^{n+1} \ee^{-tS_{j-1}}
\widehat{\eta}_j}  \biggr\} .
\end{eqnarray*}
Since $\mathbf{Q}\{ \min_{1\le j\le n} S_j>0 \} \ge c_{73}
n^{-1/2}$ [see (\ref{Bingham})], this leads to
\begin{eqnarray*}
&& \mathbf{E}_\mathbf{Q} \bigl\{ \bigl( n^{1/2}\wedge
V\bigl(w_n^{(n)}\bigr)^+\bigr) \underline{Y  }_n^{  -s}  \bigr\}
\\
&&\qquad \ge c_{74}   \varepsilon  \mathbf{E}_\mathbf{Q}^{(n)}
\biggl \{ {{\bf1}_{ \{ S_n > \varepsilon n^{1/2}\} } \over
\sum_{j=1}^{n+1} \ee^{-tS_{j-1}}
\widehat{\eta}_j}  \biggr\}
\\
&&\qquad \ge c_{74}   \varepsilon \biggl[ \mathbf{E}_\mathbf{Q}^{(n)}
 \biggl\{ {1\over \sum_{j=1}^{n+1} \ee^{-tS_{j-1}}
\widehat{\eta}_j}  \biggr\}
- \mathbf{Q}^{(n)} \{ S_n \le\varepsilon n^{1/2}\}  \biggr].
\end{eqnarray*}

Let $\rho(s)>0$ be as in Corollary \ref{c:V}. We have
\mbox{$\mathbf{E}_\mathbf{Q}\{(\sum_{u\in\mathscr{I}_1^{(1)}}
\ee^{-s V(u)})^{\rho(s)} \}<\infty$} by (\ref{existence-moment-Wn}).
Since $\rho(s)\le\rho({s\over 2})$, we also have
$\mathbf{E}_\mathbf{Q}\{ (\sum_{u\in\mathscr{I}_1^{(1)}} \ee
^{-{s/2} V(u)})^{\rho(s)} \}<\infty$. Therefore,
$\mathbf{E}_\mathbf{Q}\{\widehat{\eta}_1^{\rho(s)}\} <\infty$.
We are thus entitled to
apply Lemma \ref{l:kozlov} (stated and proved below) to see that
$\mathbf{E}_\mathbf{Q}^{(n)} \{ {1\over1+ \sum_{j=1}^{n+1}
\ee^{-tS_{j-1}} \widehat{\eta}_j} \} \ge c_{75}$ for
some $c_{75}\in(0,   \infty)$ and all
$n\ge n_0$. Since ${1\over\sum_{j=1}^{n+1} \ee^{-tS_{j-1}}
\widehat{\eta}_j} \ge{1\over1+\sum_{j=1}^{n+1} \ee^{-tS_{j-1}}
\widehat{\eta}_j}$, this yields
\[
\mathbf{E}_\mathbf{Q} \bigl\{ \bigl( n^{1/2}\wedge V\bigl(w_n^{(n)}\bigr)^+\bigr)
\underline{Y  }_n^{ -s}  \bigr\} \ge c_{74}   \varepsilon
\bigl[ c_{75} - \mathbf{Q}^{(n)} \{S_n \le\varepsilon n^{1/2}\}  \bigr] ,
\qquad n\ge n_0.
\]
On the other hand, $S_n/n^{1/2}$ under $\mathbf{Q}^{(n)}$ converges
weakly to the terminal value of a Brownian meander (see
Bolthausen \cite{bolthausen}); in particular,\break
$\lim_{\varepsilon\to 0} \lim_{n\to\infty} \mathbf{Q}^{(n)}
\{ S_n \le\varepsilon n^{1/2}\} =0$.
We can thus choose (and fix) a small $\varepsilon>0$ such that
$\mathbf{Q}^{(n)} \{ S_n \le\varepsilon n^{1/2}\} \le{c_{75}\over2}$ for all
$n\ge n_1$. Therefore, for $n\ge n_0+n_1$,
\[
\mathbf{E}_\mathbf{Q} \bigl\{ \bigl( n^{1/2}\wedge V\bigl(w_n^{(n)}\bigr)^+\bigr)
\underline{Y  }_n^{-s}  \bigr\} \ge c_{74}
\varepsilon   \biggl[ c_{75} - {c_{75}\over2}  \biggr] .
\]
As a consequence, we have proved that, for $0<s<1$,
\[
\liminf_{n\to\infty} \mathbf{E}_\mathbf{Q}
\bigl\{ \bigl( n^{1/2}\wedge V\bigl(w_n^{(n)}\bigr)^+\bigr)
\underline{Y  }_n^{  -s}  \bigr\}>0,
\]
which, in view of (\ref{tartee}), yields the first
inequality in (\ref{tension}), and thus completes the proof of the
lower bound in Theorem \ref{t:tension}.
\end{pf*}

We complete the proof of Theorem \ref{t:tension} by proving the
following lemma, which is a very simple variant of a result of
Kozlov \cite{kozlov}.
\begin{lemma}\label{l:kozlov}
Let $\{ (X_k,  \eta_k),   k\ge1\}$ be a
sequence of i.i.d.\ random vectors defined
on $(\Omega, \mathscr{F}, {\mathbb P})$ with
${\mathbb P}\{ \eta_1 \ge0\} =1$, such that
${\mathbb E}\{ \eta_1^\theta\}<\infty$ for some
$\theta>0$. We assume ${\mathbb E}(X_1)=0$ and
$0<{\mathbb E}(X_1^2)<\infty$.
Let $S_0:=0$ and $S_n := X_1+\cdots+ X_n$, for $n\ge1$. Then
%
\begin{equation}\label{kozlov}
\lim_{n\to\infty}   {\mathbb E} \biggl\{ {1 \over
1+ \sum_{k=1}^{n+1} \eta_k\ee^{-S_{k-1}}}
   \bigg|   \min_{1\le k\le n} S_k >0
 \biggr\} = c_{76} \in(0,   \infty).
\end{equation}
\end{lemma}
\begin{pf}
The lemma is an analogue of the identity (26) of
Kozlov \cite{kozlov}, except that the distribution of our $\eta_1$ is
slightly different from that of Kozlov's, which explains the moment
condition ${\mathbb E}\{ \eta_1^\theta\} < \infty$: this condition
will be seen to guarantee
%
\begin{eqnarray}\label{kozlov-troncature}
&& \lim_{j\to\infty} \limsup_{n\to\infty}
{\mathbb E} \biggl\{ {1\over 1+ \sum_{k=1}^j \eta_k\ee^{-S_{k-1}}}
\nonumber\\[-8pt]
\\[-8pt]
&&\hspace*{70pt}{}
- {1\over 1+ \sum_{k=1}^{n+1} \eta_k\ee^{-S_{k-1}}}
\bigg   |   \min_{1\le k\le n} S_k >0 \biggr\} =0 .
\nonumber
\end{eqnarray}
The identity (\ref{kozlov-troncature}), which plays the role
of Kozlov's Lemma 1 in \cite{kozlov}, is the key ingredient in the
proof of (\ref{kozlov}). Since the rest of the proof goes along the
lines of~\cite{kozlov} with obvious modifications, we only prove
(\ref{kozlov-troncature}) here.

Without loss of generality, we assume $\theta\le2$ (otherwise, we can
replace $\theta$ by 2). We observe that, for $n> j$, the integrand in
(\ref{kozlov-troncature}) is nonnegative, and is
\[
\le{\sum_{k=j+1}^{n+1} \eta_k\ee^{-S_{k-1}}\over
1+ \sum_{k=1}^{n+1} \eta_k\ee^{-S_{k-1}}}
\le \biggl( {\sum_{k=j+1}^{n+1} \eta_k\ee^{-S_{k-1}} \over1+
\sum_{k=1}^{n+1} \eta_k\ee^{-S_{k-1}}}  \biggr)^{\theta/2}
\le \biggl( \sum_{k=j+1}^{n+1} \eta_k\ee^{-S_{k-1}} \biggr)^{\theta/2} ,
\]
which is bounded by $\sum_{k=j+1}^{n+1} \eta_k^{\theta/2}
\ee^{-{\theta/2} S_{k-1}}$. Since
${\mathbb P}\{ \min_{1\le k\le n} S_k >0 \} \sim c_4/n^{1/2}$
[see (\ref{Bingham})], we only need to check that
%
\begin{equation}\label{kozlov-troncature1}
\lim_{j\to\infty} \limsup_{n\to\infty}
n^{1/2} \sum_{k=j+1}^{n+1} {\mathbb E}
\bigl\{\eta_k^{\theta/2}\ee^{-{\theta/2} S_{k-1}}
  {\bf1}_{ \{ \min_{1\le i\le n} S_i >0\} } \bigr\} =0 .
\end{equation}

Let $\mathrm{LHS}_{(\ref{kozlov-troncature1})}$ denote the
$n^{1/2}\sum_{k=j+1}^{n+1} {\mathbb E}\{ \cdot\}$ expression on the
left-hand side. Let $\widehat{S}_i = \widehat{S}_i (k) :=
S_{i+k}-S_k$, $i\ge0$. It is clear that $(\widehat{S}_i,   i\ge0)$
is independent of $(\eta_k,   X_1, \ldots , X_k)$, and is
distributed as $(S_i,   i\ge0)$. Write $\underline{S}_{k-1}
:= \min_{1\le j \le k-1} S_j$ and $\underline{\widehat{S}}_{n-k}
:= \min_{1\le i\le n-k} \widehat{S}_i$. Then
\[
\mathrm{LHS}_{(\ref{kozlov-troncature1})} \le n^{1/2}
\sum_{k=j+1}^{n+1} {\mathbb E}
\bigl\{ \eta_k^{\theta/2} \ee^{-{\theta/2}
S_{k-1}}   {\bf1}_{ \{ \underline{S}_{k-1}>0,   \underline
{\widehat{S}}_{n-k} > -S_{k-1} -X_k\} }  \bigr\} .
\]
To estimate ${\mathbb E}\{\cdot\}$ on the right-hand side,
we first condition upon $(\eta_k,   S_{k-1},   \underline{S}_{k-1},\break
  X_k)$, which leaves us to estimate the tail probability of
$\underline{\widehat{S}}_{n-k}$. At this stage, it is convenient to
recall (see (13) of Kozlov \cite{kozlov}) that
${\mathbb P}\{\underline{\widehat{S}}_{n-k} > -y\} \le c_{54}
{1+y^+\over (n-k+1)^{1/2}}$ for some $c_{54}>0$ and all
$y\in\mathbb{R}$. Accordingly,
\begin{eqnarray*}
\mathrm{LHS}_{(\ref{kozlov-troncature1})}
&\le& c_{54}  n^{1/2} \sum_{k=j+1}^{n+1}
{\mathbb E} \biggl\{ \eta_k^{\theta/2}
\ee^{-{\theta/2} S_{k-1}}
{\bf1}_{ \{ \underline{S}_{k-1}>0\} }
{1+ (S_{k-1} + X_k)^+ \over(n-k+1)^{1/2}} \biggr\}
\\
&\le& c_{54}  n^{1/2} \sum_{k=j+1}^{n+1}
{\mathbb E}\biggl \{ \eta_k^{\theta/2}
\ee^{-{\theta/2} S_{k-1}}
{\bf1}_{ \{ \underline{S}_{k-1}>0\} }
{1+ S_{k-1} + X_k^+ \over(n-k+1)^{1/2}} \biggr\} .
\end{eqnarray*}

On the right-hand side, $(\eta_k,   X_k)$ is independent of
$(\underline{S}_{k-1},   S_{k-1})$. We condition upon
$(\underline{S}_{k-1},   S_{k-1})$: for any $z\ge1$, an application of the
Cauchy--Schwarz inequality gives
\[
{\mathbb E}\{ \eta_k^{\theta/2} (z + X_k^+)\} \le
[{\mathbb E}( \eta_k^\theta)]^{1/2}  [{\mathbb E}\{(z + X_k^+)^2\}]^{1/2} .
\]
Of course, ${\mathbb E}( \eta_k^\theta) = {\mathbb E}( \eta_1^\theta)<\infty$
by assumption, and ${\mathbb E}\{(z + X_k^+)^2\}\le2
{\mathbb E}(z^2 + X_k^2)= 2[z^2+ {\mathbb E}(X_1^2)]$. Thus,
${\mathbb E}\{ \eta_k^{\theta/2} (z + X_k^+)\} \le c_{77}   z$ for
$z\ge1$. Consequently, with $c_{78} := c_{54} c_{77}$,
\begin{eqnarray*}
\mathrm{LHS}_{(\ref{kozlov-troncature1})}
&\le& c_{78}   n^{1/2} \sum_{k=j+1}^{n+1}
{\mathbb E} \biggl\{ \ee^{-{\theta/2} S_{k-1}}
{\bf1}_{ \{ \underline{S}_{k-1}>0\} }
{1+ S_{k-1} \over(n-k+1)^{1/2}}  \biggr\}
\\
&\le&c_{79}   n^{1/2} \sum_{k=j+1}^{n+1}
{\mathbb E} \biggl\{ \ee^{-{\theta/3} S_{k-1}}
{\bf1}_{ \{ \underline{S}_{k-1}>0\} }
{1\over(n-k+1)^{1/2}}  \biggr\} ,
\end{eqnarray*}
the last inequality following from the fact that $\sup_{x>0}
(1+x) \ee^{-{\theta/6} x} <\infty$.

We use once again the estimate (\ref{Bingham}), which implies
${1\over(n-k+1)^{1/2}} \le c_{80}  {\mathbb P}\{ S_i > S_{k-1},   \forall
k\le i\le n\}$. Since $(S_i-S_{k-1},   k\le i\le n)$ is independent of
$(S_{k-1},   \underline{S}_{k-1})$, this implies, with
$c_{81} :=c_{79} c_{80}$,
\begin{eqnarray*}
\mathrm{LHS}_{(\ref{kozlov-troncature1})}
&\le& c_{81}   n^{1/2} \sum_{k=j+1}^{n+1}
{\mathbb E} \bigl\{ \ee^{-{\theta/3} S_{k-1}}
{\bf1}_{ \{ \underline{S}_{k-1}>0,
S_i > S_{k-1},   \forall k\le i\le n\} } \bigr\}
\\
&\le&c_{81}   n^{1/2} \sum_{k=j+1}^{n+1}
{\mathbb E} \bigl\{ \ee^{-{\theta/3} S_{k-1}}
{\bf1}_{ \{ \underline{S}_n>0\} } \bigr\} ,
\end{eqnarray*}
where $\underline{S}_n := \min_{1\le i\le n} S_i$. It
remains to check that
%
\begin{equation}\label{kozlov-2}
\lim_{j\to\infty} \limsup_{n\to\infty}
n^{1/2} \sum_{k=j+1}^{n+1}
{\mathbb E} \bigl\{ \ee^{-{\theta/3} S_{k-1}}
{\bf1}_{ \{ \underline{S}_n>0\} }  \bigr\}=0.
\end{equation}
This would immediately follow from Lemma 1 of Kozlov
\cite{kozlov}, but we have been kindly informed by Gerold Alsmeyer (to whom
we are grateful) of a flaw in its proof, on page 800, line 3 of
\cite{kozlov}, so we need to proceed differently. Since
${\mathbb E}\{\ee^{-{\theta/3} S_{k-1}}   {\bf1}_{ \{ \underline{S}_n>0\} } \}
\le n^{-(3/2)+o(1)} (n-k+2)^{-1/2}$ (for $n\to\infty$) uniformly in
$k\in[{n\over2},   n+1]$, we have $n^{1/2} \sum_{k=\lfloor
n/2\rfloor}^{n+1}{\mathbb E}\{ \ee^{-{\theta/3} S_{k-1}}
{\bf1}_{ \{ \underline{S}_n>0\} } \} \to0$, $n\to\infty$. On the
other hand, (36) of Kozlov \cite{kozlov} (applied to $\delta={1\over
2}$ and $\eta_i =1$ there) implies that $\lim_j \limsup_n n^{1/2}
\sum_{k=j+1}^{\lfloor n/2\rfloor} {\mathbb E}\{\ee^{-{\theta/3}
S_{k-1}}   {\bf1}_{ \{ \underline{S}_n>0\} } \} =0$. Therefore,
(\ref{kozlov-2}) holds: Lemma \ref{l:kozlov} is proved.
\end{pf}

\section[Proof of Theorem 1.3 and (1.14)--(1.15)
of Theorem 1.4]{Proof of Theorem \protect\ref{t:Mn} and
(\protect\ref{Wnbeta-liminf-as})--(\protect\ref{Wnbeta-proba})
of Theorem \protect\ref{t:derrida-spohn}}\label{s:Mn}

In this section we prove Theorem \ref{t:Mn}, as well as parts
(\ref{Wnbeta-liminf-as})--(\ref{Wnbeta-proba}) of Theorem
\ref{t:derrida-spohn}. We assume~(\ref{hyp2}), (\ref{hyp3}) and
(\ref{hyp}) throughout the section.
\begin{pf*}{Proof of Theorem \ref{t:Mn} and
(\ref{Wnbeta-liminf-as}) and (\ref{Wnbeta-proba}) of Theorem \ref{t:derrida-spohn}}
\textit{Upper bounds.} Let $\varepsilon>0$. By Theorem
\ref{t:derrida-spohn-moment} and Chebyshev's inequality,
$\mathbf{P}\{W_{n,\beta} >n^{-(3\beta/2) + \varepsilon} \} \to0$.
Therefore, $W_{n,\beta}\le n^{-(3\beta/2)+o(1)}$ in probability,
yielding the upper bound in (\ref{Wnbeta-proba}).

The upper bound in (\ref{Wnbeta-liminf-as}) follows trivially from the
upper bound in (\ref{Wnbeta-proba}).

It remains to prove the upper bound in Theorem \ref{t:Mn}. Fix
$\gamma\in(0,   1)$. Since $W_n^\gamma$ is a nonnegative supermartingale,
the maximal inequality tells that, for any $n\le m$ and any $\lambda>0$,
\[
\mathbf{P} \biggl\{ \max_{n\le j\le m} W_j^\gamma\ge\lambda \biggr\}
\le{\mathbf{E}
(W_n^\gamma)\over\lambda} \le{c_{82}\over\lambda n^{\gamma/2}},
\]
the last inequality being a consequence of Theorem \ref
{t:tension}. Let $\varepsilon>0$ and let $n_k := \lfloor
k^{2/\varepsilon}\rfloor$. Then $\sum_k \mathbf{P}\{ \max_{n_k\le
j\le
n_{k+1}} W_j^\gamma\ge n_k^{-(\gamma/2)+\varepsilon}\} <\infty$. By
the Borel--Cantelli lemma, almost surely for all large $k$, $\max
_{n_k\le j\le n_{k+1}} W_j < n_k^{-(1/2)+(\varepsilon/\gamma)}$.
Since ${\varepsilon\over\gamma}$~can be arbitrarily small, this
yields the desired upper bound: $W_n \le n^{-(1/2)+o(1)}$ a.s.
\end{pf*}
\begin{pf*}{Proof of Theorem \ref{t:Mn} and
(\ref{Wnbeta-liminf-as}) and (\ref{Wnbeta-proba}) of Theorem
\ref{t:derrida-spohn}}
\textit{Lower bounds}. To prove the lower bound in
(\ref{Wnbeta-liminf-as}) and (\ref{Wnbeta-proba}), we use the Paley--Zygmund
inequality and Theorem \ref{t:derrida-spohn-moment} to see that
%
\begin{equation}\label{P(Wnbeta>)>}
\mathbf{P}\bigl\{ W_{n,\beta} > n^{-(3\beta/2)+o(1)} \bigr\} \ge
n^{o(1)}, \qquad n\to\infty.
\end{equation}
This is the analogue of (\ref{McDiarmid1}) for $W_n$. {F}rom
here, the argument follows the lines in the proof of the upper bound in
(\ref{V-liminf-as}) of Theorem \ref{t:leftmost} (Section
\ref{s:proof(1.7)}), and goes as follows: let $\varepsilon>0$ and let
$\tau_n := \inf\{ k\ge1\dvtx  \#\{ u\dvtx |u|=k\} \ge n^{2\varepsilon}\}$.
Then
\begin{eqnarray*}
&& \mathbf{P} \biggl\{ \tau_n <\infty,
\min_{k\in[{n/2},   n]} W_{k+\tau_n, \beta}
\le n^{-(3\beta/2)-\varepsilon}
\exp\biggl[-\beta\max_{|x|=\tau_n} V(x)\biggr]  \biggr\}
\\
&&\qquad \le \sum_{k\in[{n/2},   n]}
\mathbf{P} \biggl\{ \tau_n <\infty,   W_{k+\tau_n,\beta}
\le n^{-(3\beta/2)-\varepsilon}
\exp\biggl[-\beta\max_{|x|=\tau_n} V(x)\biggr]  \biggr\}
\\
&&\qquad \le \sum_{k\in[{n/2},   n]}
\bigl ( \mathbf{P} \bigl\{ W_{k,\beta}
\le n^{-(3\beta/2)-\varepsilon}  \bigr\}
 \bigr)^{\lfloor n^{2\varepsilon}\rfloor} ,
\end{eqnarray*}
which, according to (\ref{P(Wnbeta>)>}), is bounded by
$n\exp( - n^{-\varepsilon} \lfloor n^{2\varepsilon}\rfloor)$ (for all
sufficiently large $n$), thus summable in $n$. By the Borel--Cantelli
lemma, almost surely for all sufficiently large $n$, we have either
$\tau_n=\infty$, or $\min_{k\in[{n/2},   n]} W_{k+\tau
_n,\beta} > n^{-(3\beta/2)-\varepsilon} \exp[-\beta\max_{|x|=\tau
_n} V(x)]$. Conditionally on the system's ultimate survival, we have
${1\over n} \max_{|x|=n} V(x) \to c_{21}$ a.s., $\tau_n \sim
{2\varepsilon\log n\over\log m}$ a.s., $n\to\infty$, and
$W_{n,\beta}\ge\min_{k\in[{n/2},   n]} W_{k+\tau_n,\beta}$
for all sufficiently large $n$. This readily yields lower bounds in
(\ref{Wnbeta-liminf-as}) and (\ref{Wnbeta-proba}): conditionally on
the system's survival, $W_{n,\beta} \ge n^{-(3\beta/2)+o(1)}$ almost
surely (and a fortiori, in probability).

The lower bound in Theorem \ref{t:Mn} is along exactly the same lines,
but using Theorem \ref{t:tension} instead of Theorem
\ref{t:derrida-spohn-moment}.
\end{pf*}

\section[Proof of Theorem 1.2]{Proof of Theorem \protect\ref{t:leftmost}}\label{s:leftmost}

Assume (\ref{hyp2}), (\ref{hyp3}) and
(\ref{hyp}). Let $\beta>1$. We trivially have $W_{n,\beta} \le W_n
\exp\{ -(\beta-1)\inf_{|u|=n} V(u)\}$ and $W_{n,\beta} \ge\exp\{
- \beta\times\break \inf_{|u|=n} V(u)\}$. Therefore, ${1\over\beta} \log{1\over
W_{n,\beta}} \le\inf_{|u|=n} V(u) \le{1\over\beta-1} \log{W_n
\over W_{n,\beta}}$ on $\mathscr{S}_n$. Since $\beta$ can be as
large as
possible, by means of Theorem \ref{t:Mn} and of parts
(\ref{Wnbeta-liminf-as}) and (\ref{Wnbeta-proba}) of Theorem
\ref{t:derrida-spohn}, we immediately get (\ref{V-limsup-as})
and (\ref{V-proba}).

Since $W_n \ge\exp\{ - \inf_{|u|=n} V(u)\}$, the lower bound in
(\ref{V-liminf-as}) follows immediately from Theorem \ref{t:Mn},
whereas the upper bound in (\ref{V-liminf-as}) was already proved in
Section~\ref{s:proof(1.7)}.

\section[Proof of part (1.13) of
Theorem 1.4]{Proof of part (\protect\ref{Wnbeta-limsup-as}) of
Theorem \protect\ref{t:derrida-spohn}}\label{s:derrida-spohn}

The upper bound follows from Theorem \ref{t:Mn} and the elementary
inequality $W_{n,\beta} \le W_n^\beta$, the lower bound from
(\ref{V-liminf-as}) and the relation $W_{n,\beta}
\ge\exp\{ -\beta\inf_{|u|=n} V(u)\}$.

\section[Proof of Theorem 1.1]{Proof of Theorem \protect\ref{t:main}}\label{s:proof-t:main}

The proof of Theorem \ref{t:main} relies on Theorem \ref{t:tension}
and a preliminary result, stated below as Proposition
\ref{p:regularity}. Theorem \ref{t:tension} ensures the tightness of
$(n^{1/2} W_n,   n\ge1)$, whereas Proposition \ref{p:regularity}
implies that ${W_{n+1}\over W_n}$ converges to 1 in probability
(conditionally on the system's survival).
\begin{proposition}\label{p:regularity}
Assume (\ref{hyp2}), (\ref{hyp3}) and (\ref{hyp}).
For any $\gamma>0$, there exists
$\gamma_1>0$ such that, for all sufficiently large $n$,
%
\begin{equation}\label{regularity}
\mathbf{P} \biggl\{  \biggl|{W_{n+1}\over W_n} -1  \biggr|
\ge n^{-\gamma}    \big|   \mathscr{S} \biggr\}\le n^{-\gamma_1}.
\end{equation}
\end{proposition}
\begin{pf}
Let $1<\beta\le\min\{ 2,   1+\rho(1)\}$,
where $\rho(1)$ is the constant in Corollary~\ref{c:V}.

We use a probability estimate of Petrov \cite{petrov}, page~82: for
centered random variables $\xi_1, \ldots, \xi_\ell$ with
$\mathbf{E}(|\xi_i|^\beta)<\infty$ (for $1\le i\le\ell$), we have
$\mathbf{E}\{ |\sum_{i=1}^\ell\xi_i|^\beta\} \le2
\sum_{i=1}^\ell\mathbf{E}\{|\xi_i|^\beta\}$.

By definition, on the set $\mathscr{S}_n$, we have
\[
{W_{n+1}\over W_n} - 1 = \sum_{|u|=n} {\ee^{-V(u)} \over W_n}
\Biggl (\sum_{x\in\mathbb{T}^{\mathrm{GW}}_u:   |x|_u=1}
\ee^{-V_u(x)} - 1 \Biggr) ,
\]
where $\mathbb{T}^{\mathrm{GW}}$ and $|x|_u$ are as in
(\ref{TGW}) and (\ref{|x|_u}), respectively. Conditioning on
${\mathscr F}_n$, and applying Proposition
\ref{p:change-proba} and Petrov's
probability inequality recalled above, we see that, on $\mathscr{S}_n$,
%
\begin{eqnarray}\label{M{n+1}/Mn}
\mathbf{E} \biggl\{  \biggl| {W_{n+1}\over W_n} - 1
 \biggr|^\beta   \big|   {\mathscr F}_n
 \biggr\} &\le& 2 \sum_{|u|=n} {\ee^{-\beta V(u)}
\over W_n^\beta} \mathbf{E}
\Biggl \{  \Biggl|\sum_{|y|=1} \ee^{-V(y)} - 1
\Biggr |^\beta \Biggr\}
\nonumber\\[-8pt]
\\[-8pt]
&=& c_{83}{W_{n,\beta} \over W_n^\beta} ,
\nonumber
\end{eqnarray}
where $c_{83} := 2 \mathbf{E}\{ | \sum_{|v|=1} \ee^{-V(v)} - 1
|^\beta\}<\infty$ [see (\ref{existence-moment-Wn})],
and $W_{n,\beta}$ is as in (\ref{Wnbeta}).

Let $\varepsilon>0$ and $b>0$. Let $s\in({\beta-1\over\beta},
1)$. Define $D_n := \{ W_n \ge n^{-(1/2)-\varepsilon} \} \cap
\{W_{n,\beta} \le n^{-(3\beta/2)+b}\}$. By Proposition
\ref{p:tail-Mn}, $\mathbf{P}\{ W_n < n^{-(1/2)-\varepsilon},
\mathscr{S}\} \le n^{-\vartheta}$ for some $\vartheta>0$ and all large $n$,
whereas, by Theorem \ref{t:derrida-spohn-moment},
$\mathbf{P}\{ W_{n,\beta} >n^{-(3\beta/2)+b} \} \le n^{3\beta(1-s)/2 -(1-s)b}
\mathbf{E}\{W_{n,\beta}^{1-s}\} = n^{-(1-s)b +o(1)}$. Therefore,
\[
\mathbf{P} \{ \mathscr{S}\setminus D_n  \}
\le n^{-\vartheta} + n^{-(1-s)b +o(1)}, \qquad n\to\infty.
\]
On the other hand, since $\mathscr{S}\subset\mathscr{S}_n$, it follows from
(\ref{M{n+1}/Mn}) and Chebyshev's inequality that, for $n\to\infty$,
\begin{eqnarray*}
\mathbf{P} \biggl\{  \biggl| {W_{n+1}\over W_n} - 1  \biggr|
\ge n^{-\gamma},   D_n,   \mathscr{S} \biggr\}
&\le& n^{\gamma\beta}   \mathbf{E}
\biggl\{ c_{83}{W_{n,\beta} \over W_n^\beta}   {\bf1}_{ D_n
\cap\mathscr{S}_n}  \biggr\}
\\
&\le& c_{83}   n^{\gamma\beta-(3\beta/2) + b +
[(1/2)+ \varepsilon]\beta}.
\end{eqnarray*}
As a consequence, when $n\to\infty$,
\[
\mathbf{P} \biggl\{  \biggl| {W_{n+1}\over W_n} - 1  \biggr|
\ge n^{-\gamma},   \mathscr{S} \biggr\}
\le n^{-\vartheta}+ n^{-(1-s)b +o(1)} + c_{83}
n^{\gamma\beta-\beta+ b + \varepsilon\beta}.
\]
We choose $\varepsilon$ and $b$ sufficiently small such that
$\gamma\beta-\beta+ b + \varepsilon\beta<0$.
Proposition~\ref{p:regularity} is proved.
\end{pf}

We now have all of the ingredients needed for the proof of
Theorem \ref{t:main}.
\begin{pf*}{Proof of Theorem \ref{t:main}}
Once Proposition \ref{p:regularity} is established, the
proof of Theorem \ref{t:main} follows
the lines of Biggins and Kyprianou \cite{biggins-kyprianou97}.

Assume (\ref{hyp2}), (\ref{hyp3}) and (\ref{hyp}). Let
$\lambda_n>0$ satisfy $\mathbf{E}\{ (\lambda_n W_n)^{1/2}\} =1$.
That is,
\[
\lambda_n :=  \{ \mathbf{E}(W_n^{1/2})  \} ^{-2}.
\]
By Theorem \ref{t:tension}, we have
$0<\liminf_{n\to\infty} {\lambda_n\over n^{1/2}}
\le\limsup_{n\to\infty}{\lambda_n\over
n^{1/2}} <\infty$, and $(\lambda_n W_n)$ is tight. Let $\widehat
{\mathscr{W}}$ be any possible (weak) limit of $(\lambda_n W_n)$
along\vspace*{1pt} a subsequence. By Theorem \ref{t:tension} and dominated
convergence, $\mathbf{E}(\widehat{\mathscr{W}}^{1/2})=1$. We now
prove the uniqueness of $\widehat{\mathscr{W}}$.

By definition,
\[
W_{n+1} = \sum_{|v|=1} \ee^{-V(v)}
\sum_{x\in\mathbb{T}^{\mathrm{GW}}_v,
 |x|_v= n} \ee^{-V_v(x)} .
\]
By assumption, $\lambda_n W_n \to\widehat{\mathscr{W}}$
in distribution when $n$ goes to infinity along a certain subsequence.
Thus, $\lambda_n W_{n+1}$ converges weakly (when $n$ goes along the
same subsequence) to $\sum_{|v|=1} \ee^{-V(v)}
\widehat{\mathscr{W}}_v$, where, conditionally on $(v, V(v),   |v|=1)$,
$\widehat{\mathscr{W}}_v$ are independent copies of $\widehat{\mathscr{W}}$.

On the other hand, by Proposition \ref{p:regularity},
$\lambda_nW_{n+1}$ also converges weakly (along the same subsequence) to
$\widehat{\mathscr{W}}$. Therefore,
\[
\widehat{\mathscr{W}}   \stackrel{\mathit{law}}{=}   \sum_{|v|=1}
\ee^{-V(v)} \widehat{\mathscr{W}}_v .
\]
This is the same equation for $\xi^*$ in (\ref{xi*}).
Recall that (\ref{xi*}) has a unique solution up to a scale change
(Liu \cite{liu00}), and since $\mathbf{E}(\widehat{\mathscr{W}}^{1/2})=1$,
we have $\widehat{\mathscr{W}}   \stackrel{\mathit{law}}{=}   c_{84}
\xi^*$, with $c_{84} := [ \mathbf{E}\{ (\xi^*)^{1/2}\}]^{-2}$. The
uniqueness (in law) of $\widehat{\mathscr{W}}$ shows that $\lambda_n
W_n$ converges weakly to $\widehat{\mathscr{W}}$ when $n\to\infty$.

By (\ref{survival-survival}), $\mathbf{P}\{ W_n>0\}
= \mathbf{P}\{\mathscr{S}_n\}\to\mathbf{P}\{ \mathscr{S}\} = \mathbf{P}\{ \xi^*>0\}$.
Let $\mathscr{W}>0$ be a random variable such that
%
\begin{equation}\label{wcv}
\mathbf{E}(\ee^{-a \mathscr{W} }) =
\mathbf{E}(\ee^{-a \widehat{\mathscr{W}} }   |
\widehat{\mathscr{W}}>0) , \qquad \forall a\ge0.
\end{equation}
It follows that, conditionally on the system's survival,
$\lambda_n W_n$ converges in distribution to $\mathscr{W}$.
\end{pf*}

\section*{Acknowledgments}

We are grateful to Bruno Jaffuel, Gerold Alsmeyer, John
Biggins and Alain Rouault for pointing out errors in the first versions
of the manuscript, and to John Biggins for fixing Lemma \ref{l:phi*}.
We also wish to thank two referees for their careful reading and for
very helpful suggestions, without which we would not have had enough
courage to remove the uniform ellipticity condition in the original manuscript.


\printaddresses

\end{document}